\documentclass[11pt,a4paper]{article}
\usepackage{amsmath,amssymb,amsfonts}
\usepackage{a4wide}
\usepackage{epsfig,epic,eepic,graphicx}
\newtheorem{theorem}{Theorem}
\newtheorem{lemma}{Lemma}
\newtheorem{proposition}{Proposition}

\newtheorem{remark}{Remark}
\makeatletter

\@addtoreset{equation}{section}
\makeatother

\newcommand{\T}{{\mathbb T}}
\newcommand{\N}{{\mathbb N}}

\newcommand{\R}{{\mathbb R}}

\newcommand{\pa}{{\partial}}
\newcommand{\om}{{\omega}}
\newcommand{\Om}{{\Omega}}
\newcommand{\na}{{\nabla}}

\newcommand{\eps}{{\varepsilon}}

\def\u{\mathbf{u}}
\def\x{\mathbf{x}}

\def\sspace{\smallskip \noindent}
\def\mspace{\medskip \noindent}

\numberwithin{equation}{section}

\title{Well-posedness for  the Prandtl system \\ without analyticity or monotonicity}
\author{\footnote{Institut de Math\'ematiques de Jussieu et Universit\'e Paris 7, 175 rue du Chevaleret, 75013 Paris France} David G\'erard-Varet, 
 \footnote{Courant Institute, NYU, 251 Mercer Street, New-York 10012. Partially supported by NSF  grant  
DMS-1211806 } Nader Masmoudi}

\date{}

\begin{document}
\maketitle

\begin{abstract}
It has been thought for a while that the  Prandtl system is only  well-posed 
under the Oleinik monotonicity  assumption or  under an analyticity assumption. 
We show that the Prandtl system is actually  locally well-posed for data  that belong to the 
 Gevrey class  $7/4$   in the horizontal  variable $x$. Our result improves the classical local well-posedness result for data that are analytic in $x$ (that is Gevrey class $1$). The proof uses 
new estimates, based on non-quadratic energy functionals.

\end{abstract}

\section{Introduction}

Our concern in this  paper is the well-posedness of the Prandtl system. This system, by now classical, was introduced by Prandtl in 1904 to describe an incompressible flow near a wall,  at high Reynolds number. Formally, it is derived from the Navier-Stokes equation with no-slip condition: 
\begin{equation} \label{NS}
\left\{
\begin{aligned}
\pa_t \u  + \u \cdot \na \u + \na p - \eps \Delta \u \:
 = \: 0, & 
\quad  \x \in \Om, \\
\na \cdot \u \:  = \: 0, & \quad \x \in \Om,\\
\u\vert_{\pa \Om}  =  0, &
\end{aligned}
\right. 
\end{equation}
that we consider for simplicity in $\Om := \T \times \R_+$. We recall  that $\u(t,\x) = (u,v)(t,x,y)$ is the velocity field of the fluid, and $p$ its pressure field. The  parameter $0 < \eps \ll 1$ is the inverse of the Reynolds number. In the limit case  $\eps = 0$, one is left formally  with the Euler equation,  for which only the impermeability condition $u \cdot n\vert_{\pa \Om} = 0$ can be prescribed. Mathematically,  this singular change of boundary condition generates strong gradients of the Navier-Stokes solution $\u^\eps$, as $\eps \rightarrow 0$. These gradients correspond to a concentration of the fluid flow in a  thin zone near the wall $\pa \Om$: the so-called  {\em boundary layer}.  The understanding of the boundary layer is a great mathematical challenge, that makes the convergence of Navier-Stokes  solutions to Euler ones a big open problem, even for smooth data.     
 
\medskip
To tackle this problem, Prandtl proposed in 1904 an asymptotic model for the flow, based on two different asymptotic expansions of  $\u^\eps$, resp.  outside and  inside the boundary layer:  
\begin{itemize}
\item outside the boundary layer,  no concentration should occur: one should have 
$$ \u^\eps(t,\x)  \sim \u^0(t,\x) , \quad \mbox{ the solution of the Euler equation.} $$ 
\item inside the boundary layer, $\u^\eps$ should exhibit strong gradients, transversally to the boundary: more precisely, the asymptotics suggested by Prandtl is  
$$ u^\eps(t,x,y) \sim u(t,x,y/\sqrt{\eps}), \quad v^\eps(t,x,y) \sim \sqrt{\eps} v(t,x,y/\sqrt{\eps}) $$
where $u = u(t,x,Y)$ and $v = v(t,x,Y)$ are {\em boundary layer profiles}, depending on a rescaled variable $Y=y/\sqrt{\eps}$, $\: Y > 0$. 
Note that the scale $\sqrt{\eps}$ is coherent with the parabolic part of (\ref{NS}a). 
\end{itemize}
If we plug the expansion above in  (\ref{NS})  and keep the leading order terms, we derive the famous Prandtl system (denoting $Y$ instead of $y$): 
\begin{equation} \label{Prandtl1}
\left\{
\begin{aligned}
\pa_t u + u \pa_x u + v \pa_y u + \pa_x p   - \pa^2_y u = 0, &  \\
\pa_y p = 0, & \\
\pa_x u + \pa_y v = 0, &  \\
 u\vert_{y=0} = v\vert_{y=0} = 0, & \quad \lim_{y \rightarrow +\infty} u = U, \quad  \lim_{y \rightarrow +\infty} p = P,
 \end{aligned}
\right.
\end{equation} 
where $U(t,x) := u^0(t,x,0)$ and $P(t,x) := p^0(t,x,0)$ are the Euler tangential velocity and pressure at the boundary. We refer to \cite{Masmoudi07hand} for  the formal  derivation 
of the  Prandtl  system. 
   The condition at $y =+\infty$  in (\ref{Prandtl1}d) is a matching condition near the boundary between the boundary layer flow and the Euler flow  ({\em matched asymptotics}). Note that,  combining (\ref{Prandtl1}b) with the boundary condition on  $p$, we get $p \equiv P$. Hence, the pressure is not an unknown in the Prandtl model: $v$ is obtained in terms of $u$ by integrating  the divergence-free  condition  (\ref{Prandtl1}c), so that  (\ref{Prandtl1}a) is a scalar evolution equation on $u$, which is {\it a priori} much simpler than the original Navier-Stokes equation. 

\medskip
However, this appealing formal asymptotics raises strong mathematical issues: well-posedness of the limit Prandtl system on one hand, justification of the Prandtl asymptotics of $\u^\eps$ on the other hand. The difficulty comes from numerous underlying fluid instabilities, that can invalidate the Prandtl model: we refer to \cite{Guy:2001} for a basic presentation of these aspects.

\medskip
The aim of the present paper is to investigate this stability problem, from a mathematical viewpoint. We shall focus on the limit  Prandtl system, namely on its well-posedness. For simplicity, we shall restrict to homogeneous data: $U = P = 0$. Extension of our results to the case of constant $U$ would not raise any problem. Extension to some $U = U(t,x)$ would require some  modifications, see  \cite{KV13} for a similar problem. Hence, we consider here the following system:      
\begin{equation} \label{Prandtl}  
\left\{
\begin{aligned}
\pa_t u + u \pa_x u + v \pa_y u  - \pa^2_y u = 0, &  \\
\pa_x u + \pa_y v = 0, &  \\
u\vert_{y=0} = v\vert_{y=0} = 0, \:\lim_{y=+\infty} u = 0 & 
\end{aligned}
\right.
\end{equation} 
with initial condition $u\vert_{t=0} = u_0$. 

\medskip
Before stating our theorem, let us review briefly known results on  the existence theory for \eqref{Prandtl1}. So far,  well-posedness has been established in two settings:
\begin{itemize}
\item The first results go back to Oleinik \cite{OS99}, who obtained some local well-posedness for {\em initial data that are monotonic with respect to $y$}: $U > 0$, $\pa_y u > 0$.  For such data, one can use the  Crocco transform: in short, using $u$ as   an independent  variable instead of $y$ and $w:= \pa_y u$ as an unknown instead of $u$, one is left with a nonlinear parabolic  equation on $w$, for which maximum principles are available: see \cite{OS99} for details. Note that under the extra condition  $\pa_x P  \le 0$ ({\em favorable pressure gradient}), one can go from local to global well-posedness, {\it cf} \cite{XZ04}. From the point of view of physics, this monotonicity assumption is known to be stabilizing: it avoids  the  {\em boundary layer separation}, see \cite{Guy:2001}. 
\item Without monotonicity, well-posedness has been established only locally in time, for {\em initial data that are analytic with respect to $x$}. We refer to \cite{SC98b,LCS03}, and to the recent extension \cite{KV13}. The assumption of analyticity can be understood as follows. By the divergence-free condition, one obtains 
$v = - \int_0^y \pa_x u$. Thus, the term $v \, \pa_y u$ in (\ref{Prandtl}a) (seen as a functional of $u$) is first order in $x$. Moreover, it is not hyperbolic. For instance, let us consider the linearization of  the Prandtl equation around a shear flow $\u = (U_s(y), 0)$:  
\begin{equation} \label{linearization}
\pa_t u + U_s \pa_x u  +  U'_s v - \pa^2_y u = 0, \quad \pa_x u + \pa_y v = 0.
\end{equation}      
If we freeze  the coefficients at some $y_0$ and compute the  dispersion relation, we obtain the  growth rate
$$ \sigma(k_x,k_y) =  U'_s(y_0) \frac{k_x}{k_y}  - k_y^2  $$
that increases linearly with the wavenumber $k_x$. This kind of growth rate would prevent any well-posedness result outside the analytic setting. 
\end{itemize}
However, as discussed in \cite{HH03}, this dispersion relation, formally obtained by freezing the coefficients, is misleading: for instance, the inviscid version of  Prandtl (that is removing the $\pa^2_y u$ term ) is locally well-posed in $C^k$, through the method of characteristics. 

\medskip
In the case of the full Prandtl system \eqref{Prandtl}, the situation is even more complex, and was addressed  recently by the first author and Emmanuel Dormy in  article \cite{GD10} (see also \cite{GVN12}). This article  contains a careful study of the linearized system \eqref{linearization}, {\em in the case of a  non-monotonic base flow $U_s$}: 
$$\exists a, \quad  U'_s(a) = 0, \quad  U''_s(a) \neq 0. $$
In short, it is shown  in  \cite{GD10}\footnote{As mentioned recently by S. Cowley to the first author, the instability mechanism used in \cite{GD10} had already been described at a formal level in \cite{Cowley}.} that the linear system  \eqref{linearization} admits approximate solutions  with growth rate    
\begin{equation} \label{growthrate}
 \sigma(k_x) \: \sim \: \delta \sqrt{k_x}, \quad \delta > 0, \quad k_x \gg 1. 
 \end{equation}
Let us stress that such growing solutions result from an interplay between the lack of monotonicity of $U_s$ and the diffusion term $\pa^2_y u$.  It is therefore coherent with the well-posedness results  obtained in the monotonic case, and for the inviscid Prandtl equation.  

\medskip
Of course, the growth rate indicated in \eqref{growthrate} yields ill-posedness of the Prandtl equation in Sobolev spaces. Still, it leaves room for well-posedness below analytic regularity. {\em Indeed, the aim of the present paper is to  prove that the Prandtl equation is locally well-posed for data that are of Gevrey  class $7/4$ in variable $x$.} Precise statements will be given in the next section.  Let us already point out  that the Gevrey classes $m$, $m > 1$, contain compactly supported functions, as opposed to the Gevrey  class 1 (analytic functions). Hence,  stability in the Gevrey context has more physical insight that stability in the analytic context: see \cite{MoVi:2011} for a similar issue.

\section{Statements of  the result} \label{statements}
Let $m \ge   1$. We recall that  the Gevrey space $G^m(\T)$ is the set of functions $f$   satisfying: $\exists C, \tau > 0 $ such that 
$$ | f^{(j)}(x) | \: \le \: C\, \tau^{-j} (j!)^m, \quad \forall j \in \N, \: x \in \T. $$ 
For a reminder on  Gevrey spaces, we  refer to the 
 paper  \cite{FT89}  by  Foias  and Temam   
as well as to the papers  \cite{LO97,FT98,OT98} where these spaces are used. 
One has  in particular $G^m(\T) = \cup_{\tau > 0} G^m_\tau(\T)$, where 
$$G^m_\tau(\T) \: := \: \left\{ f \in C^\infty(\T), \quad \sup_{j} \, \tau^j \, (j!)^{-m} \, (j+1)^{10} \| f^{(j)} \|_{L^2} < \infty\right\}$$ 
is a Banach space, stable by multiplication. Note that the extra factor   $(j+1)^{10}$ is useful in 
 proving  the stability by multiplication  \cite{LO97,FT98}
(the exponent $10$ is  arbitrary, any power greater than 1 works as well).  

\mspace
In the context of the Prandtl equation, functions depend not only on $x$, but also on $y$. We  just require Sobolev regularity and polynomial decay with respect to the  $y$ variable. 
For $s \in \N, \gamma \ge 0$, we define the   spaces 
$$ H^s_\gamma \: := \: \left\{ g \in H^s(\R_+), \quad  (1+y)^{\gamma+k} g^{(k)} \in  L^2(\R_+), \: k=0...s \right\},    \quad \| g \|_{H^s_\gamma}^2 \: := \: \sum_{k=0}^s \| (1+y)^{\gamma+k} g^{(k)}  \|_{L^2}^2. $$  
We write $L^2_\gamma$ instead of $H^0_\gamma$.           
Accordingly, we introduce the space 
$$ G^m(\T; \, H^s_\gamma ) \: := \: \cup_{\tau > 0}   G^m_\tau(\T; \, H^s_\gamma ) $$
where 
$$G^m_\tau(\T; \, H^s_\gamma ) :=   \left\{ f \in C^\infty(\T; \, H^s_\gamma), \quad \sup_{j} \, \tau^j 
\, (j!)^{-m} \, (j+1)^{10} \| \pa^j_x f \|_{L^2(H^s_\gamma)} < \infty\right\} . $$ 
We shall consider initial data $u_0$ satisfying
$$u_0 \in G^m_{\tau_0}(\T; \, H^{s+1}_{\gamma-1}),   \quad   \omega_0 \: := \: \pa_y u_0 \in   G^m_{\tau_0}(\T; \, H^s_{\gamma})  $$
for some $m, s, \tau_0$ to be fixed later. 
 Our focus will  be on data that are non-monotonic with respect to $y$.  More precisely,  we shall assume  the existence of a single  curve of non-degenerate critical points: 

\mspace
(H) {\em  \quad  $\omega_0(x,y) = 0\: $ iff  $\: y = a_0(x), \: $ for some curve 
 $a_0(x)  > 0$, with $\: \:  \pa_y \omega_0(x,a_0(x)) > 0,$ $ \: \forall x$. }   

\mspace
Without loss of generality, we will assume that $a_0(x) < 3$. 
Besides, to control the behaviour of the flow at large $y$, we shall need  uniform lower and upper bounds on vorticity.  We assume the existence of  $\sigma >0$ and $\delta > 0$  such that  

\mspace
(H') {\em  \quad  For all  $y > 3$,
 for all  $x$,  for all  $\alpha \in \N^2 \: \: |\alpha| \le 2 $,
$$  |\omega_0(x,y)| \ge \frac{ 2 \delta}{(1+y)^\sigma}, \quad| \pa^\alpha \omega_0(x,y) | \le \frac{1}{ 2 \delta (1+y)^{\sigma+\alpha_2}}.$$} 

\mspace
We can now state our main result: 
\begin{theorem} \label{theorem1}
Let $\tau_0 > 0$, $s \ge  8  $ even, $\gamma \ge 1  $,   $\sigma \ge  
  \gamma + \frac12$, $\delta > 0$. Let $u_0$ satisfying 
$$  u_0 \in G^{7/4}_{\tau_0}(\T; \, H^{s+1}_{\gamma-1}),   \quad   \omega_0 \: := \: \pa_y u_0 \in   G^{7/4}_{\tau_0}(\T; \, H^s_{\gamma}),  $$
the compatibility condition: $u_0\vert_{y=0}  = 0$, and assumptions {\em (H), (H')} above. Then, there exists $T > 0$, $0 < \tau \le \tau_0$ and a unique solution 
$$ u \in L^\infty\bigl(0,T;   G^{7/4}_\tau(\T; \, H^{s+1}_{\gamma-1})\bigr), \quad 
\omega \in L^\infty\bigl(0,T;   G^{7/4}_\tau(\T; \, H^s_\gamma)\bigr), $$
of the Prandtl equation \eqref{Prandtl}, with initial data $u_0$. 
\end{theorem}
Remark that  a solution $u$ of \eqref{Prandtl} with the regularity above automatically satisfies 
$$ \pa_t u \in L^\infty\bigl(0,T;   G^{7/4}_{\tau'}(\T; \, H^{s-1}_{\gamma-1})\bigr), \quad \forall \tau' < \tau. $$
This yields continuity in time of $u$ with values in  $G^{7/4}_{\tau'}(\T; \, H^{s-1}_{\gamma-1})$, $\tau' < \tau$, giving a meaning to the initial condition. 

\mspace
Let us point out  that  the solution  $u$ of the theorem  remains in the Gevrey space $G^{7/4}$  in $x$, but   does not  stay {\it a priori} in $G^{7/4}_{\tau_0}$: it is likely that the exponent $\tau$ deteriorates with time. 

\mspace 
The theorem will be the consequence of  the control of some well-chosen Gevrey type norms 
 that evolve in time. More precisely, introducing the vorticity $\omega \: := \: \pa_y u$,  we will control  an  energy of the form $E(t, \tau(t))$, with
$$ E(t,\tau) \: := \: \sum_{j=0}^{+\infty} \left( \tau^j \,    (j!)^{-7/4} \, (j+1)^{10}\right)^2 \, \| \pa_x^j \omega(t,\cdot) \|_{L^2(H^s_\gamma)}^2  $$
 for some positive function $\tau(t)$ decreasing linearly and fast enough with $t$. Actually, it will be slightly better to consider a variant of the previous energy, namely
 \begin{equation} \label{Eomega}
 E_{\omega}(t,\tau) \: := \: \sum_{j=0}^{+\infty} \left( \tau^{j} \,    (j!)^{-7/4} \, (j+1)^{10}\right)^2 \, \| \omega \|_{{\cal H}^j_\gamma}^2, 
 \end{equation}
 with 
$$ \| \omega \|_{{\cal H}^j_\gamma}^2 \: := \:   \sum_{\substack{J = (j_1,j_2) \in \N^2 \\\, |J| = j, \, 0 \le j_2 \le s}}  \| (1+y)^{\gamma+j_2} \pa^J \omega(t,\cdot) \|_{L^2(\T\times \R_+)}^2.$$
 We leave it  to the reader to check that: $\exists  \, c, C > 0$ such that 
 $$E(t,\tau) \: \le \: E_{\omega}(t,\tau), \quad   E_{\omega}(t,\tau) \: \le \: C_{\tau,\tau'} E(t,\tau'), \quad \forall  \, 0 < \tau < \tau' .  $$
Roughly, the choice of $E_\omega$ upon $E$ is connected with estimates on the hyperbolic part of \eqref{Prandtl}, for which derivatives with respect to $x$ and $y$ have the same cost.

 \mspace
Still, the energy $E_\omega$ can not be controlled in a direct manner. As mentioned in  \cite{MW13prep}, there is a problem with estimating $\pa^j_x \omega$ . If the $\pa^j_x  $ hits the $v$ in the transport term $v \pa_y \omega$, 
we end up with  $\pa^j_x  v \pa_y \omega$, which  gives an {\it a priori} loss of one full derivative 
and requires analytic data. This is the main trouble with the Prandtl system, already pointed out in the introduction.  

\mspace
To overcome this difficulty, we shall rely on new estimates, notably inspired by the recent paper  
  \cite{MW13prep}   by  Wong and  the second author.   This paper  is about the well-posedness of the Prandtl equations  for monotonic data, that is when $u_0$ is increasing with $y$. As we recalled in the introduction, in this monotonic case, well-posedness was obtained by Oleinik in the 60's, using the   Crocco transform.  The novelty  in   \cite{MW13prep}   is to obtain such well-posedness result 
 without using  the  Crocco transform,  namely 
 {\em  performing Sobolev estimates  in the original Eulerian formulation}  
 (see also
\cite{AWXY12} where estimates on the linearized problem were  done). 
The main point in the proof is to avoid the loss of derivative generated by the $v$-term. 
One key idea is the following:  combining properly the velocity formulation  \eqref{Prandtl}
and the vorticity formulation ($\omega = \pa_y u$)
\begin{equation}\label{vort}
 \pa_t \omega + u \pa_x \omega + v \pa_y \omega - \pa^2_y \omega = 0, 
\end{equation}
one is left with  an equation of the form  
$$ \pa_t g^{mw} + u \pa_x g^{mw} - \pa^2_y g^{mw} \: = \: \mbox{commutators }  $$
on the new nonlinear quantity: 
$$  g^{mw} \: := \: \omega - \frac{\pa_y \omega}{\omega} u \: = \: \omega \pa_y \frac{u}{\omega}.  $$
The main point with this new equation is that the function $v$, responsible for the loss of one $x$-derivative,  does not appear. Similarly, one can write down equations on 
\begin{equation} \label{gjorigin} 
g_j^{mw} \: := \: \pa_x^j \omega - \frac{\pa_y \omega}{\omega} \pa_x^j u 
\end{equation}
that do not involve the bad term $\pa_x^j v$. Moreover, broadly speaking, one can show that the control of the family $g_j$ in $L^2$  amounts to the control of the family $\pa_x^j \omega$ in $L^2$.  Hence, one can expect to derive a Gronwall type inequality, at the Sobolev level, using the $g_j$'s
(see  \cite{MW13prep}).

\mspace
Of course, the main difference between the present context and the one in  
  \cite{MW13prep}   is assumption (H), that is the existence of a curve of critical points. This critical curve, that reads initially $y = a_0(x)$, with $0 < a_0 < 3$, should  evolve into some $y = a(t,x)$, with $0 < a < 3$  
   for small times.
Differentiating the relation $ \omega(t,x,a(t,x))=0$, 
one finds that the  function $a$ should be governed by the ordinary differential equation in  
the variable $t$:   
\begin{equation}
 \pa_t a(t,x) + \frac{\pa_t \omega(t,x,a(t,x))}{\pa_y \omega(t,x,a(t,x))} = 0, \quad a(0,x) \: = \: a_0(x). 
 \end{equation}
The quantities $g_j^{mw}$'s are not suitable in a neighborhood of this curve. As we shall see,  even away from it, the $\pa^j_x \omega$'s are  not controlled by  the $g_j^{mw}$'s in a suitable way, due to nonlocal phenomena. 

\mspace
Therefore,  we need to introduce an  additional quantity,  that somehow recovers the information lost near the critical curve.  We stress that in this region, the flow should  not be monotonic anymore. 
 However, due to our non-degeneracy assumption, it 
 should remain convex: $\pa^2_y u = \pa_y \omega  > 0$, for $t$ and $y-a(t,x)$ small enough. It turns out that this kind of convexity assumption has been used in a close context, namely in the study of the hydrostatic Euler equations. These equations, set in $\T \times (0,1)$, read
\begin{equation*}
\left\{
\begin{aligned}
\pa_t u + u \pa_x u + v \pa_y u +\pa_x p   =  0, \pa_y p & = 0,  \\
\pa_x u + \pa_y v & = 0, \\
v\vert_{y=0} = v\vert_{y=1} & = 0.   
\end{aligned}
\right.
\end{equation*}
Again, for this system, there is a possible loss of $x$-derivative through the function $v$. As shown in Brenier \cite{Brenier99}, Grenier \cite{Grenier99}, the well-posedness of the hydrostatic equations requires some  convexity assumption.
 This fact was emphasized in another recent paper by  Wong and the second author 
 \cite{MW12}. Starting again from the vorticity equation, 
$$ \pa_t \omega + u \pa_x \omega + v \pa_y \omega   = 0, $$
and dividing by $\sqrt{\pa_y \omega}$, one is left with an equation of the type 
$$ \pa_t h^{mw} + u \pa_x h^{mw}  + \sqrt{\pa_y \omega}\, v =   \mbox{ commutators }$$
on $h^{mw}  \: := \: \frac{\omega}{\sqrt{\pa_y \omega}}$. One can then take advantage of the cancellation 
\begin{equation}  \label{cancelhj} 
\int_{\T \times (0,1)} \sqrt{\pa_y \omega} \, v  h^{mw}  =  \int_{\T \times (0,1)}  v \pa_y u = -\int_{\T \times (0,1)}  \pa_y v u = \int_{\T \times (0,1)} \pa_x \frac{|u|^2}{2}  = 0 
\end{equation}
to get rid of the bad terms in $v$. Let us note that the same idea was used by Grenier in \cite{Grenier00jde} to establish the stability of some characteristic boundary layers.  Similar cancellations hold with higher order derivatives in $x$, through the introduction of 
\begin{equation} \label{hjorigin}
h_j^{mw}  \: := \: \frac{\pa^j_x \omega}{\sqrt{\pa_y \omega}}. 
\end{equation}
Thanks to these quantities, one can obtain local in time Sobolev estimates, like for the Prandtl equations in the monotonic case.

\mspace
In our context, with regards to the previous remarks, it is tempting to replace the original energy $E_\omega$ by a modified one, based on functions $g_j$ and $h_j$ like in \eqref{gjorigin}  and 
 \eqref{hjorigin}. To be more specific, one could think of  combining  two local energies: one away from the critical curve, based on the $g_j$'s, and one in a neighborhood of the critical curve, based on the $h_j$'s. With regards to the recent works \cite{MW12,MW13prep}, one could even expect to obtain stability in Sobolev like spaces.  However, such localization process is not straightforward. Indeed, the Prandtl equation, like the hydrostatic one (see recent developments by Renardy \cite{Renardy09b})
 is highly non-local. This non-locality explains the ill-posedness of the Prandtl equation in the Sobolev setting, as can be seen from the mechanism described by the first author and 
Dormy  \cite{GD10}.

\mspace
At the level of the energy estimates that we will perform, this non-locality will be reflected by some annoying commutator terms. 
To control such bad commutators, we will use  a functional of the following type:
\begin{equation} \label{calE}
{\cal E}(\alpha,t,\tau) \: := \:  \dot{E}_\omega(t,\tau) \: + \:   E_h(t,\tau) \: + \: E^1_g(t,\tau) \: + \: \alpha \,  E^2_g(t,\tau)
\end{equation}
where $\alpha > 0$ is a parameter  to be chosen later. This functional splits into four parts: 
\begin{itemize}
\item The first part is a {\em vorticity energy} 
\begin{equation} \label{vorticityenergy}
\dot{E}_{\omega}(t,\tau) \: := \: \sum_{j\in \N} \left( \tau^{j} \,    (j!)^{-7/4} \, (j+1)^{10}\right)^2 \, \| \omega \|_{\dot{{\cal H}}^j_\gamma}^2\end{equation}
with 
$$ \| \omega \|_{\dot{{\cal H}}^j_\gamma}^2 \: := \:   \sum_{\substack{J = (j_1,j_2) \in \N^2 \\\, |J| = j, \, 0 <  j_2 \le s}}  \| (1+y)^{\gamma+j_2} \pa^J \omega(t,\cdot) \|_{L^2(\T\times \R_+)}^2.$$
The difference with the original energy $E_\omega(t,\tau)$ is the restriction $j_2 > 0$ which 
means that the derivatives $\pa_x^j \om$ are not included here.  For 
each $j$  this energy will provide a  control of all $(x,y)$ derivatives of order $j$ but $\pa_x^{j} \omega$. For $j=0$, one has $\| \omega \|_{\dot{{\cal H}}^j_\gamma} = 0$, but we keep it in the sum for unity. 
\item The second part is a  {\em hydrostatic energy} 
\begin{equation} \label{hydroenergy}
E_h(t,\tau) \: := \: \sum_{j\in \N} \left( \tau^{j} \,    (j!)^{-7/4} \, (j+1)^{10}\right)^2 \, \| h_j(t,\cdot) \|_{L^2(\T\times \R_+)}^2 
\end{equation}
where 
\begin{equation} \label{hj}
 h_j(t,x,y) \: = \: \chi(y-a(t,x)) \, \frac{\pa_x^j \omega}{\sqrt{\pa_y \omega}}(t,x,y),
 \end{equation}
with $\chi = \chi(p) \in C^\infty_c(\R)$   equal to 1 in a neighborhood of $p=0$. We take  $\chi$ with small enough support, so that $\chi(y-a)$ is  compactly supported in  $(0,3)$, and $\pa_y \omega > 0$ over the 
support of $\chi$. This truncation function corresponds to the localization near the critical curve mentioned above. \item The third part is a (first) {\em monotonicity energy}
\begin{equation} \label{monoenergy1}
 E^1_g(t,\tau)  \: = \:  \sum_{j\in \N} \left( \tau^j \,    (j!)^{-7/4} \, (j+1)^{10}\right)^2 \,  \| g_j(t,\cdot) \|_{L^2(L^2_\gamma)}^2 
 \end{equation}
where
\begin{equation} \label{gj}
g_j(t,x,y) \: :=  \: \Big(\psi (y) \omega(t,x,y) + 1-  \psi(y)  \Big)  \left( \pa_x^j \omega - \frac{\pa_y \omega}{\omega} \pa_x^j u \right)(t,x,y)
\end{equation} 
with $\psi = \psi(y) \in C^\infty_c(\R)$, equal to $1$ in an open  neighborhood of  $[0,3]$. 
\mspace 
Note that the truncation $\psi$ makes the quantity $g_j$ well-defined for all $y$, even in the neighborhood of $y = a(t,x)$. Indeed, for large $y$, it amounts to the original definition  \eqref{gjorigin}, whereas near the critical curve, it reads 
$$ g_j(t,x,y) =  \omega \pa_x^j \omega  -   \pa_y\omega \pa_x^j u.  $$
Note that the first term at the right-hand side  vanishes near the critical curve, which leads to a loss of control of $\pa^j_x \omega$ in terms of $g_j$. As explained before, this is why we add the hydrostatic energy to the energy functional. More precisely,  we will show that  the sum of the vorticity energy, the hydrostatic energy, and the first monotonicity energy controls the original functional $E_\omega(t,\tau)$. 
However, we are not able to obtain a closed estimate on  this sum.  
\item Hence, we need to add a second {\em monotonicity energy}
\begin{equation} \label{monoenergy2}
 E^2_g(t,\tau)  \: = \:  \sum_{j\in \N} \left( \tau^j \,    (j!)^{-7/4} \, (j+1)^{10}\right)^2  \: (j+1)^{3/2}  \| \tilde g_j(t,\cdot) \|_{L^2(\T \times \R_+)}^2 
 \end{equation}
 where
\begin{equation}
\tilde g_j(t,x,y) \: := \:  \pa_x^{j-5}\left(  \omega \pa_x^5 \omega  -   \pa_y\omega \pa_x^5 u \right) 
\end{equation}
(with convention $\pa_x^{k} =0$ for $k < 0$). 
Let us remark that  $\tilde g_j$ is  close to $g_j$: for instance, in the region $\psi =1$, one has $\tilde g_j = g_j$ up to commutator terms. Indeed, the replacement of $g_j$ by $\tilde g_j$ is only  a technical issue, that will be explained in due course.  The real key-point in the definition of this second monotonicity energy is  the  extra factor $(j+1)^{3/2}$, creating an anisotropy in the total  energy ${\cal E}(\alpha,t,\tau)$. {\em Such a  choice of  anisotropic energy is the main feature that will allow us to prove stability estimates in the Gevrey setting, below  the analytic case.} 
\end{itemize}

\section{A priori estimates}
The key for Gevrey well-posedness is some {\it a priori} estimate on the anisotropic energy ${\cal E}(\alpha,t,\tau)$ introduced in the previous section. 
This energy involves notably the  functions $g_j$, see  \eqref{gj}, that contains  the factor $\pa_y \omega/\omega$. To control their behaviour  for large $y$,  a lower bound on $\omega$ and an upper bound on $\pa_y \omega$ are needed. More precisely, we shall work with vorticities $\omega$  satisfying: {\em  for all  $y > 3$,   for all   $t,x$, for all  $\alpha \in \N^2, |\alpha| \le 2$,
\begin{equation} \label{lowerboundomega}
 |\omega(t,x,y)| \ge \frac{\delta}{(1+y)^\sigma}, \quad | \pa^\alpha \omega(t,x,y) | \le \frac{1}{\delta (1+y)^{\sigma+\alpha_2}}.
\end{equation} }
Let us point out that this condition is the "all-time" version of condition (H'),  the latter dealing only with  the initial time. Of course, in the end,  we will only assume (H') and 
   we will need to show that such upper and lower bounds are preserved with time (up to a choice of 
 a smaller $\delta$). 

\mspace
\begin{theorem} {\bf (Main a priori estimate)} \label{theorem2}

\sspace
Let $T > 0$, $\tau \ge 1$, $s, \gamma, \sigma$ as in Theorem \ref{theorem1}. Let $u$ be a smooth solution of the Prandtl equation over $]0,T]$, with a single curve of non-degenerate critical points: $y = a(t,x), \: 0 < a < 3 $. Assume that the  vorticity $\omega$ satisfies  \eqref{lowerboundomega} over $]0,T]$, and that 
\begin{equation} \label{EM}
 E_\omega(t,\tau) \: \le \: M, \quad M > 0, \quad \forall t \in ]0,T]. 
\end{equation} 
Then, there exists $\alpha > 0$, $C > 0$ such that for all $t \in ]0,T]$: 
$$ \pa_t  {\cal E}(\alpha,t,\tau) \: \le \: C \, \pa_\tau {\cal E}(\alpha,t,\tau).   $$
\end{theorem}
A close look at the proof will show that  the constants $\alpha$ and  $C$ depend  on $\tau$,  $M$, $\inf a$, $\sup a$,  $\inf_{\{y=a\}} |\pa_y \omega|$, $\inf_{ | y - a | \ge  \eps, y \le 3} |\omega|$  (where, for instance, $\eps := \frac{\inf_{\{y=a\}} |\pa_y \omega|}{2 \sup_{\T \times \R_+} |\pa^2_y \omega|}$), and  the $\delta$ in \eqref{lowerboundomega}.  

\mspace
Let us mention again that in the upper bound $a < 3$ or in the condition $y  > 3$ in \eqref{lowerboundomega}, the choice of the value $3$ is purely arbitrary. 

\mspace
The estimate of the theorem will yield the Gevrey stability, playing on the radius $\tau$ of Gevrey regularity. Indeed, taking some time-dependent $\tau(t)$, we
 observe that 
$$ \pa_t {\cal E}(\alpha,t,\tau(t)) = \pa_t {\cal E}(\alpha,t,\tau(t)) + \tau'(t) \pa_\tau  {\cal E}(\alpha,t,\tau(t))  \leq  (C + \tau'(t))   {\cal E}(\alpha,t,\tau(t)) < 0 $$
if $\tau$ is decreasing fast enough with time. 

\mspace
We insist that all this section is about {\it a priori} estimates. The construction of solutions will require  a further approximation scheme, on which similar estimates will be shown to hold. This will be detailed in later sections. 
 
\subsection{Preliminaries}
 In order to prove Theorem \ref{theorem2}, we need  some extra  notations. 
For any $p\in [1,+\infty[$,  we introduce the weighted $l^p$ space 
 $$ l^p(\tau) \: :=\:  \left\{  (a_j)_{j \in \N}, \quad \sum_{j=0}^{+\infty} \left(\tau^j  (j!)^{-7/4} \, (j+1)^{10}\right)^2    |a_j|^p  \, < +\infty \right\},  $$
with norm
$$ \| a_j \|_{l^p(\tau)} \: := \: \left( \sum_{j=0}^{+\infty}   \left(\tau^j  (j!)^{-7/4} \, (j+1)^{10}\right)^2    |a_j|^p \right)^{1/p}.$$
In particular, we have 
$\dot{E}_\omega(t,\tau) = \Big\|    \| \omega \|_{\dot{{\cal H}}^j_\gamma}   \Big\|^2_{l^2(\tau) } $
and  $ \pa_\tau \dot{E}_\omega(t,\tau) \sim  \Big\|   
  j^{1/2} \| \omega \|_{\dot{{\cal H}}^j_\gamma}   \Big\|^2_{l^2(\tau) }.  $

\medskip
We shall make repeated use  of the following  inequality: 
\begin{lemma} \label{binom}
For all $m \le 5$, 
\begin{equation*}
\| \sum_{k=0}^{\frac{j}{2}}   \binom{j}{k} a_{k+m} \, b_{j-k} \|_{l^2(\tau)} \: \le \: C_\tau \, \| a_j \|_{l^2(\tau)} \,   \| b_j \|_{l^2(\tau)} 
\end{equation*}
and symmetrically, 
\begin{equation*}
\| \sum_{k=\frac{j}{2}}^j   \binom{j}{k} a_k \, b_{j-k+m} \|_{l^2(\tau)} \: \le \: C_\tau \, \| a_j \|_{l^2(\tau)} \,   \| b_j \|_{l^2(\tau)} 
\end{equation*}
\end{lemma}
\begin{remark}
In this lemma and in all the text,  the notation $\frac{j}{2}$ that appears as an index in the sums is slightly abusive:  it stands for the integer part of $\frac{j}{2}$. 
\end{remark}

\mspace
{\em Proof:} Denoting $\alpha_j(\tau) := \tau^j (j!)^{-7/4} (j+1)^{10}$, we find
\begin{align*}
 \| \sum_{k=0}^{\frac{j}{2}}   \binom{j}{k} a_{k+m} \, b_{j-k} \|_{l^2(\tau)}  & \lesssim  
\| \sum_{k=0}^{\frac{j}{2}} \binom{j}{k}^{-3/4} \left(\frac{(k+m)!}{k!}\right)^{7/4} \, (k+m+1)^{-10} \alpha_{k+m}(\tau) a_{k+m} \,  \alpha_{j-k}(\tau)  b_{j-k}  \|_{l^2} \\
 & \lesssim   \| \sum_{k=0}^{+\infty} (k+1)^{\frac{7}{4}m-10} \alpha_{k+m}(\tau) a_{k+m} \,  \alpha_{j-k}(\tau)  b_{j-k}  \|_{l^2} 
 \end{align*}
By a standard convolution inequality for discrete sums, we get 
$$  \| \sum_{k=0}^{\frac{j}{2}}   \binom{j}{k} a_{k+m} \, b_{j-k} \|_{l^2(\tau)}
\: \le \: \| (j+1)^{\frac{7}{4}m-10} \left( \alpha_j(\tau) a_j \right) \|_{l^1} \,   \|  \alpha_j(\tau) b_j \|_{l^2} $$
As $m \le 5$, we can use Cauchy Schwartz inequality to bound the first factor, which yields the result.

\mspace
Let us also emphasize in a lemma some important relations between our energy functionals: 
\begin{lemma} \label{relations}
\begin{equation} \label{relation_energies1}
E_\omega(t,\tau) - \dot{E}_\omega(t,\tau)  \: \lesssim \: E^1_g(t,\tau) + E_h(t,\tau) \: \lesssim \:   E_\omega(t,\tau) - \dot{E}_\omega(t,\tau) 
\end{equation}
and 
\begin{equation} \label{relation_energies2}
 \pa_\tau E_\omega(t,\tau) - \pa_\tau \dot{E}_\omega(t,\tau)  \: \lesssim \:  \pa_\tau  E^1_g(t,\tau) +  \pa_\tau E_h(t,\tau)  \: \lesssim \: \pa_\tau E_\omega(t,\tau) - \pa_\tau \dot{E}_\omega(t,\tau).
 \end{equation}
 \end{lemma}
These relations, to be proved in the appendix, follow from a direct but crucial representation of $\pa^j_x u$ in terms of $g_j$ (see \eqref{gj} for the definition of $g_j$):
\begin{lemma} \label{decompo}
One can write:
\begin{equation} \pa^j_x u(t,x,y) = \left\{ \begin{aligned} & \omega(t,x,y)  \int_3^y \left(\psi + \frac{(1-\psi)}{\omega}\right)^{-1} \frac{g_j}{\omega^2} \: \: + \:C_j(t,x) \omega(t,x,y), \quad y > a(t,x),
  \\ &  \omega(t,x,y)  \int_0^y \left(\psi + \frac{(1-\psi)}{\omega}\right)^{-1} \frac{g_j}{\omega^2} , \quad y < a(t,x)  \end{aligned} \right. \end{equation} 
with $ C_j \: := \: -\pa_x^j u(t,x,3)/\omega(t,x,3)$. 
\end{lemma}
{\em Proof:}  The proof is trivial: integrate the relation 
\begin{equation} \label{gjbis}
\pa_y \frac{\pa^j_x u}{\omega}  \: = \:  \left(\psi + \frac{(1-\psi)}{\omega}\right)^{-1} \frac{g_j}{\omega^2} \end{equation} 
from $0$ to $y$ and  from 
 $3$ to $y$ respectively. We  point out that  the choice $3$ is arbitrary. 

\mspace
We stress  again that   $u$ is assumed to be a smooth solution of the Prandtl equation, notably smooth across the critical curve $y = a(t,x)$.  Note also that  \eqref{gjbis}  formally yields the relation
$$  \pa^j_x u(t,x,y) \: = \:   \omega(t,x,y)  \int_0^y \left(\psi + \frac{(1-\psi)}{\omega}\right)^{-1} \frac{g_j}{\omega^2}  $$
 {\em for all y}.  However, the right-hand side does not make sense
 for $y > a(t,x)$:  as $\omega$ degenerates near $y = a$, the integral is not properly defined. One can not use it to  bound  nicely $\pa^j_x u$ in terms of $g_j$, even away from $y = a$.  This is why we substitute to this formula the decomposition 
 $$ \pa^j_x u(t,x,y) = \omega(t,x,y)  \int_3^y \left(\psi + \frac{(1-\psi)}{\omega}\right)^{-1} \frac{g_j}{\omega^2}  \: + \:C_j(t,x) \omega(t,x,y) $$
 with a first term that depends nicely on $g_j$ away from $y=a$, and a second term that, broadly speaking, will be controlled by the hydrostatic energy. 

\mspace
Let us conclude this paragraph by a trivial remark:  for any energy functional $E = E_\omega, E_h \dots$ that  will be used in the text, one has $E \le \pa_\tau E$. This will be used many times without mentioning.   

\subsection{Estimate of the vorticity energy} \label{subsecvorticity}
We focus here on the  energy $\dot{E}_\omega(t,\tau)$ defined in \eqref{vorticityenergy}. We also introduce  the functional
\begin{equation} \label{vorticitydissipation}
\dot{D}_{\omega}(t,\tau) \: := \: \sum_{j=0}^{+\infty} \left( \tau^j \,    (j!)^{-7/4} \, (j+1)^{10}\right)^2 \, \|  \pa_y \omega(t,\cdot) \|_{\dot{{\cal H}}^j_\gamma}^2. 
\end{equation}
which will appear due to to the viscous term in the Prandtl equation. We shall establish
\begin{proposition} {\bf (Estimate on the vorticity energy)} \label{propvorticity}

\sspace
Let $T > 0$, $\tau \ge 1$, $s, \gamma, \sigma$ be  as in Theorem \ref{theorem1}. Let $u$ be a smooth solution of the Prandtl equation over $]0,T]$. Assume that the  vorticity $\omega$ satisfies  \eqref{lowerboundomega} over $]0,T]$, and that 
$$ E_\omega(t,\tau)  \: \le \: M, \quad M > 0, \quad \forall t \in ]0,T]. $$
Then, one has for some $C > 0$ and all $t \in ]0,T]$: 
$$ \pa_t \dot{E}_\omega(t,\tau) \: + \: \dot{D}_\omega(t,\tau) \: \le \: C \, \pa_\tau E_\omega(t,\tau).  $$
\end{proposition}
 
\mspace
The starting point is to write down an equation on $\omega_J \: := \: \pa^{J} \omega$, $ J = (j_1,j_2) \in \N^2$, \quad $0 < j_2 \le s$. Differentiating the vorticity equation \eqref{vort}, 
we find 
$$ \pa_t \omega_J  + u \pa_x  \omega_J + v \pa_y \omega_J - \pa^2_y \omega_J =  -[\pa^J, u] \pa_x \omega -[\pa^J, v] \pa_y \omega.  $$
After multiplication by  $ (1+y)^{2\gamma + 2j_2} \omega_J$, integration over $\T \times \R_+$  and standard integration by parts, we get
\begin{align*}
& \frac{1}{2} \, \frac{d}{dt} \| (1+y)^{\gamma+j_2} \omega_J(t,\cdot) \|_{L^2(\T\times \R_+)}^2 \: + \: \| (1+y)^{\gamma+j_2} \pa_y \omega_J(t,\cdot) \|_{L^2(\T\times \R_+)}^2 \\
 = &\: -\int_{\T \times \R_+} (1+y)^{2\gamma + 2j_2}  [\pa^J, u] \pa_x \omega  \, \omega_J \: - \: \int_{\T \times \R_+} (1+y)^{2\gamma + 2j_2} [\pa^J, v] \pa_y \omega \, \omega_J \\
 &\: -   (2\gamma + 2j_2) \int_{\T \times \R_+} (1+y)^{2\gamma+2j_2-1} \left( \pa_y \omega_J - v \omega_J \right) \, \omega_J \: + \: \int_{\T \times \{0\}} \pa_y \omega_J \, \omega_J. 
\end{align*} 
Multiplying by $ \left(\tau^{|J|} (|J|!)^{-7/4} \, |J|^{10}\right)^2$  and summing over $J \in \N \times [|1,s|]$, we obtain 
\begin{align*}
 \pa_t \dot{E}_\omega(t,\tau) \: &+ \: \dot{D}_\omega(t,\tau) \:  \\ 
& \lesssim \:  \left( \| A_j \|_{l^2(\tau)} + \| B_j \|_{l^2(\tau)}  \: + \:  \| C_j \|_{l^2(\tau)}+  
\| D_j \|_{l^2(\tau)}  \right) \| j^{1/2} \| \omega \|_{\dot{{\cal H}}^j_\gamma} \|_{l^2(\tau)}   \:  + \: \| E_j \|_{l^1(\tau)}   \\
&   \lesssim   \:  \left( \| A_j \|_{l^2(\tau)} + \| B_j \|_{l^2(\tau)}  \: + \:  \| C_j \|_{l^2(\tau)} +  \| D_j \|_{l^2(\tau)}  \right)\sqrt{ \pa_\tau \dot{E}_\omega(t,\tau)}  
 \:  + \: \| E_j \|_{l^1(\tau)} 
\end{align*}
where 
\begin{align}
A_j^2 \: &:= \:  \frac{1}{j}  \sum_{\substack{J = (j_1,j_2) \in \N^2 \\\, |J| = j, \, 0 <  j_2 \le s}}  \| (1+y)^{\gamma + j_2}  [\pa^J, u] \pa_x \omega \|_{L^2(\T \times \R_+)}^2, \\
B_j^2 \: &:= \:  \frac{1}{j}  \sum_{\substack{J = (j_1,j_2) \in \N^2 \\\, |J| = j, \, 0 <  j_2 \le s}}   \| (1+y)^{\gamma + j_2}  [\pa^J, v] \pa_y \omega \|_{L^2(\T \times \R_+)}^2, \\
C_j^2 \: &:= \:   \frac{1}{j} \sum_{\substack{J = (j_1,j_2) \in \N^2 \\\, |J| = j, \, 0 <  j_2 \le s}}   \| (1+y)^{\gamma + j_2-1} v \,   \omega_J  \|_{L^2(\T \times \R_+)}^2, \\
D_j^2 \: &:= \:   \frac{1}{j} \sum_{\substack{J = (j_1,j_2) \in \N^2 \\\, |J| = j, \, 0 <  j_2 \le s}}   \| (1+y)^{\gamma + j_2-1} \pa_y \omega_J  \|_{L^2(\T \times \R_+)}^2, \\
E_j \: & =  \:  \sum_{\substack{J = (j_1,j_2) \in \N^2 \\\, |J| = j, \, 0 <  j_2 \le s}}   \int_{\T \times \{0\}} \pa_y \omega_J \, \omega_J. 
\end{align}
Note that   the first term at the right-hand side. comes from Cauchy-Schwarz inequality, first in $L^2(\T \times \R_+)$, then in  $l^2(\N \times [|1, s|])$.  

\subsubsection*{Estimate on $A_j$}
We first write, for  $J \in \N \times [|1,s|]$, $|J|= j$: 
$$  (1+y)^{\gamma + j_2} [\pa^J, u] \pa_x \omega \: = \: \sum_{\substack{K \le J \\ |K| > 0}} \binom{J}{K} (1+y)^{k_2} \pa^K u \  (1+y)^{\gamma+j_2-k_2} \, \pa^{J-K}\pa_x \omega.  $$
As $\gamma \ge 1$,  we obtain  
$$  (1+y)^{\gamma + j_2}  \bigl|   [\pa^J, u] \pa_x \omega \bigr| \le \sum_{\substack{K \le J \\ |K| > 0}} \binom{J}{K}  \bigl| (1+y)^{\gamma+k_2-1}   \pa^{K} u \bigr| \,  \bigl|  (1+y)^{\gamma+j_2-k_2} \pa^{J-K} \pa_x \omega\bigr|.  $$ 
Then, we write 
\begin{multline}
\| (1+y)^{\gamma + j_2}  [\pa^J, u] \pa_x \omega \|_{L^2(\T \times \R_+)} \\
 \le  \sum_{\substack{K \le J \\ \frac{|J|}{2} \ge |K| > 0}}  \binom{J}{K}  a^1_K \, b^1_{J-K}  
 \:  + \: \sum_{\substack{K \le J \\  J \ge |K| > \frac{|J|}{2}}}  \binom{J}{K}  a^2_K \, b^2_{J-K}  
 \end{multline}
where 
\begin{align*}
& a^1_K \: :=  \: \|  (1+y)^{\gamma+k_2-1} \pa^K u \|_{L^\infty(\T \times \R_+)}, \quad b^1_{J-K}   \: :=  \: \|  (1+y)^{\gamma+j_2- k_2} \pa^{J-K} \pa_x \omega  \|_{L^2(\T \times \R_+)} \\ 
& a^2_K \: :=  \: \|  (1+y)^{\gamma+k_2-1} \pa^K u \|_{L^2(\T \times \R_+)}, \quad \: b^2_{J-K}   \: :=  \: \|  (1+y)^{\gamma+j_2- k_2} \pa^{J-K} \pa_x \omega  \|_{L^\infty(\T \times \R_+)} 
\end{align*}

\mspace
{\em Treatment of the first sum.} 

\sspace
We first notice  that $\binom{J}{K} \lesssim \binom{|J|}{|K|}$, uniformly for $J \in \N \times [|1,s|]$. 
Then, using  Sobolev  and Hardy  inequalities \eqref{sobolev}, \eqref{hardy1},  we obtain easily that 
$$ a^1_K \: \lesssim  \: \| \omega \|_{{\cal H}^{|K|}_\gamma} \: + \:  \| \omega \|_{{\cal H}^{|K|+1}_\gamma}. $$
Also, clearly 
$$ b^1_{J-K} \: \lesssim \:  \| \omega \|_{{\cal H}^{|J-K|+1}_\gamma}.  $$
Hence, we  obtain ($k = |K|, \, j = |J|$): 
\begin{align*} 
&  \frac{1}{j^{1/2}} \sum_{\substack{K \le J \\ \frac{|J|}{2} \ge |K| > 0}}  \binom{J}{K}  a^1_K \, b^1_{J-K}     \:  \lesssim \: \frac{1}{j^{1/2}} \sum_{k=1}^{\frac{j}{2}} \binom{j}{k}  \| \omega \|_{{\cal H}^k_\gamma}  \,    \| \omega \|_{{\cal H}^{j-k+1}_\gamma}   \: + \:  \frac{1}{j^{1/2}}  \sum_{k=1}^{\frac{j}{2}} \binom{j}{k}  \| \omega \|_{{\cal H}^{k+1}_\gamma} \,  \| \omega \|_{{\cal H}^{j-k+1}_\gamma}  \\
& \lesssim  \:  \frac{1}{j^{1/2}} \sum_{k=0}^{\frac{j}{2}-1} \binom{j}{k+1}  \| \omega \|_{{\cal H}^{k+1}_\gamma}  \,    \| \omega \|_{{\cal H}^{j-k}_\gamma}   \: + \:  \frac{1}{j^{1/2}}  \sum_{k=0}^{\frac{j}{2}-1} \binom{j}{k+1}  \| \omega \|_{{\cal H}^{k+2}_\gamma} \,  \| \omega \|_{{\cal H}^{j-k}_\gamma}  \\
&  \lesssim  \:  \sum_{k=0}^{\frac{j}{2}} \binom{j}{k}  \| \omega \|_{{\cal H}^{k+1}_\gamma}  \,  \left( (j-k)^{1/2}  \| \omega \|_{{\cal H}^{j-k}_\gamma}  \right)  \: + \:    \sum_{k=0}^{\frac{j}{2}} \binom{j}{k}  \| \omega \|_{{\cal H}^{k+2}_\gamma} \, \left(  (j-k)^{1/2}  \| \omega \|_{{\cal H}^{j-k}_\gamma} \right). 
\end{align*}
We finally control the $l^2(\tau)$ norm of the right-hand side using Lemma \ref{binom} (we take $m=1$ for the first term, $m=2$ for the second term). We end up with 
\begin{equation}  \label{estimAj1} 
\|  \frac{1}{j^{1/2}}  \sum_{\substack{J = (j_1,j_2) \in \N^2 \\\, |J| = j, \, 0 <  j_2 \le s}} \: \sum_{\substack{K \le J \\ \frac{|J|}{2} \ge |K| > 0}}  \binom{J}{K}  a^1_K \, b^1_{J-K}  \|_{l^2(\tau)} \:  \lesssim  \: \sqrt{E_\omega(t,\tau)} \, \sqrt{\pa_\tau E_\omega(t,\tau)} \: \lesssim  \:  \sqrt{\pa_\tau E_\omega(t,\tau)} 
\end{equation}

\mspace
{\em Treatment of the second sum.} 

\sspace
We proceed as for the first sum, reversing the role of factors $a$ and $b$: more precisely, we apply the Hardy inequality \eqref{hardy1} to $a^2_K$:
$$ a^2_K \: \lesssim  \:  \| \omega \|_{{\cal H}^{|K|}_\gamma}  $$
 and the Sobolev bound \eqref{sobolev} to  $b^2_{J-K}$. Note that, due to the last two terms at the right-hand side of \eqref{sobolev}, the quantities $ \pa_y \pa^{J-K} \pa_x \omega$ and $ \pa_y \pa^{J-K} \pa^2_x \omega$  are needed to control  $b^2_{J-K}$. In the special case where $j_2 = s$, $k_2 = 0$, they involve $s+1$ $\ \pa_y$-derivatives. Thus,  they can not be controlled by $\| \omega \|_{{\cal H}^{|J-K|+m}_\gamma}$, similarly to what we did for $a^1_K$ (the definition of ${\cal H}^j_\gamma$ spaces involves only $j_2 \le s$   $\, \pa_y$-derivatives).  
Hence, in this special case, we rather use the inequality 
$$ \| (1+y)^{\gamma+j_2-k_2} \pa_y \pa^{J-K}  \pa_x^m \omega \|_{L^2} \: \lesssim \: \| \pa_y \omega \|_{\dot{\cal H}^{|J-K|+m}}, \quad m=1,2.$$   
We insist that the r.h.s involves the homogeneous space $\dot{\cal H}^{|J-K|+m}$, because $j_2 - k_2 \neq 0$ ($j_2 = s, k_2=0$). This is important to use the dissipation term $\dot{D}_\omega$.  Eventually, we find 
$$  b^2_{J-K} \: \lesssim \: \sum_{m=1}^3    \| \omega \|_{{\cal H}^{|J-K|+m}_\gamma} \: + \:    \| \pa_y \omega \|_{\dot{{\cal H}}^{|J-K|+m}_\gamma}  . $$
Then,  
\begin{equation*}
 \frac{1}{j^{1/2}}  \sum_{\substack{K \le J\\  J \ge |K| > \frac{|J|}{2}}}  \binom{J}{K}  a^2_K \, b^2_{J-K}  
\: \le \:  \frac{1}{j^{1/2}} \sum_{m=1}^3 \sum_{k=\frac{j}{2}}^j   \binom{j}{k}  \| \omega \|_{{\cal H}^k_\gamma}  \, \left( \| \omega \|_{{\cal H}^{j-k+m}_\gamma} \: + \:  \| \pa_y \omega \|_{\dot{{\cal H}}^{j-k+m}_\gamma} \right)
\end{equation*} 
We crudely bound $j^{-1/2}$ by $1$, and apply Lemma \ref{binom}, to find 
\begin{equation}  \label{estimAj2} 
\begin{aligned}
\|  \frac{1}{j^{1/2}}  \sum_{\substack{J = (j_1,j_2) \in \N^2 \\ |J| = j, \, 0 <  j_2 \le s}} \: 
 \sum_{\substack{K \le J \\  J \ge |K| > \frac{|J|}{2}}}  \binom{J}{K}  a^2_K \, b^2_{J-K} \|_{l^2(\tau)} & \lesssim \sqrt{E_\omega(t,\tau)} \left(  \sqrt{E_\omega(t,\tau)} \: + \:  \sqrt{\dot{D}_\omega(t,\tau)} \right) \\
 & \lesssim   \sqrt{E_\omega(t,\tau)} \: + \:  \sqrt{\dot{D}_\omega(t,\tau)} 
\end{aligned}
\end{equation}
Gathering of the bounds \eqref{estimAj1}-\eqref{estimAj2} leaves us with 
\begin{equation} \label{estimAj}  \| A_j \|_{l^2(\tau)} \:  \lesssim  \:   \, \sqrt{E_\omega(t,\tau)} +\sqrt{\pa_\tau E_\omega(t,\tau)} \: + \: \sqrt{\dot{D}_\omega(t,\tau)} \:  \lesssim \:   \sqrt{\pa_\tau E_\omega(t,\tau)} \: + \: \sqrt{\dot{D}_\omega(t,\tau)} 
\end{equation} 
\subsubsection*{Estimate on $B_j$}
Let us first point out  that    $v  = -\int_0^{y} \pa_x u$ does not decay at $y=\infty$. One has   
\begin{equation} \label{v_infty}
\begin{aligned}
 \| \pa^k_x v \|_{L^\infty} &  \lesssim \| \int_0^{+\infty} |\pa^{k+1}_x u| \, dy \|_{L^\infty(\T)} \\
&  \lesssim \| \pa^{k+1}_x u \|_{L^\infty(L^2_{\gamma-1})}  
\: \lesssim \: \| \pa^{k+1}_x u \|_{L^2(L^2_{\gamma-1})} \: + \:  \| \pa^{k+2}_x u \|_{L^2(L^2_{\gamma-1})} \\
& \lesssim \: \| \pa^{k+1}_x \omega \|_{L^2(L^2_\gamma)} \: + \:  \| \pa^{k+2}_x \omega \|_{L^2(L^2_\gamma)}
\end{aligned}
\end{equation}
where the last line comes  from \eqref{hardy2}. For $K = (k_1,k_2)$ with $k_2 > 0$, we find 
\begin{equation} \label{Djv_infty}
\|(1+y)^{\gamma + k_2 - 1} \pa^K v \|_{L^\infty} \:  =  \: \|(1+y)^{\gamma + k_2 - 1} \pa^{k_1+1}_x \pa^{k_2 -1}_y u \|_{L^\infty} 
\end{equation}
which, as seen before, can be controlled through \eqref{sobolev} and \eqref{hardy1}.

\mspace
Proceeding (almost) as for $A_j$, we get
\begin{multline}
\| (1+y)^{\gamma + j_2} [\pa^J, v] \pa_y \omega  \|_{L^2(\T \times \R_+)} \\
 \le  \sum_{\substack{K \le J \\ \frac{|J|}{2} \ge |K| > 0}}  \binom{J}{K}  a^1_K \, b^1_{J-K}  
 \:  + \: \sum_{\substack{K \le J, \\ j_2 - k_2 \le s-2 \\  J \ge |K| > \frac{|J|}{2}}}  \binom{J}{K}  a^2_K \, b^2_{J-K}  \:  + \: \sum_{\substack{K \le J, \\ s-1 \le  j_2 - k_2 \le s \\  J \ge |K| > \frac{|J|}{2}}}  \binom{J}{K}  a^3_K \, b^3_{J-K}.  
 \end{multline}
where this time
\begin{align*}
& a^1_K \: :=  \: \|  (1+y)^{k_2} \pa^K v \|_{L^\infty(\T \times \R_+)}, \quad b^1_{J-K}   \: :=  \: \|  (1+y)^{\gamma+j_2- k_2} \pa^{J-K} \pa_y \omega  \|_{L^2(\T \times \R_+)} \\ 
& a^2_K \: :=  \: \|  (1+y)^{k_2-1} \pa^K v \|_{L^2(\T \times \R_+)},  \: b^2_{J-K}   \: :=  \: \|  (1+y)^{\gamma+j_2- k_2+1} \pa^{J-K} \pa_y \omega  \|_{L^\infty(\T \times \R_+)}, \\
& a^3_K \: :=  \: \|  (1+y)^{k_2} \pa^K v \|_{L^\infty(\T \times \R_+)}, \quad  b^3_{J-K}   \: :=  \: \|  (1+y)^{\gamma+j_2- k_2} \pa^{J-K} \pa_y \omega  \|_{L^2(\T \times \R_+)}. 
\end{align*}

\mspace
{\em Treatment of the first sum}

\sspace
Thanks to  \eqref{v_infty} and \eqref{Djv_infty}, we find 
$$ a^1_K \: \lesssim \: \sum_{m=1}^2 \| \omega \|_{{\cal H}^{|K|+m}_\gamma}  . $$
As regards $b^1_K$, one must take care of the special case: $j_2 = s$, $k_2 = 0$,  for which  $\pa^{J-K} \pa_y \omega$ involves $s+1$ derivatives  with respect to $y$. Hence, we write 
$$ \| b^1_K \| \: \le \: \|  \omega \|_{{\cal H}^{|J-K|+1}_\gamma} \: + \: \| \pa_y \omega \|_{\dot{{\cal H}}^{|J-K|}_\gamma} $$
where the last term at the right-hand side accounts for the special case.  
Proceeding as for $A_j$, relying on Lemma \ref{binom}, we find 
\begin{equation} \label{estimBj1}
\|  \frac{1}{j^{1/2}}  \sum_{\substack{J = (j_1,j_2) \in \N^2 \\\, |J| = j, \, 0 <  j_2 \le s}} \: \sum_{\substack{K \le J \\ \frac{|J|}{2} \ge |K| > 0}}  \binom{J}{K}  a^1_K \, b^1_{J-K}  \|_{l^2(\tau)} \lesssim \sqrt{\pa_\tau E_\omega(t,\tau)} \: + \: \sqrt{\dot{D}_\omega(t,\tau)}.
\end{equation}

\mspace
{\em Treatment of the second sum}

\sspace
We handle the second sum as we did for $A_j$: we apply  Hardy inequality to $a^2_K$, Sobolev inequality to $b^2_K$. Remark that we put factor  
$(1+y)^{k_2-1}$ in front of $\pa^K v$, in the definition of $a^2_K$. It allows the $L^2$ norm to be finite even in the case $k_2 = 0$.  Of course, this forces a factor $(1+y)^{\gamma + j_2-k_2+1}$ in front of $\pa^{J-K} \pa_y \omega$, which is harmless thanks to the extra $\pa_y$. 

\mspace
Moreover, as we restrict here to indices satisfying $j_2 - k_2 \le s-2$, the control of $b^2_K$ by the inequality \eqref{sobolev} (which involves for instance  $\pa^{J-K} \pa^2_y \omega$) only requires terms with less  than $s$ derivatives in $y$. We leave the details to the reader. We get 
\begin{equation} \label{estimBj2}
\|  \frac{1}{j^{1/2}}  \sum_{\substack{J = (j_1,j_2) \in \N^2 \\\, |J| = j, \\ 0 <  j_2 \le s}} \:  \sum_{\substack{K \le J, \\ j_2 - k_2 \le s-2 \\  J \ge |K| > \frac{|J|}{2}}}   \binom{J}{K}  a^2_K \, b^2_{J-K}  \|_{l^2(\tau)} \: \lesssim \:  E_\omega(t,\tau) \: \lesssim \: \sqrt{E_\omega(t,\tau)}
\end{equation}
 
\mspace
{\em Treatment of the third sum}

\sspace
Note that the third sum is empty except when $j_2 = s$ (in which case $k_2 \in \{0,1\}$) or when $j_2 =s-1$ (in which case $k_2=0$).  In both cases, the important thing to notice is that 
\begin{equation} \label{remarque_K}
|K| \: \le \:  k_1 \: + \: 1  \: \le \: j_1 - 1 \: \le \: |J| - s, \quad \mbox{ and  }\:  \binom{J}{K} \le \binom{|J|}{|K|+m} \:\: \forall m \le s  
\end{equation}
Using once again \eqref{v_infty}, \eqref{Djv_infty} and \eqref{sobolev}, we end up with 
$$  a^3_K \: \lesssim \: \sum_{m=1}^2 \| \omega \|_{{\cal H}^{|K|+m}_\gamma} $$
whereas 
$$ b^3_K \: \lesssim \:  \| \pa_y \omega \|_{\dot{{\cal H}}^{|J-K|}_\gamma}. $$
We  get, taking \eqref{remarque_K} into account: 
\begin{equation*}
 \frac{1}{j^{1/2}}  \sum_{\substack{K \le J, \\  s-1 \le j_2 - k_2 \le s \\  J \ge |K| > \frac{|J|}{2}}} \binom{J}{K}  a^3_K \, b^3_{J-K}  
 \lesssim \:  \frac{1}{j^{1/2}} \sum_{m=1}^2 \sum_{k=\frac{j}{2}}^{j-s}   \binom{j}{k+m}  \| \omega \|_{{\cal H}^{k+m}_\gamma}  \,  \| \pa_y \omega \|_{\dot{{\cal H}}^{j-k}_\gamma} 
\end{equation*} 
We can bound  $\frac{1}{j^{1/2}}$ by $1$, and make a  change of index: $k' = k+m$. Then, application of Lemma  \ref{binom} leads to 
\begin{equation} \label{estimBj3}
\| \frac{1}{j^{1/2}}  \sum_{\substack{K \le J, \\  s-1 \le j_2 - k_2 \le s \\  J \ge |K| > \frac{|J|}{2}}} \binom{J}{K}  a^3_K \, b^3_{J-K}  \|_{l^2(\tau)}
\\ \lesssim \: \sqrt{E_\omega(t,\tau)} \sqrt{\dot{D}_\omega(t,\tau)}  \: \lesssim \: \sqrt{\dot{D}_\omega(t,\tau)}
\end{equation}
Gathering of the bounds \eqref{estimBj1}-\eqref{estimBj2}-\eqref{estimBj3} leads to
\begin{equation} \label{estimBj}  \| B_j \|_{l^2(\tau)} \:  \lesssim  \: \sqrt{E_\omega(t,\tau)} +\sqrt{\pa_\tau E_\omega(t,\tau)} \: + \: \sqrt{\dot{D}_\omega(t,\tau)} \: \lesssim \: \sqrt{\pa_\tau E_\omega(t,\tau)} \: + \: \sqrt{\dot{D}_\omega(t,\tau)}
\end{equation}

\subsubsection*{Estimate on $C_j$}
Clearly, 
$$  | C_j | \: \le \: \| \frac{v}{(1+y)} \|_{L^\infty} \| \omega \|_{\dot{{\cal H}}^j_\gamma} \: \lesssim \sqrt{E_\omega(t,\tau)} \| \omega \|_{\dot{{\cal H}}^j_\gamma},  $$ 
see \eqref{v_infty}. Thus, 
\begin{equation}  \label{estimCj}
\| C_j \|_{l^2(\tau)} \: \lesssim \sqrt{E_\omega(t,\tau)}  \sqrt{\dot{E}_\omega(t,\tau)} \: \lesssim \:  \sqrt{E_\omega(t,\tau)}
\end{equation}

\subsubsection*{Estimate on $D_j$}
Clearly, 
\begin{equation} \label{estimDj}
 \| D_j \|_{l^2(\tau)} \:  \lesssim  \: \sqrt{\dot{D}_\omega(t,\tau)}
\end{equation}

\subsubsection*{Estimate on $E_j$}
To handle the boundary term, a simple application of the trace theorem is not enough. One shall adapt ideas from \cite{MW13prep}. The main point is to reduce the number of derivatives in the boundary  term $\pa_y \omega_I \, \omega_I\vert_{y=0}$, thanks to the equation. For instance, one can observe that 
$$ \pa_y \omega\vert_{y=0} = \left( \pa_t u + u \pa_x u + v \pa_y u \right)\vert_{y=0} = 0. $$
 Then, 
 \begin{equation} \label{reductiond3omega}
  \pa^3_y \omega\vert_{y=0}  = \pa_y \left( \pa_t \omega + u \pa_x \omega + v \pa_y \omega \right)\vert_{y=0} = \omega \pa_x \omega\vert_{y=0}.  
  \end{equation}
 For higher derivatives, one has
\begin{lemma}  \label{reductionlemma} {\bf(from \cite[Lemma 5.9]{MW13prep})}
 
 \sspace
 For $j_2 \ge 4$ an even number, $\pa_y \omega_J$ is a linear combination of terms of the form 
 $$ \pa_x^{j_1} \left(  \prod_{l=1}^N \pa_x^{\alpha_l} \pa_y^{\beta_l} \omega\vert_{y=0} \right) $$
 with
 $$ 2 \le N \le \frac{j_2}{2}, \quad 
 \left\{\begin{aligned}
 &\sum_{l=1}^N 3 \alpha_l + \beta_l = j_2 + 1, \\
& \sum_{l=1}^N \alpha_l \le \frac{j_2}{2} - 1, \quad   \sum_{l=1}^N \beta_l \le j_2 - 2, \\
& \alpha_l + \beta_l \le j_2 - 1, \quad \forall l=1...N.  
\end{aligned}
\right.
$$
\end{lemma}

\mspace
Besides this lemma, we need a slight generalization of Lemma \ref{binom}, whose proof is left to the reader: 
\begin{lemma} \label{multinom}
Let $N \ge 2$, $m_2, \dots, m_N \le 5$. Let 
$$ K_1(j) \: := \: \left\{ (k_1,\dots,k_N) \in \N^N, \quad \mbox{ s.t. } \: k_1 + \dots + k_n = j, \quad k_1  \ge \frac{j}{n} \right\} $$
for all $j \in \N$. 
Then, for sequences $a^l_j$, $j \in \N$,    $l=1...N$, one has 
$$ \left\| \sum_{(k_1,...,k_N) \in K_1(j)}     \frac{j!}{\prod_{l=1}^N k_l !} \, a^1_{k_1} \, \prod_{l=2}^N a^l_{k_l + m_l} \right\|_{l^2(\tau)}  \: \lesssim \: \prod_{l=1}^N \| a_j^l \|_{l^2(\tau)}. 
$$
\end{lemma}

\mspace
We now write 
\begin{align*}
E_j \: & = \: \sum_{\substack{J  \in \N^* \times 2\N, \\ |J| = j, \, 0 <  j_2 \le s}}   \int_{\T \times \{0\}} \pa_y \omega_J \, \omega_J \: + \: \sum_{\substack{J  \in \N^* \times (2\N+1), \\  |J| = j, \, 0 <  j_2 \le s}}    \int_{\T \times \{0\}} \pa_y \omega_J \, \omega_J \\
& \quad + {\bf 1}_{[| 1, \dots, s|]}(j)  \int_{\T \times \{0\}} \pa_y^{j+1} \omega \, \pa_y^j \omega \: := \: E^1_j \: + \: E^2_j \: + \: E^3_j 
\end{align*}

\mspace
{\em Study of $E^1_j$}. Let $J  \in \N^* \times 2\N, \, |J| = j, \, 0 <  j_2 \le s$. We write 
\begin{align*}
 \left| \int_{\T \times \{0\}} \pa_y \omega_J \, \omega_J \right| \: & \le \: \| \pa_y \omega_J\vert_{y=0} \|_{L^2(\T)} \,    \| \omega_J\vert_{y=0} \|_{L^2(\T)} \\ & \le \:  \| \pa_y \omega_J\vert_{y=0} \|_{L^2(\T)}  \, \| \omega_J \|_{L^2(\T \times \R_+)}^{1/2} \, \| \pa_y \omega_J \|_{L^2(\T \times \R_+)}^{1/2} 
   \le \: \| \pa_y \omega_J\vert_{y=0} \|_{L^2(\T)} \, \| \pa_y \omega \|_{\dot{{\cal H}}^j_{\gamma}} 
 \end{align*}
 applying the trace theorem and then the Hardy inequality \eqref{hardy1} to the last factor. As regards the first factor, one must use the reductions seen above. We focus on the case  $j_2 \ge 4$ (and even), which is treated thanks to Lemma \ref{reductionlemma}. The case $j_2= 2$, involving \eqref{reductiond3omega}, is simpler.  The $L^2$ norm of $\| \pa_y \omega_J\vert_{y=0} \|_{L^2(\T)}$ can be bounded by a finite number of terms of the type  $\|  \pa_x^{j_1} \left(  \prod_{l=1}^N \pa_x^{\alpha_l} \pa_y^{\beta_l} \omega\vert_{y=0} \right) \|_{L^2(\T)}$, where $N$ and the $(\alpha_l,\beta_l)$'s satisfy the conditions of the lemma. 
 
 \mspace
 As usual, we use the Lebnitz formula to write 
 \begin{multline*} 
   \|  \pa_x^{j_1} \left(  \prod_{l=1}^N \pa_x^{\alpha_l} \pa_y^{\beta_l} \omega\vert_{y=0} \right) \|_{L^2(\T)} 
   = \| \sum_{k_1 + \dots + k_n = j_1}       \frac{j_1!}{\prod_{l=1}^N k_l !}  \prod_{l=1}^N \pa_x^{\alpha_l +k_l} \pa_y^{\beta_l} \omega\vert_{y=0} \|_{L^2(\T)}  \\
   \le \sum_{l'=1}^N  \| \sum_{(k_1, \dots, k_n) \in K_{l'}(j_1)}      \frac{j_1!}{\prod_{l=1}^N k_l !}  \prod_{l=1}^N \pa_x^{\alpha_l +k_l} \pa_y^{\beta_l} \omega\vert_{y=0} \|_{L^2(\T)}  \: := \: \sum_{l'=1}^N {\cal N}_{l'}  
  \end{multline*}
  where 
  $$ K_l(j_1) \: := \: \left\{ (k_1,\dots,k_N) \in \N^N, \quad \mbox{ s.t. } \: k_1 + \dots + k_N = j_1, \quad k_l  \ge \frac{j_1}{N} \right\}. $$
 We now bound ${\cal N}_1$, the other terms being treated in the same way. We write
 $${\cal N}_1 \: \lesssim \: \sum_{(k_1,\dots,k_N) \in K_1(j_1)}    \frac{j_1!}{\prod_{l=1}^N k_l !}  \, 
 \| \pa_x^{\alpha_1 +k_1} \pa_y^{\beta_1} \omega\vert_{y=0} \|_{L^2(\T)} \, \prod_{l=2}^N \|  \pa_x^{\alpha_l +k_l} \pa_y^{\beta_l} \omega  \|_{L^\infty(\T \times \R_+)} $$
 Note that 
\begin{equation*}
\left\{
\begin{aligned}
& \beta_l + 2 \le j_2 \le s \quad \mbox{for all} \: l,    \\
& \alpha_1 + k_1 + \beta_1 +1 \: \le \: j, \quad   \alpha_l + k_l + \beta_l +2 \: \le \: j, \quad l=2\dots N 
\end{aligned}
\right.
\end{equation*}
($k_l \le j_1 - 1$ for $l \ge 2$, as $k_1 \ge j_1/n$). We use the  the trace theorem with the first factor:
\begin{align*}
\| \pa_x^{\alpha_1 +k_1} \pa_y^{\beta_1} \omega\vert_{y=0} \|_{L^2(\T)} \: & \lesssim \: \|   \pa_x^{\alpha_1 +k_1} \pa_y^{\beta_1} \omega \|_{L^2(\T \times \R_+)}^{1/2} \, \|   \pa_x^{\alpha_1 +k_1} \pa_y^{\beta_1+1} \omega \|_{L^2(\T \times \R_+)}^{1/2} \\ & \lesssim \:   \|  (1+y) \pa_x^{\alpha_1 +k_1} \pa_y^{\beta_1+1} \omega \|_{L^2(\T \times \R_+)}
\end{align*}
The last bound comes from \eqref{hardy1}. Using the  Sobolev imbedding \eqref{sobolev}  with  the second factor. We get:
\begin{align*}
 {\cal N}_1 \:  & \lesssim \: \sum_{(k_1,\dots,k_N) \in K_1(j_1)}   \frac{j_1!}{\prod_{l=1}^N k_l !}  \, \| \omega \|_{\dot{{\cal H}}^{\alpha_1 + k_1 + \beta_1 + 1}_\gamma} \, \prod_{l=2}^N \sum_{m=0}^2 \| \omega \|_{\dot{{\cal H}}^{\alpha_l + k_l + \beta_l + m}_\gamma}  \\
 & \lesssim \:   \sum_{(k'_1,\dots,k'_N) \in K_1(j')}    \frac{j' !}{\prod_{l=1}^N k'_l !} \, \| \omega \|_{\dot{{\cal H}}^{k'_1}_\gamma}  \, \prod_{l=2}^N \sum_{m=0}^2 \| \omega \|_{\dot{{\cal H}}^{k'_l+m}_\gamma}
 \end{align*}
The last inequality comes from the change of index 
$$k'_1 :=  k_1 +\alpha_1 + \beta_1 + 1,  \quad k'_l :=  k_l + \alpha_l + \beta_l + m, \quad l \ge 2$$
noticing that $ j' := j_1 + \sum_{l=1}^N (\alpha_l + \beta_l) + 1 \le j$. {\em By symmetry, the same bound applies to ${\cal N}_2, \: {\cal N}_N$}. Finally, we deduce from Lemma  \ref{multinom} that  
$$\| E^1_j \|_{l^1(\tau)} \:  \lesssim \: \sqrt{\dot{E}_\omega(t,\tau)} \, \sqrt{\dot{D}_\omega(t,\tau)} \: \lesssim \: \eta \dot{D}_\omega(t,\tau) \: + \: C_\eta \dot{E}_\omega(t,\tau). $$

\mspace
{\em Study of $E^2_j$}. Let $J  \in \N^* \times (2\N+1), \, |J| = j, \, 0 <  j_2 \le s$. Note that $j_2 \le s-1$ (because $s$ is even). We integrate by parts with respect to $x$: 
$$ \int_{\T \times \{0 \}} \pa_y \omega_J \, \omega_J = -\int_{\T \times \{0 \}}  \pa^{(j_1-1)}_x  \pa_y^{j_2+1} \omega \, \pa_x^{j_1+1}  \pa_y^{j_2} \omega. $$
Then, we apply the boundary reduction lemma to the second factor in the integrand. From there, the treatment is exactly the same as in the first case, 
 and leads to 
 $$ \| E^2_j \|_{l^1(\tau)} \: \lesssim \:  \eta \dot{D}_\omega(t,\tau) \: + \: C_\eta \dot{E}_\omega(t,\tau). $$
 
 \mspace
 {\em Study of $E^3_j$}. Note that   $E^3_j$ is non zero only if  $j \le s$. When $j=s$, one uses  Lemma \ref{reductionlemma}.  Otherwise,  \begin{align*}
  \: \sum_{j=0}^{s-1} \left| \int_{\T\times \{0\}} \pa_y^{j+1} \pa_y^j \omega \right| \: \le \:   \sum_{j=0}^{s-1}  \| \pa_y^{j+1} \omega \|_{L^2(\T \times \R_+)} \| \pa_y^j \omega \|^{1/2}_{L^2(\T \times \R_+)} \, \| \pa_y \pa_y^{j+1} \omega \|^{1/2}_{L^2(\T \times \R_+)} \\ 
\end{align*}
We end up with 
$$\| E^3_j \|_{l^1(\tau)}   \lesssim \: \sqrt{\dot{E}_\omega(t,\tau)} \, \sqrt{\dot{D}_\omega(t,\tau)} \: \lesssim \: \eta \dot{D}_\omega(t,\tau) \: + \: C_\eta \dot{E}_\omega(t,\tau). 
$$

\mspace
Eventually, 
\begin{equation} \label{estimEj}
\| E_j \|_{l^1(\tau)} \: \lesssim \:  \eta \dot{D}_\omega(t,\tau) \: + \: C_\eta \dot{E}_\omega(t,\tau). 
\end{equation}
 Combining this last inequality with \eqref{estimAj}-\eqref{estimBj}-\eqref{estimCj}-\eqref{estimDj}, we obtain 
\begin{align*}
 \pa_t \dot{E}_\omega(t,\tau) \: + \: \dot{D}_\omega(t,\tau) \:  
&  \lesssim \:  \left( \sqrt{\pa_\tau E_\omega(t,\tau)} + \sqrt{\dot{D}_\omega(t,\tau)} \right) \, \sqrt{\dot{E}_\omega(t,\tau)} \: + \: \eta  \dot{D}_\omega(t,\tau) + C_\eta \dot{E}_\omega(t,\tau) \\
& \lesssim \:  \eta  \dot{D}_\omega(t,\tau)  \: + \: C_\eta \, \pa_\tau E_\omega(t,\tau) 
\end{align*}
Taking $\eta$ small enough yields Proposition \ref{propvorticity}.

\subsection{Estimate of the hydrostatic energy} \label{subsechydro}
This section is devoted to the hydrostatic  energy $E_h(t,\tau)$ defined in \eqref{hydroenergy}, with its "viscous" counterpart 
\begin{equation} \label{hydrodissipation}
D_h(t,\tau) \: := \: \sum_{j=0}^{+\infty} \left( \tau^j \,    (j!)^{-7/4} \, (j+1)^{10}\right)^2 \, \|  \pa_y h_j(t,\cdot) \|_{L^2(\T \times \R_+}^2. 
\end{equation}
\begin{proposition} {\bf (Estimate on the hydrostatic energy)} \label{prophydro}

\sspace
Under the same assumptions as in Theorem \ref{theorem2},  one has for some $C > 0$ and all $t \in ]0,T]$: \begin{align} \label{Eh-prop}
 \pa_t E_h(t,\tau) \: + \: D_h(t,\tau) \: \le \: C \, \left( \pa_\tau E_\omega(t,\tau)  \: + \:  \pa_\tau E^2_g(t,\tau) \right). 
\end{align}  
\end{proposition}
The first step is to establish the equation on the function $h_j$ defined by \eqref{hj}: 
$$ h_j \: = \: \chi^a \, \frac{\pa_x^j \omega}{\sqrt{\pa_y \omega}}, \quad \chi^a(t,x,y) = \chi(y-a(t,x)). $$
Starting from the vorticity equation
$$ \pa_t \omega + u \pa_x \omega + v \pa_y \omega - \pa^2_y \omega  = 0, $$
we get 
\begin{align*}
\pa_t h_j + u \pa_x h_j + v \pa_y - \pa^2_y h_j \: & = \: -\frac{\chi^a}{\sqrt{\pa_y \omega}} [ \pa^j_x, u ] \pa_x \omega 
\:  - \:  \frac{\chi^a}{\sqrt{\pa_y \omega}}  \left( [ \pa^j_x, v ] \pa_y \omega \: - \: \pa^j_x v \pa_y \omega \right) \\
&  - \: [  \frac{\chi^a}{\sqrt{\pa_y \omega}}, \pa_t + u \pa_x  + v \pa_y - \pa^2_y] \pa^j_x \omega  \: - \:   \frac{\chi^a}{\sqrt{\pa_y \omega}}  \pa^j_x v \pa_y \omega.
\end{align*}
Note that we have singled out the term $\pa^j_x v \pa_y \omega = - \int^0_y \pa^{j+1}_x u \pa_y \omega$  in the commutator with $v \pa_y$. This is the hardest term to control, as it involves $j+1$ derivatives with respect to $x$. We shall use a cancellation property similar to \eqref{cancelhj}. However, such cancellation will not be enough, and we shall rely on the extra regularity 
offered by  the monotonicity energy $E^2_g$. 

\mspace
Performing a standard energy estimate on the previous equation, multypling by $ \left(\tau^{j} (j!)^{-7/4} \, (j+1)^{10}\right)^2$,  summing over $j$, we end up with 
\begin{align} \label{Ehd} 
 \frac{1}{2} \pa_t E_h(t,\tau) \: + \: D_h(t,\tau) \: \le \: \left( \| A_j \|_{l^2(\tau)} + \| B_j \|_{l^2(\tau)} +   \| C_j \|_{l^2(\tau)} \right) \sqrt{\pa_\tau E_h(t,\tau)} \: +\:  \| D_j \|_{l^1(\tau)} \: + \:   \| E_j \|_{l^1(\tau)} 
 \end{align} 
where 
\begin{align*} 
 A_j \: & := \: \frac{1}{j^{1/2}}  \| \frac{\chi^a}{\sqrt{\pa_y \omega}} [ \pa^j_x, u ] \pa_x \omega\|_{L^2}, \quad B_j \: := \:  \frac{1}{j^{1/2}}  \|  \frac{\chi^a}{\sqrt{\pa_y \omega}}  \left( [ \pa^j_x, v ] \pa_y \omega \: - \: \pa^j_x v \pa_y \omega \right) \|_{L^2}, \\
 C_j \: & := \:  \frac{1}{j^{1/2}}  \| (\pa_t + u \pa_x + v \pa_y - \pa^2_y)  \bigl(  \frac{\chi^a}{\sqrt{\pa_y \omega}} \bigr) \, \pa^j_x \omega  \|_{L^2}, \quad  D_j \:  := \:  2 \int_{\T \times \R_+} \pa_y \frac{\chi^a}{\sqrt{\pa_y \omega}} \, \pa_y \pa_x^j  \omega  \; h_j   \\
 E_j \: & := \: \int_{\T \times \R_+} \frac{\chi^a}{\sqrt{\pa_y \omega}}  \pa^j_x v \, \pa_y \omega \; h_j. 
 \end{align*}
Note that $h_j$ is compactly supported in $y$, so that no boundary integral is present on  the right-hand side. 
 
 \mspace
 \subsubsection*{Estimate on $A_j$}
 We proceed here as for the term $A_j$  of subsection \ref{subsecvorticity}. The treatment is actually simpler, as only $x$  derivatives of $u$ and $\omega$ are  involved in the expression of $A_j$. By  use of  the Sobolev inequality \eqref{sobolev}, one never encounters more than one $y$-derivative, in particular never more than $s$.  For brevity, we skip the details. We find 
 \begin{equation} \label{Ajhydro}
  \| A_j\|_{l^2(\tau)} \: \lesssim \:   \sqrt{\pa_\tau E_\omega(t,\tau)}
 \end{equation}
 
   \subsubsection*{Estimate on $B_j$}
   Again, the treatment is parallel to the one of subsection \ref{subsecvorticity}. Denoting by $K^a$ the support of $\chi^a$, we find that 
\begin{align*}
B_j^1 \: & \lesssim \: \frac{1}{j^{1/2}} \sum_{k=1}^{j-1} \binom{j}{k} \| \pa^k_x v \, \pa_y \pa_x^{j-k} \omega \|_{L^2(\T \times K^a)} \\
\: & \lesssim \: \frac{1}{j^{1/2}} \sum_{k=1}^{\frac{j}{2}} \binom{j}{k}  \| \pa^k_x v \|_{L^\infty(\T \times K^a)} \, 
 \| \pa_y \pa_x^{j-k} \omega \|_{L^2(\T \times K^a)} \\
 \:  & + \:  \frac{1}{j^{1/2}}   \sum_{k=\frac{j}{2}}^{j-1}  \binom{j}{k}  \| \pa^k_x v \|_{L^2(\T \times K^a)} \,  \| \pa_y \pa_x^{j-k} \omega \|_{L^\infty(\T \times K^a)} \: := \: B_j^1  \: + \: B_j^2
 \end{align*}
 For the first term, we set $k' := k-1$ and use \eqref{v_infty} to get (we drop the prime)
 \begin{align*}
 B^1_j \: & \lesssim \:   \frac{1}{j^{1/2}} \sum_{k=0}^{\frac{j}{2}-1}   \binom{j}{k+1} \sum_{m=2}^3 \| \omega \|_{{\cal H}^{k+m}_\gamma} \, \| \omega \|_{{\cal H}^{j-k}_\gamma} \\
 & \lesssim \: \sum_{k=0}^{\frac{j}{2}}  \sum_{m=2}^3 \binom{j}{k} \| \omega \|_{{\cal H}^{k+m}_\gamma} \left( (j-k)^{1/2} \| \omega \|_{{\cal H}^{j-k}_\gamma} \right)
 \end{align*}
 The second term $B_j^2$ is handled in a symmetric way, thanks to the fact that  index $k$ stops at $j-1$. 
 Applying Lemma \ref{binom}, we eventually derive the bound: 
 \begin{equation} \label{Bjhydro}
 \| B_j \|_{l^2(\tau)} \: \lesssim \: \sqrt{E_\omega(t,\tau)} \, \sqrt{\pa_\tau E_\omega(t,\tau)} \: \lesssim  \sqrt{\pa_\tau E_\omega(t,\tau)}
 \end{equation}
 
  \subsubsection*{Estimates on $C_j, D_j$}
  Clearly, 
\begin{equation} \label{Cjhydro}
  C_j \: \lesssim \: \frac{1}{j^{1/2}} \| \pa_x^j \omega \|_{L^2(\T \times K^a)}, \quad \| C_j \|_{l^2(\tau)} \: \lesssim \: \sqrt{E_\omega(t,\tau)} 
\end{equation}  
 and after integration by parts
 \begin{equation} \label{Djhydro}
 \begin{aligned}
&D_j   \:  \lesssim \:  \frac{1}{j^{1/2}} \| \pa^j_x \omega\|_{L^2(\T \times K^a)} \left( \| h_j \|_{L^2(\T \times \R_+)} + \| \pa_y h_j \|_{L^2(\T \times \R_+)} \right), \\
 &\| D_j \|_{l^1(\tau)} \:  \lesssim \: \sqrt{E_\omega(t,\tau)} \left(  
 \sqrt{E_h(t,\tau)} \: + \: \sqrt{D_h(t,\tau)} \right)  \\
 &  \lesssim \: \eta D_h(t,\tau) \: + \: C_\eta E_\omega(t,\tau). 
 \end{aligned}
 \end{equation}
 Let us stress that Lemma \ref{relations} allowed us  to control $E_h$ by $E_\omega$.    
 \subsubsection*{Estimates on $E_j$}
It remains to handle the bad term: 
\begin{align*}
E_j & \: =  \: \int_{\T \times \R_+} (\chi^a)^2 \pa_x^j v \pa_x^j \omega = -\int_{\T \times \R_+} \pa_y (\chi^a)^2 \pa_x^j v \pa_x^j  u - \int_{\T \times \R_+} (\chi^a)^2 \pa_x^j \pa_y v \pa^j_x u \\
& \: =  \: \int_{\T \times \R_+} \pa_y (\chi^a)^2 \pa_x \left( \int_0^y \pa_x^j u \right)  \pa_x^j  u \: + \:   \int_{\T \times \R_+}  (\chi^a)^2 \pa_x \frac{(\pa_x^j u)^2}{2} \\
& \: =  \: -  \int_{\T \times \R_+} \pa_y (\chi^a)^2   \int_0^y\!\!  \pa_x^j u \,\,     \pa_x^{j+1}  u \: - 
\:  \int_{\T \times \R_+} \left( \pa_{xy} (\chi^a)^2  \int_0^y \!\! \pa^j_x u \,\,  \pa^j_x u \: + \: \pa_x (\chi^a)^2  \frac{(\pa_x^j u)^2}{2} \right) \\
& \: := \: E^1_j \: + \: E^2_j
\end{align*} 
 Clearly, 
\begin{equation} \label{boundE2j} 
\| E^2_j \|_{l^1(\tau)} \: \lesssim \:  \Big \| \, \| \pa_x^j u \|_{L^2(\T \times K^a)}  \Big \|_{l^2(\tau)}^2   \: \lesssim \: \Big \| \, \| \pa_x^j \omega \|_{L^2(\T \times K^a)} \Big  \|_{l^2(\tau)}^2 \: \lesssim \: E_\omega(t,\tau) 
\end{equation}
  with  Poincar\'e inequality allowing to go from $u$ to $\omega$. 
  
  \mspace
  We are left with the treatment of $E^1_j$. {\em A crucial remark is that for  $\eps > 0$ small enough, for all $t \in [0,T]$,  the integrand in $E^1_j = E^1_j(t)$ is supported in 
  $$\{ |y -a| \ge \eps  \} :=  \{ (x,y), \: |y-a(t,x)| \ge \eps \} . $$}
 Indeed,  $\chi^a$ is one in a neighborhood of $y=a$ (moreover, $\chi^a = 0$ for $y \ge 3$). We then recall the decomposition stated in   Lemma \ref{decompo}: for all $j$,
\begin{equation} \label{decompo2}
 \pa^j_x u \: = \:  u^j_g \: + \: C_j \, 1_{\{y > a\}} \, \omega \: 
 \end{equation}
  with 
 \begin{align*} 
 u^j_g  \:  & :=  \:  \omega(t,x,y)  \int_0^y  \left(\psi + \frac{(1-\psi)}{\omega}\right)^{-1} \frac{g_j}{\omega^2}  \: \mbox{ in } \: \{ y < a\}, \\ 
  u^j_g  \: & :=  \:  \omega(t,x,y)  \int_3^y  \left(\psi + \frac{(1-\psi)}{\omega}\right)^{-1} \frac{g_j}{\omega^2}  \: \mbox{ in } \: \{ y > a\}. 
  \end{align*}
Note that 
\begin{equation} \label{ujg}
\| u^j_g \|_{L^2(\{ \eps \le |y- a|, \:  y \le  3  \})} \: \lesssim \: \| g_j \|_{L^2(\{y \le 3\})}.
 \end{equation} 
For $y < a$, we write 
 \begin{equation} \label{E1j-}
 \begin{aligned}
  \left| \int_{\{y < a\}} \pa_y(\chi^a)^2 \left( \int_0^y \pa_x^j u \right) \pa_x^{j+1} u \right| \: & = \: \left| \int_{\{y < a\}} \pa_y(\chi^a)^2 \left( \int_0^y u^j_g \right) u^{j+1}_g \right| \\
&  \lesssim \: \| g_j \|_{L^2(\{y \le 3\})}   \, \| g_{j+1} \|_{L^2(  \{y \le 3\})} \\ 
&   \lesssim \: \| g_j \|_{L^2(\{y \le 3\})}   \, \| \pa^{j+1}_x \omega \|_{L^2(  \{y \le 3\})}. 
\end{aligned}
 \end{equation}
 
 \mspace
 For $y > a$, we  rather write
 \begin{equation} 
 \begin{aligned}
&  \left| \int_{\{y > a\}} \pa_y(\chi^a)^2 \left( \int_0^y \pa_x^j u \right) \pa_x^{j+1} u \right| \: = \: \left| \int_{\{y > a\}} \pa_y(\chi^a)^2 \left( \int_0^y u^j_g \right) \pa^{j+1}_x u  +  \int_{\{y > a\}}     C_j \left( \int_a^y \omega  \right) \pa^{j+1}_x u   \right| \\   
&  \lesssim \: \| \int_0^y u^j_g   \|_{L^2(\{a + \eps \le  y \le 3\})}   \, \| \pa^{j+1}_x \omega \|_{L^2(\T \times \R_+)}  \:   +  \:   \left| \int_{\{y > a\}}     C_j \left( \int_a^y \omega \right) \pa^{j+1}_x u   \right| .   \\
\end{aligned}
 \end{equation} 
   Note that although $y$ is away from $a$,  the function $ \int_0^y u^j_g$ involves values of $u^j_g$ near the critical curve $y=a$. As the quantity $\frac{g_j}{\omega^2}$ (involved in the definition of $u^j_g$) degenerates at $y=a$, a control like \eqref{ujg} is not obvious. Still, we claim:
\begin{lemma} \label{lem-ujg}
\begin{equation} \label{intujg}
 \| \int_0^y \!\!u^j_g \;  \|_{L^2(\{a + \eps \le  y \le 3 \})} \: \lesssim \: \| g_j \|_{L^2(\{ y \le  3\})}  \: \lesssim \: \| \pa^j_x \omega \|_{L^2(\T \times \R_+)} . 
 \end{equation} 
Furthermore, we have 
\begin{equation} \label{estimCj1}
\| C_j \|_{L^2(\T)} \: \lesssim \:  \:  \| \pa^j_x \omega \|_{L^2(\T \times \R_+)}, \quad \| \pa_x C_j \|_{L^2(\T)} \: \lesssim \:  \| \pa^{j+1}_x \omega \|_{L^2(\T \times \R_+)} \: + \:  \| \pa^j_x \omega \|_{L^2(\T \times \R_+)}.
\end{equation}
\end{lemma}
Finally, using above bounds and the  decomposition  
$$ \pa^{j+1}_x u =  \pa_x u^j_g  + \pa_x C_j \omega \: + \: C_j \pa_x \omega$$
 we end up with 
 \begin{align*}
& \left|  \int_{\{y > a\}} \pa_y(\chi^a)^2 \left( \int_0^y \pa_x^j u \right) \pa_x^{j+1} u \right| \\
&  \lesssim \:   \| g_j \|_{L^2(\{ y \le  3\})} \| \pa^{j+1}_x \omega \|_{L^2(\T \times \R_+)}   \: + \:  \left|  \int_{\T \times \R_+}  \pa_y(\chi^a)^2  C_j \left( \int_a^y \omega \right)  \pa_x u^j_g \right| \\ 
& + \:  \left|  \int_{\T \times \R_+}  \pa_y(\chi^a)^2 \pa_x \frac{C_j^2}{2}\,  \omega \left( \int_0^y \omega \right) \right|  + \:  \| \pa^j_x \omega \|_{L^2(\T \times \R_+)}^2   \\ 
& \lesssim \:   \| g_j \|_{L^2(\{ y \le  3\})} \| \pa^{j+1}_x \omega \|_{L^2(\T \times \R_+)}  \: + \:  \| \pa^j_x \omega \|_{L^2(\T \times \R_+)}^2 
 \end{align*}
 after  integration  by parts of  the integral terms at the right-hand side. The structure of the second integral term is crucial. Indeed, the  term  involving 
the product $ (\pa_x C_j \omega) ( C_j \int_a^y \omega )  $  can be integrated by parts. 
Note that  a more general  bilinear term in $C_j$ and $\pa_x C_j$ would  generically involve  an upper bound like $\| C_j \|_{L^2} \, \| \pa_x C_j \|_{L^2}$. Such bound would  ruin our strategy,  based on anisotropic energy. 
 
 \mspace
 From the previous bound, from inequality  \eqref{E1j-} and Cauchy-Schwarz  in $l^2(\tau)$, we get 
 \begin{equation*}
 \| E^1_j \|_{l^1(\tau)}  \: \lesssim \:  \Big\| \frac{1}{j^{5/4}}  \| \pa^{j+1}_x \omega  \|_{L^2(\T \times \R_+)}  
 \,  \Big\|_{l^2(\tau)} \;    \Big \| \,  j^{5/4} \,\|  g_j  \|_{L^2(\{ y \le 3 \})} \Big \|_{l^2(\tau)} \: + \:  \Big\| \, \| \pa^j_x \omega \|_{L^2(\T \times \R_+)}  
  \Big\|_{l^2(\tau)}^2 . 
 \end{equation*}
  Here, we remark that 
 \begin{equation} \label{paj+1omega}
 \begin{aligned}
   \Big\| \frac{1}{j^{5/4}} \| \pa^{j+1}_x \omega  \|_{L^2(\T \times \R_+)}  
 \, \Big \|_{l^2(\tau)}^2 \: & = \: \sum_{j=0}^{+\infty}  \frac{1}{j^{5/2}} \left( \tau^j (j!)^{-7/4} (j+1)^{10} \right)^2 \| \omega \|_{{\cal H}^{j+1}_\gamma}^2 \\
& \lesssim \:  \sum_{j=1}^{+\infty} j \left( \tau^j (j!)^{-7/4} (j+1)^{10} \right)^2 \| \omega \|_{{\cal H}^{j}_\gamma}^2 \\
&   \lesssim \: \| \sqrt{j}  \| \omega \|_{{\cal H}^{j}_\gamma} \|_{l^2(\tau)}^2 \: \lesssim \: \pa_\tau E_\omega(t,\tau). 
\end{aligned}
\end{equation}
 To  conclude, it remains to evaluate the $l^2(\tau)$ norm of $  j^{5/4} \,\|  g_j  \|_{L^2(\{ y \le 3\})}$. 
 We use
 the following bound, to be proved in appendix   
\begin{lemma} \label{lemmatildegjgj}
For all $\alpha \ge 0$, 
\begin{equation}
\| \, \| j^{\alpha/4} g_j \|_{L^2(\{ y \le3\})}   \|_{l^2(\tau)}  \: \lesssim \:  \| \, \| j^{\alpha/4} \tilde{g}_j \|_{L^2(\{ y \le 3\})}    \|_{l^2(\tau)} \: + \: \|  j^{(\alpha-3)/4} \| \omega \|_{{\cal H}^j_\gamma} \|_{l^2(\tau)} . 
\end{equation}  
\end{lemma}

\mspace 
For $\alpha=5$, we find 
\begin{equation}
\| \, \| j^{5/4} g_j \|_{L^2(\{ y \le 3\})}   \|_{l^2(\tau)}  \: \lesssim \:  \sqrt{\pa_\tau E^2_g(t,\tau)} \: + \: \sqrt{\pa_\tau E_\omega(t,\tau)}. 
\end{equation}  
We inject this bound  and bound \eqref{paj+1omega} in the estimate of $\| E^1_j \|_{l^1(\tau)}$, 
  we get 
\begin{equation} \label{boundE1j}
\| E^1_j \|_{l^1(\tau)}  \: \lesssim \:  \sqrt{\pa_\tau E_\omega(t,\tau)} \sqrt{\pa_\tau E^2_g(t,\tau)} \, \: + \: \pa_\tau E_\omega(t,\tau) \: \lesssim \: \pa_\tau E_\omega(t,\tau) \: + \:   \pa_\tau E^2_g(t,\tau)
   \end{equation}
 Combination of this inequality and inequality \eqref{boundE2j} yields 
 \begin{equation} \label{Ejhydro}
 \| E_j \|_{l^1(\tau)}  \: \lesssim \: \pa_\tau E_\omega(t,\tau) \: + \:   \pa_\tau E^2_g(t,\tau)
\end{equation}
Collecting \eqref{Ajhydro}-\eqref{Bjhydro}-\eqref{Cjhydro}-\eqref{Djhydro}-\eqref{Ejhydro} leads to 
 \begin{equation*}
  \pa_t E_h(t,\tau) \: + \: D_h(t,\tau) \:  \lesssim \:  \eta D_h(t,\tau) \: + \: C_\eta \, \pa_\tau E_\omega(t,\tau) \: + \:   \pa_\tau E^2_g(t,\tau). 
 \end{equation*}
 Taking $\eta$ small enough yields Proposition \ref{prophydro}.

\subsection{Estimate of the (second) monotonicity energy} \label{subsecmono}
The control of $E_h$, performed in the previous subsection, has relied on  the so-called monotonicity energy $E^2_g$. Its time variations $\pa_t E^2_g$ need in turn to be controlled by $\pa_\tau {\cal E}$, in order to obtain a closed estimate.  The main difficulty comes from the extra factor $j^{3/2}$ in the definition of $E^2_g$: indeed, naive energy bounds would  involve the $l^2(\tau)$ norm of $\left( j^{5/2} \,  \|  \omega \|_{{\cal H}^j_\gamma}\right)_{j \in \N}$, which is not controlled by $\sqrt{\pa_\tau {\cal E}}$. 

\mspace
To get a good estimate, we must take advantage of  cancellation properties in the equation for 
$$  \tilde g_j(t,x,y) \: = \:  \pa_x^{j-5}\left(  \omega \pa_x^5 \omega  -   \pa_y\omega \pa_x^5 u \right).  $$
As mentioned in section \ref{statements}, this kind of cancellations was used for the first time in paper \cite{MW13prep}, in the case of monotonic data: this paper shows indeed some Sobolev stability, with special role played by   
\begin{equation*}  
g_j \: := \: \omega \pa_x^j \omega - \pa_y \omega \pa_x^j u.  
\end{equation*}
However, in our setting, there is a technical difficulty with using $g_j$ (or the  $g_j$ defined in \eqref{gj}, better suited to large $y$). Broadly speaking, 
the equation for $g_j$ involves the commutator term 
$ \: j \, \pa_x^{j-1} v  \, \pa_x \pa_y \omega$, with 
  a  factor $j$ in front. This factor is harmless at the Sobolev level (finite $j$), but  is annoying at the Gevrey level, for which  behaviour at large $j$ is important. This is why we rather use $\tilde g_j$ than $g_j$: ${\cal C}_j$ is somehow replaced by $5 \, \pa_x^{j-1}  v \ \pa_x \pa_y \omega$, with no bad factor. Roughly, instead of differentiating $j$ times the velocity and vorticity equations, and combining them to obtain cancellations, we go the reverse way: we first combine the equations and then differentiate, so as to benefit from "earlier" cancellations. 

\mspace
We shall prove  
\begin{proposition} \label{propmono2}
Under the assumptions of Theorem \ref{theorem2}, we have 
\begin{equation} \label{Eg-prop}
\pa_t E^2_g(t,\tau) \: + \: D^2_g(t,\tau) \: \le \: C \, \left( \pa_\tau E_\omega(t,\tau) \: + \:  \pa_\tau E_g^2(t,\tau) \: + \: D_h(t,\tau) \right), \quad C > 0.
\end{equation}
\end{proposition}
We recall  that $D_h$ was defined  in   \eqref{hydrodissipation}, whereas 
\begin{equation} \label{mono2dissipation}
D^2_g(t,\tau) \: := \:  \sum_{j\in \N} \left( \tau^j \,    (j!)^{-7/4} \, (j+1)^{10}\right)^2  \: j^{3/2}  \|\pa_y \tilde g_j(t,\cdot) \|_{L^2(\T \times \R_+)}^2
\end{equation}

\mspace
As usual, the starting point is to derive an equation for  $\tilde g_j$,  $\: j \ge 5$.  First, we apply $\omega \pa^x_5$ to the vorticity equation, $\pa_y \omega \pa_x^5$ to the velocity equation, and subtract  one from the other: we obtain  
\begin{multline} 
\pa_t \tilde g_5 + u \pa_x \tilde g_5 + v \pa_y \tilde g_5 - \pa^2_y \tilde g_5 = - \sum_{k=1}^5 \binom{5}{k} \pa^k_x u \, \bar{g}_{5-k+1} \: - \: \sum_{k=1}^4\binom{5}{k}  \pa_x^k v \,  \hat{g}_{5-k+1} \\
- (\pa_t + u \pa_x + v \pa_y - \pa^2_y ) \pa_y \omega \, \pa^5_x u + 2 \pa^2_y \omega \, \pa^5_x \omega - 2 \pa_y \omega \, \pa_y \pa^5_x \omega \:  := \: \sum_{i=1}^5 {\cal C}_i 
\end{multline}
with 
\begin{equation} \label{ggg}
\begin{aligned}
\bar{g}_k \: & := \: \omega \pa^k_x \omega  -  \pa_y \omega \pa^k_x u, \\
\hat{g}_k \: & := \: \omega \pa^{k-1}_x \pa_y \omega   - \pa_y \omega \pa_x^{k-1}  
\omega, \quad 1 \le k.
\end{aligned}
\end{equation}
We finally apply $\pa_x^{j-5}$ to the equation, which yields
\begin{equation}
\pa_t \tilde g_j + u \pa_x \tilde g_j + v \pa_y \tilde g_j - \pa^2_y \tilde g_j = [u \pa_x, \pa_x^{j-5}] \tilde g_5 + [v\pa_y, \pa_x^{j-5}] \tilde g_5 \: + \: \sum_{i=1}^5 \pa^{j-5}_x {\cal C}_i 
\end{equation}
We perform an energy estimate, multiply by $j^{3/2} \left(\tau^j (j!)^{-7/4} (j+1)^{10} \right)^2$ and sum over $j$:
\begin{equation}
\pa_t E^2_g(t,\tau) \: + \: D^2_g(t,\tau)\: \le \: \left( \| A_j \|_{l^2(\tau)} +   \| B_j \|_{l^2(\tau)} + \sum_{i=1}^4  \| C_{i,j} \|_{l^2(\tau)} \right) \sqrt{\pa_\tau E^2_g(t,\tau)} \: + \: \| \max(D_j, 0) \|_{l^1(\tau)} 
\end{equation}
with 
\begin{align*}
 A_j & := j^{1/4} \|  [u \pa_x, \pa_x^{j-5}] \tilde g_5  \|_{L^2(\T \times \R_+)}, \quad B_j := j^{1/4} \|  [v\pa_y, \pa_x^{j-5}] \tilde g_5 \|_{L^2(\T \times \R_+)}, \\
 C_{i,j} & := j^{1/4} \| \pa_x^{j-5} {\cal C}_i \|_{L^2(\T \times \R_+)}, \quad D_j \: := \: j^{3/2} \int_{\T \times \R_+} \pa_x^{j-5} {\cal C}_5 \tilde g_j.
 \end{align*}
Note that 
$$\pa_y \tilde g_j \vert_{y=0} = \pa^{j-5}_x \left( \pa_y \omega \pa^5_x \omega  - \pa^2_y \omega \pa^5_x u \right)\vert_{y=0} = 0.$$
 Indeed, $\pa^5_x u\vert_{y=0} = 0$, and $\pa_y \omega\vert_{y=0} = \pa^2_y u\vert_{y=0} = 0$, evaluating the Prandtl equation at $y=0$. In particular, there is no boundary term due to the integration by parts. 
 
 \subsubsection*{Estimate on $A_j$}
 We find 
 $$   [u \pa_x, \pa_x^{j-5}] \tilde g_5 \: = \: \sum_{k=1}^{j-5} \binom{j-5}{k} \pa_x^k u \, \tilde g_{j-k+1}.  $$
Thus,
 \begin{equation*}
 A_j \:  \le \: j^{1/4} \sum_{k=1}^{j/2}   \binom{j-5}{k} \| \pa_x^k u \|_{L^\infty} \, \| \tilde g_{j-k+1} \|_{L^2} \: + \: j^{1/4} \sum_{k=j/2}^{j-5}  \| \pa_x^k u \|_{L^2} \, \| \tilde g_{j-k+1} \|_{L^\infty} \: := \: A^1_j + A^2_j. 
 \end{equation*}
 Using like before \eqref{hardy2} and \eqref{sobolev}, we get
 $$ \| \pa_x^k u \|_{L^\infty} \: \lesssim \: \sum_{m=0}^1 \| \omega \|_{{\cal H}^{k+m}_\gamma}, \quad 
  \| \tilde g_{j-k+1} \|_{L^\infty} \: \lesssim \:  \sum_{m=1}^2 \| \tilde g_{j-k+m} \|_{L^2} \: + \:  \sum_{m=2}^3 \| \pa_y \tilde g_{j-k-1+m} \|_{L^2}. $$ 
 Hence (set $k' := k-1$, drop the prime) 
\begin{equation*}
A^1_j \:  \lesssim \: j^{1/4} \sum_{k=0}^{j/2-1} \binom{j-5}{k+1} \sum_{m=1}^2 \| \omega \|_{{\cal H}^{k+m}_\gamma}  \, \| \tilde g_{j-k} \|_{L^2} \: 
 \lesssim \:   \sum_{k=0}^{j/2} \binom{j}{k} \sum_{m=1}^2 \| \omega \|_{{\cal H}^{k+m}_\gamma}  
(j-k)^{1/4} \,  \| \tilde g_{j-k} \|_{L^2} 
\end{equation*}
using that for $0 \le k \le j/2$, one has $ \binom{j-5}{k+1} \, \lesssim \,  \, \binom{j}{k}$ and $j^{1/4} \, \lesssim \, (j-k)^{1/4}$. Hence, by Lemma \ref{binom} 
\begin{equation} \label{A1jmono2}
\begin{aligned}
\| A^1_j \|_{l^2(\tau)} \: & \lesssim \: \| \, \| \omega \|_{{\cal H}^j_\gamma} \|_{l^2(\tau)} \, \| j^{1/4} \, \| \tilde g_j \|_{L^2(\T \times \R_+)} \|_{l^2(\tau)}\\
& \lesssim \: \sqrt{E_\omega(t,\tau)} \sqrt{\pa_\tau E^2_g(t,\tau)} \: \lesssim \:  \sqrt{\pa_\tau E^2_g(t,\tau)} 
\end{aligned}
\end{equation}

\mspace
Finally, 
\begin{align*} 
A^2_j  \: &  \lesssim \: j^{1/4} \sum_{k=j/2}^{j-5} \binom{j-5}{k}  \| \pa_x^k u \|_{L^2}  \, \sum_{m=1}^2 \| \tilde g_{j-k+m} \|_{L^2} \: + \:  \sum_{m=2}^3 \| \pa_y \tilde g_{j-k-1+m} \|_{L^2} \\
 &  \lesssim \:  \sum_{k=j/2}^{j} \binom{j}{k} k^{1/4} \| \omega \|_{{\cal H}^k_\gamma} \left( \sum_{m=1}^2 \| \tilde g_{j-k+m} \|_{L^2} \: + \:  \sum_{m=2}^3 \| \pa_y \tilde g_{j-k-1+m} \|_{L^2} \right) 
 \end{align*}
 using that for $j/2 \le k \le j-5$, one has   $ \binom{j-5}{k} \, \lesssim \, \binom{j}{k}$ and $j^{1/4} \lesssim k^{1/4} $. Hence, 
\begin{equation*}
\| A^2_j \|_{l^2(\tau)} \: \le \: \left( \| \: \| \tilde g_j  \|_{L^2(\T \times \R_+)} \|_{l^2(\tau)}  \: + \:   \| \: \| \pa_y \tilde g_{j-1} \|_{L^2(\T \times \R_+)}\|_{l^2(\tau)} \right) \,   \Big\| j^{1/4} \, \| \omega \|_{{\cal H}^j_\gamma} \Big\|_{l^2(\tau)} . 
\end{equation*}
We conclude by using 
\begin{lemma} \label{lemmatildegj}
For all $\alpha \ge 0$, for all $l=0,1,2$,  
$$ \| j^\alpha \| \pa_y^l \tilde g_{j-l} \|_{L^2(\T \times \R_+)}  \|_{l^2(\tau)}  \: \lesssim \: \| j^\alpha \| \omega \|_{{\cal H}^j_\gamma} \|_{l^2(\tau)}, \quad \alpha \ge 0. $$
\end{lemma}
For a proof, see the appendix.   We get 
\begin{equation} \label{A2jmono2}
 \| A^2_j \|_{l^2(\tau)} \: \lesssim \: \sqrt{E_\omega(t,\tau)} \, \sqrt{\pa_\tau E_\omega(t,\tau)} \: \lesssim \: \sqrt{\pa_\tau E_\omega(t,\tau)} 
 \end{equation}
 and combining with \eqref{A1jmono2}, 
 \begin{equation}
  \| A_j \|_{l^2(\tau)} \: \lesssim \: \sqrt{\pa_\tau E^2_g(t,\tau)} \: + \: \sqrt{\pa_\tau E_\omega(t,\tau)} . 
  \end{equation}

 \subsubsection*{Estimate on $B_j$}
Leibniz formula gives
$$  [v \pa_y , \pa_x^{j-5} ] \tilde g_5 \: = \:  \sum_{k=1}^{j-5} \binom{j-5}{k} \pa^k_x v \, \pa_y \tilde g_{j-k}.  $$
Proceeding as before, we get easily 
\begin{align*}
 B_j \: & \le \: j^{1/4} \sum_{k=1}^{j/2}   \binom{j-5}{k} \sum_{m=1}^2 \| \omega \|_{{\cal H}^{k+m}_\gamma} \, \| \pa_y \tilde g_{j-k} \|_{L^2} \\
 & + \: j^{1/4} \sum_{k=j/2}^{j-5}  \| \omega \|_{{\cal H}^{k+1}_\gamma} \, \left( \sum_{m=1}^2 \| \pa_y \tilde g_{j-k-1+m} \|_{L^2} + \sum_{m=2}^3 \| \pa^2_y \tilde g_{j-k-2+m} \|_{L^2} \right) \: := \: B^1_j + B^2_j
 \end{align*}
 Then, using that $\binom{j-5}{k} \lesssim \binom{j}{k}$, $\: j^{1/4} \lesssim (j-k)^{1/4}$ for $0 \le k \le j/2$, we find 
 $$ B^1_j \: \le \: \sum_{k=0}^{j/2} \binom{j}{k} \sum_{m=1}^2 \| \omega  \|_{{\cal H}^{k+m}_\gamma}  (j-k)^{1/4} \| \pa_y \tilde g_{j-k} \|_{L^2} $$
 so that 
 \begin{equation} \label{B1jmono2}
 \| B^1_j \|_{l^2(\tau)} \: \lesssim \: \sqrt{E_\omega(t,\tau)} \, \sqrt{\| j^{1/4} \| \pa_y \tilde g_j \|_{L^2(\T \times \R_+)} \|_{l^2(\tau)}} \: \lesssim \: \sqrt{D^2_g(t,\tau)} 
 \end{equation}

\mspace
Also, 
\begin{align*}
B^2_j \: & \lesssim \: j^{1/4} \sum_{k=j/2+1}^{j-4} \binom{j-5}{k-1} \| \omega \|_{{\cal H}^k_\gamma}  \left( \sum_{m=2}^3 \| \pa_y \tilde g_{j-k-1+m} \|_{L^2} + \sum_{m=3}^4 \| \pa^2_y \tilde g_{j-k-2+m} \|_{L^2} \right) \\
& \lesssim \: \sum_{k=j/2}^{j} \binom{j}{k} \left( k^{1/4} \| \omega \|_{{\cal H}^k_\gamma} \right)  \left( \sum_{m=2}^3 \| \pa_y \tilde g_{j-k-1+m} \|_{L^2} + \sum_{m=3}^4 \| \pa^2_y \tilde g_{j-k-2+m} \|_{L^2} \right) 
\end{align*}
using that $\binom{j-5}{k-1} \lesssim \binom{j}{k}$, $\: j^{1/4} \lesssim k^{1/4}$ in the range $j/2 \le k \le  j-4$. Hence,
\begin{equation} \label{B2jmono2}
\begin{aligned}
\| B^2_j \|_{l^2(\tau)} \: & \lesssim \: \sqrt{\pa_\tau E_\omega(t,\tau)} \left( \| \, \| \pa_y \tilde g_{j-1} \|_{L^2(\T \times \R_+)} \|_{l^2(\tau)} \: + \:   \| \, 
\| \pa^2_y \tilde g_{j-2}  \|_{L^2(\T \times \R_+)} \,   \|_{l^2(\tau)} \right) \\
&  \lesssim \:   \sqrt{\pa_\tau E_\omega(t,\tau)}    \sqrt{E_\omega(t,\tau)} \: \lesssim  \sqrt{\pa_\tau E_\omega(t,\tau)}   
\end{aligned}
\end{equation}
where we have applied Lemma \ref{lemmatildegj} to go from the first to the second line. Gathering the estimates on $B^1_j$ and $B^2_j$, 
\begin{equation} \label{Bjmono2}
\| B_j \|_{l^2(\tau)} \: \lesssim \: \sqrt{D^2_g(t,\tau)} + \sqrt{\pa_\tau E_\omega(t,\tau)} 
\end{equation}

\subsubsection*{Estimates on $C_{i,j}$, $1 \le i \le 2$}
The quantities $C_{1,j}$ and $C_{2,j}$ involve commutators with the nonlinearities. As we have manipulated several times such commutators, we shall not detail the computation of their $l^2(\tau)$ norms. One can check that 
\begin{equation} \label{C1j}
\begin{aligned}
\| C_{1,j} \|_{l^2(\tau)} \: & \lesssim \: \| \, \| \omega \|_{{\cal H}^j_\gamma} \|_{l^2(\tau)} \, \sum_{m=1}^5 \| j^{1/4} \, \| \pa^{j-5}_x \bar{g}_m \|_{L^2(\T \times \R_+)} \|_{l^2(\tau)} \\
&  + \:  \|  j^{1/4} \, \| \omega \|_{{\cal H}^j_\gamma} \|_{l^2(\tau)}  \sum_{m=1}^5 \left( \| \: \| \pa^{j-5}_x \bar{g}_m \: \|_{L^2(\T \times \R_+)} \|_{l^2(\tau)} +  \| \: \|\pa_y \pa^{j-6}_x \bar{g}_m  \|_{L^2(\T \times \R_+)} \|_{l^2(\tau)} \right) \\
&  \lesssim \: \sqrt{\pa_\tau E_\omega(t,\tau)} 
\end{aligned}
\end{equation}
Here, we  use  implicitly the bound 
\begin{equation} \label{barg1}
 \| \, j^\alpha \| \pa_y^l \pa^{j-5-l}_x  \bar{g}_{m} \|_{L^2(\T \times \R_+)} \|_{l^2(\tau)}  \: \lesssim \: \| \: j^{\alpha} \, \| \omega \|_{{\cal H}^j_\gamma} \:
  \|_{l^2(\tau)}   
\end{equation}
which is valid for all $1 \le m \le 5$,  $\alpha \ge 0$ and $0 \le l \le 2$. Note that the case $m=5$ corresponds to  Lemma \ref{lemmatildegj}. This inequality  holds {\it a fortiori} for $m \le 4$, as  less $x$-derivatives are involved. Actually, one can even show that
\begin{equation} \label{barg2}
 \mbox{for } \: 1\le m \le 4, \quad \| \: j^\alpha \| \pa_y^l  \pa^{j-4-l}_x  \bar{g}_{m} \|_{L^2(\T \times \R_+)} \|_{l^2(\tau)}  \: \lesssim \: \| j^{\alpha} \, \| \omega \|_{{\cal H}^j_\gamma}  \|_{l^2(\tau)}   
\end{equation}
We leave the details to the reader. 

\mspace
As regards $C_{2,j}$, 
\begin{equation} \label{C2j}
\begin{aligned}
\| C_{2,j} \|_{l^2(\tau)} \: & \lesssim \: \| \, \| \omega \|_{{\cal H}^j_\gamma} \|_{l^2(\tau)} \, \sum_{m=1}^4 \| j^{1/4} \, \|\pa^{j-5}_x \hat{g}_m \|_{L^2(\T \times \R_+)} \|_{l^2(\tau)} \\
&  + \:  \|  j^{1/4} \, \| \omega \|_{{\cal H}^j_\gamma} \|_{l^2(\tau)}  \sum_{m=2}^5 \left( \| \: \|   \pa^{j-5}_x \hat{g}_m \|_{L^2(\T \times \R_+)} \|_{l^2(\tau)} +  \| \: \| \pa_y \pa^{j-6}_x \hat{g}_m \|_{L^2(\T \times \R_+)} \|_{l^2(\tau)} \right) \\
&  \lesssim \: \sqrt{\pa_\tau  E_\omega(t,\tau)} 
\end{aligned}
\end{equation}
Again, we have applied inequalities of type \eqref{barg1}-\eqref{barg2}, with $\bar{g}$ replaced by $\hat{g}$.  

\subsubsection*{Estimate on $C_{4,j}$}
Once more, one can follow closely the previous ideas, and derive
\begin{equation} \label{C4j} 
\| C_{4,j} \|_{l^2(\tau)} \: \lesssim \: \sqrt{E_\omega(t,\tau)} \| j^{1/4} \,   \| \omega \|_{{\cal H}^j_\gamma} \|_{l^2(\tau)} \: \lesssim \sqrt{\pa_\tau  E_\omega(t,\tau)}
\end{equation}

\subsubsection*{Estimate on $C_{3,j}$}
We remind that 
$$ (\pa_t + u \pa_x + v \pa_y - \pa^2_y ) \omega = 0$$
 so that 
 $$  (\pa_t + u \pa_x + v \pa_y - \pa^2_y ) \pa_y \omega = - \pa_y u \pa_x \omega - \pa_y v \pa_y \omega = - \pa_y u \pa_x \omega + \pa_x u \pa_y \omega. $$
 Hence, 
 $$ \pa_x^{j-5} {\cal C}_{3,j} = \sum_{k=0}^{j-5} \binom{j-5}{k} \pa^k_x \left( - \pa_y u \pa_x \omega + \pa_x u \pa_y \omega \right) \pa^{j-k}_x u. $$
Denoting  $a_j :=  \| \pa_x^j \left( - \pa_y u \pa_x \omega + \pa_x u \pa_y \omega \right) \|_{H^2}$, we obtain (with \eqref{sobolev}, \eqref{hardy1})`
\begin{equation*}
C_{3,j} \:  \le \: j^{1/4} \, \sum_{k=0}^{j/2} \binom{j-5}{k} a_k \| \omega \|_{{\cal H}^{j-k}_\gamma} 
 + \sum_{k=j/2}^{j-5}   \binom{j-5}{k}  a_k \sum_{m=0}^1 \| \omega \|_{{\cal H}^{j-k+m}_\gamma}  \: = \: C_{3,j}^1 \: + \: C_{3,j}^2. 
\end{equation*}
We get
\begin{equation*}
C_{3,j}^1 \: \lesssim \:   \sum_{k=0}^{j/2} \binom{j}{k} a_{(k-3)+3} \, (j-k)^{1/4} \| \omega \|_{{\cal H}^{j-k}_\gamma} 
\end{equation*} 
because $\binom{j-5}{k} \lesssim \binom{j}{k}$ and $\: j^{1/4} \lesssim (j-k)^{1/4}$ for $0 \le k \le j/2$. We then use Lemma \ref{binom} (with $m=3$) to obtain  
\begin{equation} \label{C3j1}
\| C_{3,j}^1 \|_{l^2(\tau)} \lesssim \| a_{j-3} \|_{l^2(\tau)} \,  \| j^{1/4} \|  \omega \|_{{\cal H}^j_\gamma}\|_{l^2(\tau)} 
\end{equation}
Then, we write 
\begin{align*}
C_{3,j}^2 \: & = \: j^{1/4} \sum_{k = j/2+3}^{j-2} \binom{j-5}{k-3} a_{k-3} \sum_{m=3}^4   \| \omega \|_{{\cal H}^{j-k+m}_\gamma}  \\
& \lesssim \: \sum_{k=j/2}^j \binom{j}{k} k^{1/4} a_{k-3} \sum_{m=3}^4   \| \omega \|_{{\cal H}^{j-k}_\gamma}. 
\end{align*}
because $\binom{j-5}{k-3} \lesssim \binom{j}{k}$ and $\: j^{1/4} \lesssim k^{1/4}$ for $j/2 \le k \le j$. Hence, 
\begin{equation}  \label{C3j2}
\| C_{3,j}^2 \|_{l^2(\tau)} \: \lesssim  \: \| j^{1/4} a_{j-3} \|_{l^2(\tau)} \, \| \, \| \omega \|_{{\cal H}^j_\gamma}\|_{l^2(\tau)} \: \lesssim \:  \| j^{1/4} a_{j-3} \|_{l^2(\tau)} 
\end{equation}
Finally, applying the same type of arguments, it is shown that 
$$ \| a_{j-3} \|_{l^2(\tau)} \lesssim  \| \, \| \omega \|_{{\cal H}^j_\gamma}\|_{l^2(\tau)} \: \lesssim \: \sqrt{E_\omega(t,\tau)}, \quad 
 \|j^{1/4} a_{j-3} \|_{l^2(\tau)} \lesssim  \| j^{1/4} \|  \omega \|_{{\cal H}^j_\gamma}\|_{l^2(\tau)} \: \lesssim \: \sqrt{\pa_\tau E_\omega(t,\tau)} $$
Inserting these bounds into \eqref{C3j1}-\eqref{C3j2} and summing: 
\begin{equation} \label{C3j}
\| C_{3,j} \|_{l^2(\tau)} \: \lesssim \: \sqrt{E_\omega(t,\tau)} 
\end{equation}

\subsubsection*{Estimate on $D_j$}
We now turn to the most delicate term 
\begin{align*} D_j  = & -2  j^{3/2} \int_{\T \times \R_+} \pa_x^{j-5} {\cal C}_5 \,  \tilde g_j = - 2 j^{3/2}\int_{\T \times \R_+}  \pa_y \omega \pa_x^j \pa_y \omega \,   \tilde g_j \\
& - 2 j^{3/2} \int_{\T \times \R_+} [\pa^{j-5}_x , \pa_y \omega ] \pa_y \pa_x^5 \omega \, \tilde g_j \: := \: D^1_j + D^2_j 
\end{align*}
 
 \mspace
 {\em Study of $D^2_j$.} We integrate by parts to find 
 $$ D^2_j = 2 j^{3/2} \int_{\T \times \R_+} [\pa^{j-5}_x , \pa_y \omega ]  \pa_x^5 \omega \, \pa_y \tilde g_j  \: + \: 2 j^{3/2} \int_{\T \times \R_+} [\pa^{j-5}_x , \pa^2_y \omega ]  \pa_x^5 \omega \,  \tilde g_j \: := \: D^{2'}_j \: + \: D^{2''}_j  $$
 Note that there is no boundary term, as $(\pa_x^k \pa_y \omega)\vert_{y=0} = \pa_x^k \pa^2_y u\vert_{y=0} = 0$ for any $k$, using the Prandtl equation. 
We have by Cauchy-Schwarz inequality :
$$ \|  D^{2'}_j \|_{l^1(\tau)} \: \le \:  \| \: j^{3/4} \| [\pa^{j-5}_x , \pa_y \omega ]  \pa_x^5 \omega \|_{L^2(\T \times \R_+)} \|_{l^2(\tau)} \:  \| \:  j^{3/4} \|\pa_y \tilde g_j \|_{L^2(\T \times \R_+)} \|_{l^2(\tau)}.$$
We find, by now standard arguments:
$$ \| j^{3/4} \| [\pa^{j-5}_x , \pa_y \omega ]  \pa_x^5 \omega \|_{L^2(\T \times \R_+)} \|_{l^2(\tau)}\: \lesssim \:  \| j^{3/4}  \, j \| \omega \|_{{\cal H}^{j-1}_\gamma} \|_{l^2(\tau)} 
$$
We notice that $\| j^{7/4}  \| \omega \|_{{\cal H}^{j-1}_\gamma} \|_{l^2(\tau)} \:\lesssim \: \| \| \omega \|_{{\cal H}^j_\gamma} \|_{l^2(\tau)}$ to conclude that 
$$ \|  D^{2'}_j \|_{l^1(\tau)} \: \lesssim \: \sqrt{E_\omega(t,\tau)} \, \sqrt{D^2_g(t,\tau)} \: \lesssim \: \eta D^2_g(t,\tau) \: + \: C_\eta E_\omega(t,\tau) $$ 
The treatment of $D^{2''}_j$ is similar: we state
$$ \|  D^{2''}_j \|_{l^1(\tau)} \: \lesssim \: \sqrt{E_\omega(t,\tau)} \, \sqrt{E^2_g(t,\tau)} \: \lesssim \:  E^2_g(t,\tau) \: + \:  E_\omega(t,\tau)  $$ 
Finally, 
\begin{equation}
 \|  D^2_j \|_{l^1(\tau)} \:   \lesssim \, \eta \, D^2_g(t,\tau) \: + \: C_\eta E_\omega(t,\tau) \: + \: E^2_g(t,\tau)  
 \end{equation}

\mspace
{\em Study of $D^1_j$.}
We shall use the decomposition of $\pa^j_x u$ recalled in \eqref{decompo2}: 
\begin{equation*} 
 \pa^j_x u \: = \: C_j \, 1_{\{y > a\}} \, \omega \: + \: u^j_g.
 \end{equation*}
Let $\eps > 0$, $\chi^{\eps,a} \in C^{\infty}_c(\R)$, satisfying 
$$\chi^{\eps,a} = 1 \quad \mbox{ over } \: \{ |y-a| \le \frac{\eps}{2} \}, \quad  
\chi^{\eps,a} = 0 \quad \mbox{ outside } \: \{ |y-a| \ge \eps \}. $$
We shall keep track of the $\eps$-dependence in all inequalities below. {\em In particular, notation $f \lesssim g$ will always refer to an inequality $f \le C g$ for a constant $C$ that does not depend on $\eps$.}
We write
$$ D^1_j = -2 j^{3/2} \int_{\T \times \R_+} \chi^{\eps,a} \pa_y \omega \, \pa^j_x \pa_y \omega \, \tilde g_j \: -2 j^{3/2} \: \int_{\T \times \R_+} (1- \chi^{\eps,a})  \pa_y \omega \, \pa^j_x \pa_y \omega \, \tilde g_j \: := \: D^{1'}_j + D^{1''}_j. $$
As $\omega$ does not vanish away from $y=a$, we can write 
\begin{align*} 
D^{1''}_j & = -2j^{3/2} \int_{\T \times \R_+}  (1- \chi^{\eps,a}) \frac{\pa_y \omega}{\omega} (\pa_y  g_j - \pa^2_y \omega \pa^j_x u)  \, \tilde g_j  =  2j^{3/2} \int_{\T \times \R_+}  \pa_y \left(  (1- \chi^{\eps,a}) \frac{\pa_y \omega}{\omega} \right) g_j  \, \tilde g_j \\
& + \:  2j^{3/2} \int_{\T \times \R_+}   (1- \chi^{\eps,a}) \frac{\pa_y \omega}{\omega}  g_j \: \pa_y\tilde g_j \: + \:  2j^{3/2} \int_{\T \times \R_+}  (1- \chi^{\eps,a}) \frac{\pa_y \omega}{\omega}\pa^2_y \omega \pa^j_x u \, \tilde g_j 
\end{align*}
It follows that 
\begin{equation*} 
\| D^{1''}_j \|_{l^1(\tau)} \: \lesssim \: C_\eps  \, \| j^{3/4} \| g_j \|_{L^2} \|_{l^2(\tau)} \,  \| j^{3/4} \| \pa_y \tilde g_j \|_{L^2} \|_{l^2(\tau)} \: + \:  \| j^{1/4} \| \omega \|_{{\cal H}^j_\gamma} \|_{l^2(\tau)} \, \| j^{5/4} \| \tilde g_j \|_{L^2} \|_{l^2(\tau)}
\end{equation*}
Pondering on Lemma \ref{lemmatildegjgj}, we get 
\begin{equation*}
\begin{aligned}
\| D^{1''}_j \|_{l^1(\tau)} \: & \lesssim \: C_\eps \, \left( \sqrt{E^2_g(t,\tau)} \: + \: \sqrt{E_\omega(t,\tau)} \right) \sqrt{D^2_g(t,\tau)} \: + \:  C_\eps \,  \sqrt{\pa_\tau E_\omega(t,\tau)} \,  \sqrt{\pa_\tau E^2_g(t,\tau)} \\
&  \lesssim \:  \eta D^2_g(t,\tau) \: + \: C_\eta \left( \pa_\tau E_\omega(t,\tau) \: + \: \pa_\tau E^2_g(t,\tau) \right)
 \end{aligned}
\end{equation*}

\mspace
It remains to handle $D^{1'}_j$. Therefore, we introduce 
$$ \omega^j_h \: := 1_{y > a} C_j \pa_y \omega, \quad \omega^j_g \: := \:  \pa_y u^j_g $$
We emphasize that both $\omega^j_h$ and $\omega^j_g $ are discontinuous across $y = a$, contrary to $\pa^j_x \omega = \omega^j_h + \omega^j_g$.
We find 
\begin{align*}
D^{1'}_j & = -2 j^{3/2} \int_{\{ y > a \}} \chi^{\eps,a} \pa_y \omega \, \pa_y \omega^j_h \, \tilde g_j  \: - \:  2 j^{3/2} \int_{\{ y > a \}} \chi^{\eps,a} \pa_y \omega \, \pa_y \omega^j_g \, \tilde g_j  - \: 2 j^{3/2} \int_{\{ y < a \}} \chi^{\eps,a} \pa_y \omega \, \pa_y \omega^j_g \, \tilde g_j \\
& =  -2 j^{3/2} \int_{\{ y > a \}} \chi^{\eps,a} \pa_y \omega \, \pa_y \omega^j_h \, \tilde g_j  +  2 j^{3/2} \int_{\T \times \R_+} \pa_y \left(\chi^{\eps,a} \pa_y \omega\right) \omega^j_g \, \tilde g_j \: + \: 2 j^{3/2}  \int_{\T \times \R_+} \chi^{\eps,a} \pa_y \omega \, \omega^j_g \,  \pa_y \tilde g_j  \\ 
& \quad + \: 2 j^{3/2}   \int_{\{y = a\}}  \pa_y \omega \,  [\omega^j_g] \,  \tilde g_j \: := \:  I + II + III + IV 
 \end{align*}
To go from the first to the second equality, we have integrated by parts  the last two integrals: gathering the boundary terms at $y=a^+$ and $y=a^-$    we have obtained a jump term
$$[\omega^j_g]  \:  :=  \: \omega^j_g\vert_{y=a^+} -  \omega^j_g\vert_{y=a^-}$$
 across $y=a$.  Clearly, 
 \begin{align*}
  \| I \|_{l^1(\tau)} \:  & \lesssim \: \| \, j^{1/4} \| C_j \|_{L^2(\T)} \|_{l^2(\tau)} \, \| j^{5/4} \| \tilde g_j \|_{L^2(\T \times \R_+} \|_{l^2(\tau)}\: \lesssim \:  \sqrt{\pa_\tau E_\omega(t,\tau)} \sqrt{\pa_\tau E^2_g(t,\tau)} \\
  & \lesssim \: \pa_\tau E_\omega(t,\tau) \: + \: \pa_\tau E^2_g(t,\tau)
  \end{align*}
 with the last inequality coming from \eqref{estimCj}. Also 
 \begin{align*}
  \| II \|_{l^1(\tau)} \: & \lesssim \: \| \, j^{1/4} \| \omega^j_g \|_{L^2(\T\times \R_+)} \|_{l^2(\tau)} \, \| j^{5/4} \| \tilde g_j \|_{L^2(\T \times \R_+} \|_{l^2(\tau)}
   \\
& \lesssim \: \left( \| \, j^{1/4}  \| C_j \|_{L^2(\T)} \|_{l^2(\tau)} \:+ \:  \| \, j^{1/4}  \| \pa^j_x \omega \|_{L^2(\T\times \R_+)}  \|_{l^2(\tau)} \right)   \, \| j^{5/4} \| \tilde g_j \|_{L^2(\T \times \R_+)} \|_{l^2(\tau)}   \\
& \lesssim \: \sqrt{\pa_\tau E_\omega(t,\tau)} \sqrt{\pa_\tau E^2_g(t,\tau)} \: \lesssim \: \pa_\tau E_\omega(t,\tau) \: + \: \pa_\tau E^2_g(t,\tau)
\end{align*}
As regards $IV$, we have
\begin{equation*}
IV = - 2 j^{3/2}   \int_{\{y = a\}}  \pa_y \omega \,  [\omega^j_h] \,  \tilde g_j  \: = \: -2 j^{3/2} \int_{\{y = a\}}  (\pa_y \omega)^2 \, C_j \,  \tilde g_j 
\end{equation*}
so that 
$$ | \, IV \, |  \: \lesssim \: \| j^{1/2} \, C_j \|_{L^2(\T)} \,  \| j^{5/4} \tilde g_j \|_{L^2(\T \times \R_+)}^{1/2} \, \| j^{3/4}  \pa_y \tilde{g}_j \|_{L^2(\T \times \R_+)}^{1/2} $$
and finally 
\begin{equation*}
\begin{aligned}
\| \, IV \, \|_{l^1(\tau)} \: & \lesssim \: \sqrt{\pa_\tau E_\omega(t,\tau)} \, E^2_g(t,\tau)^{1/4} \,   D^2_g(t,\tau)^{1/4} 
& \lesssim \: \eta  D^2_g(t,\tau) \: + \: C_\eta \left(  \pa_\tau E_\omega(t,\tau) \: + \:  E^2_g(t,\tau) \right) 
\end{aligned}
\end{equation*}

\mspace
To handle $III$, we split the integral as follows: 
$$  \int_{\T \times \R_+} \chi^{\eps,a} \pa_y \omega \, \omega^j_g \, \pa_y  \tilde g_j = \int_{\T \times \R_+} \chi^{\eps,a} \pa_y \omega \, \omega^j_g \,   \pa_y g_j \: + \:   \int_{\T \times \R_+} \chi^{\eps,a} \pa_y \omega \, \omega^j_g \,  \pa_y \left( \tilde g_j - g_j \right) \: = \: i \:  +  \: ii $$
Note that for $\eps$ small enough, on the support of $\chi^{\eps,a}$, one has: 
$$g_j \, = \, \omega \pa^j_x \omega - \pa_y \omega \pa^j_x u \: =  \: \omega  \omega^j_g - \pa_y \omega  u^j_g. $$  
so that 
$$ \pa_y g_j =  \omega \pa_y   \omega^j_g - \pa^2_y \omega u^j_g. $$
\mspace
We write 
\begin{align*} 
i \: & = \: \int_{\T \times \R_+}  \chi^{\eps,a}  \pa_y \omega \, \omega \,   \frac{1}{2} \pa_y  (\omega^j_g)^2 - \int_{\T \times \R_+}  \chi^{\eps,a}    \pa_y \omega \, \pa^2_y \omega \,  \omega^j_g  \,  u^j_g \\
& = \: -  \frac{1}{2}\int_{\T \times \R_+}  \chi^{\eps,a}  (\pa_y \omega)^2   (\omega^j_g)^2 \: - \:   \frac{1}{2}\int_{\T \times \R_+}  \pa_y   \chi^{\eps,a}  \pa_y  \omega \,  \omega \,  (\omega^j_g)^2   \: - \:   \frac{1}{2}\int_{\T \times \R_+}    \chi^{\eps,a}  \pa^2_y  \omega \,  \omega \,  (\omega^j_g)^2 \\
 & \quad + \: \int_{\T \times \R_+}  \chi^{\eps,a} \pa_y^2 \omega \,  \omega^j_g  \left( g_j - \omega  \omega^j_g \right)  
\end{align*}
Note that, by assumption \eqref{lowerboundomega}, one has for $\eps > 0$ small enough: 
$$ -\frac{1}{2}\int_{\T \times \R_+}  \chi^{\eps,a}  (\pa_y \omega)^2   (\omega^j_g)^2 \: \le  \: - \eta  \, \int_{\T \times \R_+}  \chi^{\eps,a}  (\omega^j_g)^2 $$
where 
$$2 \eta \: := \: \frac{1}{2} \, \inf_{t,x} \left( \pa_y \omega(t,a(t,x)) \right)^2.$$ 
On the other hand, as $\omega$ vanishes linearly in $y-a$ at the critical curve:
 $$ \frac{1}{2}  \left|  \int_{\T \times \R_+}    \chi^{\eps,a}  \pa^2_y  \omega \,  \omega \,  (\omega^j_g)^2 \right| \: \lesssim \:  \eps \,  \int_{\T \times \R_+}    \chi^{\eps,a}   (\omega^j_g)^2 $$
 For $\eps$ small enough, we can absorb the latter  by the former, resulting in 
 $$ i \: \le \: - \frac{\eta}{2} \int_{\T \times \R_+}  \chi^{\eps,a}   (\omega^j_g)^2  \: - \:  \frac{1}{2} \int_{\T \times \R_+}  \pa_y   \chi^{\eps,a}  \pa_y  \omega \,  \omega  \,  (\omega^j_g)^2  +  \: \int_{\T \times \R_+}  \chi^{\eps,a} \pa_y^2 \omega \,  \omega^j_g   g_j $$
Note that $\pa_y \chi^{\eps, a}$ vanishes in a neighborhood of $y=a$ so that, roughly speaking, $\omega^j_g$ behaves like $g_j$ over the support of    $\pa_y \chi^{\eps, a}$. More precisely, 
$$  \left|   \int_{\T \times \R_+}  \pa_y   \chi^{\eps,a}  \pa_y  \omega \,  \omega  \,  (\omega^j_g)^2 \right|
 \: \le \: C_\eps \, \| g_j \|_{L^2(\T \times \R_+)}^2   $$
 Finally, we apply  a Young's inequality to the last term: 
 $$ \left| \int_{\T \times \R_+}  \chi^{\eps,a} \pa_y^2 \omega \,  \omega^j_g   g_j \right| \: \le \: \frac{\eta}{4}  \int_{\T \times \R_+}  \chi^{\eps,a}   (\omega^j_g)^2  \: + \: C_\eta  \|  g_j \|_{L^2(\T \times \R_+)}^2.  $$
 We obtain, for $\eta > 0$ small enough 
\begin{equation*}
 i \: \lesssim \:  - \eta \int_{\T \times \R_+}  \chi^{\eps,a}   (\omega^j_g)^2 \: + \:   \| g_j \|_{L^2(\T \times \R_+)}^2 
\end{equation*}
Finally,  we apply Young's inequality to $ii$: for all $\eta' > 0$,  
\begin{equation*}
 ii \:  \lesssim  \: \eta' \int_{\T \times \R_+}  \chi^{\eps,a}   (\omega^j_g)^2 \: + \:   C_{\eta'} \|  \sqrt{\chi^{\eps,a}}  \pa_y (\tilde g_j - g_j ) \|^2_{L^2(\T \times \R_+)} 
\end{equation*}
Combining with the bound on $i$, we find (for $\eta'$ small enough)
$$ \max(III, 0) \: = \: 2 j^{3/2} \max(i + ii, 0) \: \lesssim \:   j^{3/2} \, \| g_j \|_{L^2(\T \times \R_+)}^2 \:  +  \: j^{3/2} \,  \|  \sqrt{\chi^{\eps,a}}  \pa_y (\tilde g_j - g_j ) \|^2_{L^2(\T \times \R_+)} $$
If we gather this last bound with those on $I$, $II$ and $IV$, we conclude that 
\begin{align*} 
\| \max (D^{1'}_j, 0) \|_{l^1(\tau)} \: & \lesssim \: \eta D^2_g(t,\tau) \: + \: C_\eta \left( \pa_\tau E_\omega(t,\tau) \: + \:  \pa_\tau E^2_g(t,\tau) \right)   \\
&  \quad + \: \|   j^{3/4} \, \| g_j \|_{L^2(\T \times \R_+)} \|_{l^2(\tau)}^2 \: + \:  \|  j^{3/4} \,  \| \sqrt{\chi^{\eps,a}} \pa_y (\tilde g_j - g_j ) \|_{L^2(\T \times \R_+)} \|_{l^2(\tau)}^2 
\end{align*}
We then use Lemma  \ref{lemmatildegjgj} together with 
\begin{lemma} \label{dytildegjgj}
\begin{equation*}
 \|  j^{3/4} \,  \|  \sqrt{\chi^{\eps,a}}  \pa_y (\tilde g_j - g_j ) \|_{L^2(\T \times \R_+)} \|_{l^2(\tau)}^2 \: \lesssim \: D_h(t,\tau) \: + \:  E_\omega(t,\tau). 
 \end{equation*} 
 \end{lemma}
 This last inequality will be proved in appendix  . We obtain 
 $$ \| \max (D^{1'}_j, 0) \|_{l^1(\tau)} \: \lesssim \:   \eta D^2_g(t,\tau) \: + \:  D_h(t,\tau)   \:   +    \: C_\eta \left( \pa_\tau E_\omega(t,\tau) \: + \:  \pa_\tau E^2_g(t,\tau) \right) $$
 In turn, we can combine this  bound with those on $D^{1''}_j$, $D^2_j$, to obtain 
 \begin{equation}
  \| \max (D_j, 0) \|_{l^1(\tau)} \: \lesssim \:   \eta D^2_g(t,\tau) \: + \:  D_h(t,\tau)    \:  + \: C_\eta \left( \pa_\tau E_\omega(t,\tau) \: + \:  \pa_\tau E^2_g(t,\tau) \right)  
\end{equation} 
 Then, we collect the estimates on $A_j$, $B_j$, $C_{i,j}$. We find 
 \begin{align*}
 \pa_t E^2_g(t,\tau) \: + \: D^2_g(t,\tau) \: & \lesssim \: \left(\sqrt{D^2_g(t,\tau)} \: + \:   \sqrt{\pa_\tau E^2_g(t,\tau)} \: + \: \sqrt{\pa_\tau E_\omega(t,\tau)} \right) \sqrt{E^2_g(t,\tau)} \\ \quad  & + \:  \eta D^2_g(t,\tau) \: + \:  D_h(t,\tau)   \:  + \: C_\eta \left( \pa_\tau E_\omega(t,\tau) \: + \:  \pa_\tau E^2_g(t,\tau) \right) \\
 & \lesssim \: \eta D^2_g(t,\tau) \: + \:  D_h(t,\tau)   \:  + \: C_\eta \left( \pa_\tau E_\omega(t,\tau) \: + \:  \pa_\tau E^2_g(t,\tau) \right)
 \end{align*}
 This concludes the proof of the proposition. 

\subsection{Estimate of the (first) monotonicity energy}
In the previous paragraphs, we have derived bounds on the time variations of $\dot{E}_\omega$, $E_h$ and $E^2_g$. These bounds involve the total energy $E_\omega$ defined in \eqref{Eomega}. By Lemma \ref{relations}, this energy  is controlled by the sum of $\dot{E}_\omega$, $E_h$ and $E^1_g$. Hence, to establish our main a priori estimate (Theorem \ref{theorem2}), it remains to control  $\pa_t E^1_g$. This is the purpose of 
\begin{proposition} \label{propmono1}
Under the same assumptions as in Proposition \ref{propmono1}, one has for some $C > 0$ and all $t \in ]0,T]$:
 \begin{align} \label{Eg1-prop}
 \pa_t E_g^1(t,\tau) \: + \: D^1_g(t,\tau) \: \le \: C \,  \pa_\tau E_\omega(t,\tau).
 \end{align}
\end{proposition}
Naturally, 
$$ D^1_g(t,\tau) \: := \:  \sum_{j\in \N} \left( \tau^j \,    (j!)^{-7/4} \, (j+1)^{10}\right)^2 \,  \| \pa_y g_j(t,\cdot) \|_{L^2(L^2_\gamma)}^2 $$

\mspace
To prove this proposition, we first write the equation on 
$$ g_j := \tilde \psi \left( \pa^j_x \omega - \frac{\pa_y \omega}{\omega}  \pa^j_x u\right), \quad \tilde \psi := \psi \omega + (1-\psi). $$
 We apply the operator $\tilde \psi \pa^j_x$ to the vorticity equation, apply 
 $\tilde \psi \frac{\pa_y \omega}{\omega} \pa^j_x$ to the velocity equation, and subtract one from another. A straightforward computation gives 
  \begin{align*}
& \pa_t g_j  + u \pa_x g_j + v \pa_y g_j - \pa^2_y g_j  \\
& = -\sum_{k=1}^j \binom{j}{k} \pa^k_x u \, g_{j-k+1} - \sum_{k=1}^{j-1} \binom{j}{k} \pa^k_x v \,  \tilde \psi \left( \pa_y \pa^{j-k}_x \omega  - \frac{\pa_y \omega}{\omega}  \pa^{j-k}_x \omega \right) \\
 & + \left( \pa_t + u \pa_x + v \pa_y - \pa^2_y \right)\tilde \psi \, \pa^j_x \omega   -   \left( \pa_t + u \pa_x + v \pa_y - \pa^2_y \right)\left( \tilde \psi \frac{\pa_y \omega}{\omega} \right) \pa^j_x u \\
& - 2 \pa_y \tilde \psi \pa_y \pa^j_x \omega \: + \: 2 \pa_y \left( \tilde \psi \frac{\pa_y \omega}{\omega} \right) \pa^j_x \omega \: := \: \sum_{i=1}^6 {\cal C}_i  
\end{align*}
 Then, one multiplies by $\left( \tau^j \,    (j!)^{-7/4} \, (j+1)^{10}\right)^2 (1+y)^{2\gamma} g_j$ and perform a standard energy estimate:  
 \begin{equation} \label{eqE1g}
\pa_t E^1_g(t,\tau) \: + \: D^1_g(t,\tau)\: \le \: \left( \| A_j \|_{l^2(\tau)} +   \| B_j \|_{l^2(\tau)} + \sum_{i \neq 5}  \| C_{i,j} \|_{l^2(\tau)} \right) \sqrt{\pa_\tau E^1_g(t,\tau)} \: + \: \| D_j \|_{l^1(\tau)} 
 \end{equation}
where: 
\begin{align*}
 A_j \: &  := \: \frac{1}{j^{1/2}}  (2\gamma-1) \| (1 + y)^{\gamma-1} v g_j \|_{L^2(\T \times \R_+)}, \quad 
 B_j \: := \: \frac{1}{j^{1/2}}  (2\gamma-1) \| (1+y)^{\gamma-1} \pa_y g_j \|_{L^2(\T \times \R_+)} \\
 C_{i,j} \: & :=  \:  \frac{1}{j^{1/2}} \| {\cal C}_i \|_{L^2(L^2_\gamma)}, \quad i=1,2,3,4,6, \quad 
 D_j \:  := \:   \int_{\T \times \R_+} 2 (1+y)^{2\gamma} \pa_y \tilde \psi \pa_y \pa^j_x \omega \, g_j. 
 \end{align*}
Again, there is no boundary term at $y=0$, as $\pa_y g_j$ vanishes there.  

\mspace
To bound above terms, one can follow the calculations already performed, notably those related to Proposition \ref{propvorticity}.  Note that to control the behavior at infinity of the  ${\cal C}_i$'s, one must have information on the behavior of the function $\frac{\pa_y \omega}{\omega}$ at infinity. Such information is provided by assumption \eqref{lowerboundomega}.  For the sake of brevity, we do not detail the computations. We have
\begin{equation*}
\| A_j \|_{l^2(\tau)} \: \lesssim \: \sqrt{E_\omega(t,\tau)}, \quad \| B_j \|_{l^2(\tau)} \: \lesssim \: \sqrt{D^2_g(t,\tau)}
\end{equation*}
As regards  the nonlinear terms, we have
\begin{equation*}
\| C_{i,j} \|_{l^2(\tau)} \: \lesssim \: \sqrt{E_\omega(t,\tau)} \, \sqrt{\pa_\tau E_\omega(t,\tau)} \: \lesssim  \sqrt{\pa_\tau E_\omega(t,\tau)}, \quad i=1,2, 
\end{equation*} 
whereas 
\begin{equation*}
 \| C_{i,j} \|_{l^2(\tau)} \: \lesssim \:  E_{\omega}(t,\tau), \quad i=3,4,6.  
\end{equation*}
Finally, after integration by parts, we find 
\begin{equation*}
\| D_j \|_{l^1(\tau)} \: \le \: C \sqrt{E_\omega(t,\tau)} \, \left( \sqrt{E^1_g(t,\tau)} + \sqrt{D^1_g(t,\tau)} \right) \: \lesssim \: \eta  \, D^1_g(t,\tau) \: + \: C_\eta E_\omega(t,\tau), \quad \forall \eta > 0
\end{equation*}
Inserting  all these bounds  in \eqref{eqE1g} (and taking parameter $\eta$ small enough) yields the estimate in Proposition \ref{propmono1}.  

\subsubsection*{Conclusion: Main a priori estimate}
On one hand, from propositions \ref{propvorticity}, \ref{prophydro} and \ref{propmono1}, we infer that 
\begin{align*}
\pa_t \left( \dot{E}_\omega + E_h  + E^1_g \right)(t,\tau) \: + \: \left( \dot{D}_\omega + D_h  + D^1_g \right)(t,\tau) &  \lesssim \:  \pa_\tau E_\omega(t,\tau) + \pa_\tau E^2_g(t,\tau) \\
&  \lesssim \:   \pa_\tau \left(\dot{E}_\omega + E_h  + E^1_g \right)(t,\tau)  + \pa_\tau E^2_g(t,\tau),
\end{align*}
where we have used \eqref{relation_energies2} to go from the first to the second line. 
On the other hand, Proposition \ref{propmono2} yields (still with relation \eqref{relation_energies2}): 
\begin{equation*} 
\pa_t E^2_g(t,\tau) \: + \: D^2_g(t,\tau) \: \lesssim \:\pa_\tau \left(\dot{E}_\omega + E_h  + E^1_g \right)(t,\tau)   \: + \:  \pa_\tau E_g^2(t,\tau) \: + \: D_h(t,\tau).
\end{equation*}
Thus, we can multiply  the last equation by a small factor $\alpha$, and add it to the first equation, to end up with 
\begin{multline*}
\pa_t \left( \dot{E}_\omega + E_h  + E^1_g + \alpha E^2_g  \right)(t,\tau) \: + \: \left( \dot{D}_\omega + D_h  + D^1_g + D^2_g \right)(t,\tau) \\
 \lesssim \:  \pa_\tau \left( \dot{E}_\omega +  E_h  +  E^1_g + \alpha E^2_g  \right)(t,\tau). 
\end{multline*} 
This concludes the proof of Theorem \ref{theorem2}. We insist that the $\alpha$ depends on $M$, but can be (and in fact needs to be) taken small.  In particular, we insist that  we can choose $\alpha \le 1$, that is with an upper bound independent of $M$. This will be crucial in the effective construction of solutions.

%

\section{Existence of a Gevrey solution}
We deal here with the existence part of Theorem \ref{theorem1}. The core of the argument is the {\it  energy estimate} stated in Theorem \ref{theorem2}. Actually, such estimate is derived on a sequence of approximations of the Prandtl system, that we now present.

\subsection{Approximate system}
We start with an initial data satisfying 
 \begin{equation}  u_0 \in G^{7/4}_{\tau_0}(\T; \, H^{s+1}_{\gamma-1}),   \quad   \omega_0 \: := \: \pa_y u_0 \in   G^{7/4}_{\tau_0}(\T; \, H^s_{\gamma}),  
 \end{equation}
as well as assumptions (H) and (H'), see section \ref{statements}.  Parameters  $\tau_0$, $s$, $\sigma$, $\delta$ are as in Theorem \ref{theorem1}.  
 To fix ideas, we assume that the curve of non-degenerate critical points satisfies the constraint $0 < a_0 < 3$. The choice of the value $3$ is of course arbitrary. 
 
\mspace
To construct solutions with all the a priori bounds of 
the previous section, we need to find a good approximating scheme:  
it should  not destroy the energy estimates, and should  guarantee that the point-wise bounds \eqref{lowerboundomega}  (which echo those in (H')) are preserved in small time, uniformly  with respect to the  approximation parameters.

\mspace 
The same scheme as in  \cite[section 4]{MW13prep}  works here. It goes through the following approximate system  ({\it regularized Prandtl system})  
\begin{equation} \label{Prandtl-eps}
\left\{
\begin{aligned}
\pa_t u^\eps + u^\eps \pa_x u^\eps + v^\eps \pa_y u^\eps  - \pa^2_y u^\eps -\eps 
  \pa^2_x u^\eps  = 0 , &  \\
\pa_x u^\eps + \pa_y v^\eps = 0, &  \\
u^\eps\vert_{y=0} = v^\eps\vert_{y=0} = 0, \:\lim_{y=+\infty} u^\eps = 0. & 
\end{aligned}
\right.
\end{equation} 
The big difference with the original Prandtl system is the tangential diffusion $-\eps \pa_x^2 u$, which allows to control the loss of $x$-derivative generated by the $v$-term.  It restores well-posedness in the Sobolev setting. A detailed construction of  solutions of \eqref{Prandtl-eps} is performed in \cite{MW13prep}. Let us denote  
$$H^{s,\gamma}(\T \times \R_+) \: := \: \cap_{k=0}^s H^k_x(\T, H^{s-k}_\gamma(\R_+)).$$
 As a  result of \cite{MW13prep}, for  any initial data
$$ u_0 \in H^{s, \gamma-1}(\T \times \R_+), \quad  \omega_0 : = \pa_y u_0 \in H^{s,\gamma}(\T \times \R_+)$$ 
there is a unique  $T^\eps_*  \in \R_+^* \cup \{+\infty\}$,  and a unique  maximal  solution  
$$ u^\eps \in L^\infty_{loc}([0,T^\eps_*); H^{s,\gamma-1}(\T \times \R_+)), \quad  \omega^\eps := \pa_y u^\eps  \in L^2_{loc}([0,T^\eps_*);  H^{s,\gamma}(\T \times \R_+)).   $$
of \eqref{Prandtl-eps}. By maximal, we mean that  
\begin{equation}
\mbox{ either }  \: T^\eps_* = +\infty, \quad \mbox{ or } \:  \limsup_{t \rightarrow T^\eps_*} \| \omega^\eps(t) \|_{H^{s,\gamma}(\T \times \R_+)} = +\infty 
 \end{equation}
 Of course, our  Gevrey data $u_0$ is regular enough to apply this result. {\em Furthermore, up to some small time (possibly depending on $\eps$), the solution $u^\eps$ (resp. $\omega^\eps$) remains  in $\displaystyle G^{7/4}_{\tilde\tau_0}(H^{s+1}_{\gamma-1})$ (resp. $\displaystyle G^{7/4}_{\tilde\tau_0}(H^s_\gamma)$), for any $\tilde \tau_0 < \tau_0$.}     

\mspace
The fact that system \eqref{Prandtl-eps} preserves Gevrey regularity in small time is somehow classical. One  way to show it is to establish Gevrey bounds on $u^{\eps,n}$, solution of  the Galerkin type approximation 
\begin{equation*} 
\left\{
\begin{aligned}
\pa_t u^{\eps,n} + P^n \left( u^{\eps,n} \pa_x u^{\eps,n} + v^{\eps,n} \pa_y u^{\eps,n} \right)  - \pa^2_y u^{\eps,n} -\eps 
  \pa^2_x u^{\eps,n}  = 0 , &  \\
\pa_x u^{\eps,n} + \pa_y v^{\eps,n} = 0, &  \\
u^{\eps,n}\vert_{y=0} = v^{\eps,n}\vert_{y=0} = 0, \:\lim_{y=+\infty} u^{\eps,n} = 0. & 
\end{aligned}
\right.
\end{equation*} 
$P_n$ is the projection over the Fourier modes $|k| \le n$ in variable $x$. One can show for this system uniform (in $n$) Gevrey bounds for small time $T^\eps$ (independent of $n$). As $n$ goes to infinity, we get Gevrey bounds for the (unique) solution $u^\eps$ of \eqref{Prandtl-eps}.  We insist that these estimates are much simpler than those of the previous section. No special structure or tricky norm is needed: as the system is fully parabolic, there is no loss of $x$-derivative.   For the sake of brevity, we do not give further details. 
 
  \mspace
Now, let  $M \ge M_0$ large enough,  and  $t \mapsto \tau(t) \in [\frac{\tilde \tau_0}{2}, \tilde \tau_0]$ a function over $\R$. We denote     
 $$ T^\eps(\tau) := \sup \left\{ \, T \in [0,T^\eps_*[, \quad  E_\omega(t,\tau(t)) \le M \: \mbox{ for all } \: t \in [0,T]  \right\}, $$ 
   where $E_\omega$ is defined in \eqref{Eomega}, replacing $u$ by $u^\eps$). By standard arguments, $T^\eps(\tau) > 0$, and 
   $ \limsup_{t \rightarrow T^\eps(\tau)}  E_\omega(t,\tau(t)) = M.$  
   
\mspace
{\em The point is to show that for a good choice of $M$ and $\tau$, there exists $T^0 > 0$ independent of $\eps$, such that $T^\eps(\tau) \ge T^0$}. Then,  the solution $u^\eps$ will exist  and have uniform Gevrey bounds on a time independent of $\eps$. Finally, standard compactness arguments  will allow to let $\eps \rightarrow 0$, and obtain a solution of the Prandtl equation with initial data $u_0$ at the limit.     
      
\mspace
We now give a few hints on how to prove the inequality $T^\eps(\tau) \ge T^0$. Let $M \ge M_0$. The key-point is to establish the following: there exists $T_0 > 0$, depending on $M$, on the initial data and on $\inf \tau = \frac{\tilde \tau_0}{2}$, {\em but  independent of $\eps$},  such that:
\begin{itemize}
\item For all  $t \mapsto \tau(t) \in [\frac{\tilde \tau_0}{2}, \tilde \tau_0]$, for all $ t \in [0, \min(T^0, T^\eps(\tau))]$,  
\begin{equation}  \label{approxbound1}
|\omega^\eps(t,x,y)| \ge \frac{\delta}{(1+y)^\sigma}, \quad | \pa^\alpha \omega^\eps(t,x,y) | \le \frac{1}{\delta (1+y)^{\sigma+\alpha_2}}.
 \end{equation}
 for the same $\delta$ as in $(H')$. 
\item  There exists $\alpha \le 1$,  and a function  $t \mapsto \tau(t)  \in [\frac{\tilde \tau_0}{2}, \tilde \tau_0]$  such that:     
 \begin{equation} \label{approxbound2}
\mbox{ For all } \: t \in [0, \min(T^0, T^\eps(\tau))], \quad  \frac{d}{dt} {\cal E}(t,\tau(t)) \: \le \: 0, 
\end{equation}
with 
$$ {\cal E}(t,\tau) = \dot{E}_\omega(t,\tau) \: + \:   E_h(t,\tau) \: + \: E^1_g(t,\tau) \: + \: \alpha \,  E^2_g(t,\tau),  $$
where the energy functionals $\dot{E}_\omega$, $\: E_h(t,\tau)$, $E^i_g$ are defined respectively in \eqref{vorticityenergy}, \eqref{hydroenergy}, \eqref{monoenergy1} and \eqref{monoenergy2}, replacing $u$ by $u^\eps$.
\end{itemize} 

\mspace
We shall discuss the bounds \eqref{approxbound1} and \eqref{approxbound2} in the next paragraph. We explain first how it allows to conclude. From \eqref{approxbound2}, we deduce that for all $t \in [0, \min(T^0, T^\eps)]$, 
$$\: {\cal E}(t,\tau(t)) \: \le \:   {\cal E}(0,\tau(0)).$$
As $\alpha \in (0,1)$, this implies in turn 
$$  \dot{E}_\omega(t,\tau(t)) \: + \:   E_h(t,\tau(t)) \: + \: E^1_g(t,\tau(t)) \: \le \:   \dot{E}_\omega(0,\tau(0)) \: + \:   E_h(0,\tau(0)) \: + \: E^1_g(t,\tau(0)) + E^2_g(0,\tau(0)). $$
By Lemma \ref{relations}, we finally get 
$$ E_\omega(t,\tau(t)) \: \le C \: \left( \dot{E}_\omega(t,\tau(t)) \: + \:   E_h(t,\tau(t)) \: + \: E^1_g(t,\tau(t)) \right)   \: \le C \: \left(   E_\omega(0,\tau(0)) \: + \:  E^2_g(0,\tau(0)) \right).  $$
A closer look at the proof of Lemma \ref{relations} shows that the constant $C$ depends on the initial data.  {\em But,  up take a smaller $T_0$,  it does not depend on the constant $M$}. For instance, it depends on $\inf_{t, x \in \T} \pa_y \omega(t,a(t,x))$: by  Sobolev imbeddings and the bound $E_\omega(t,\tau(t)) \le M$, one has for $T_0 =T_0(M)$ small enough: 
$$    \inf_{t, x \in \T} \pa_y \omega(t,a(t,x)) \: \ge \: \frac{1}{2}  \inf_{x \in \T} \pa_y \omega(x,a_0(x)). $$
We leave the details to the reader.  Hence, we can take 
$$ M >  2 C  \dot{E}_\omega(0,\tau(0)) \: + \:   E_h(0,\tau(0)) \: + \: E^1_g(t,\tau(0)) + E^2_g(0,\tau(0)).$$ 
This ensures that ${\cal E}(t,\tau(t)) \le \frac{M}{2}$ over $[0, \min(T^0, T^\eps(\tau))]$, and from there that $T^\eps(\tau) \ge T^0$.  

\subsection{Uniform bounds}

\subsubsection{Maximum principle: proof of \eqref{approxbound1}}
This paragraph is devoted to the proof of \eqref{approxbound1}. One must show the upper and lower bounds over a time interval $[0,\min(T^0,T^\eps)]$ for some $T^0$ independent of $\eps$. Note that the Sobolev imbeddings and the bound $E_\omega(t,\tau(t)) \le M$ allow to bound some weighted $L^\infty$ norms of $\omega^\eps$ and some derivatives. However, the bounds in \eqref{approxbound1} can not be directly deduced from there, due to the weight $(1+y)^\sigma$, $\sigma \gg \gamma$.
Nevertheless, as the regularized Prandtl equation is parabolic, it can be deduced from a maximum principle. This is described in details in article \cite{MW13prep}. The proof of the upper bounds
$$  | \pa^\alpha \omega^\eps(t,x,y) | \le \frac{1}{\delta (1+y)^{\sigma+\alpha_2}}, \quad \forall |\alpha| \le 2, $$
is exactly the same as in \cite{MW13prep}: see inequality (5.5),  and section 5.2. As regards the lower bound on $(1+y)^\sigma \, \omega^\eps$, there is only a slight modification.   In \cite{MW13prep},  $\omega \neq 0$ for all $y \ge 0$, so that the minimum principle can be applied over domains of the type $[0,T] \times \T \times \R_+$. In our situation, we must place ourselves above the curve of critical points. More precisely, as $\omega_0(x,a_0(x))$ = 0, $\pa_y \omega_0(x,a_0(x)) \neq 0$, there exists $\eta > 0$, such that $\omega_0(x,a_0(x)+\eta)   \neq 0$. One then applies the minimum principle over $[0,T] \times 
\{ (x,y), \quad y = a_0(x) + \eta \}$, $T$ small enough. We refer to \cite{MW13prep}, inequality (5.6) and section 5.2.      

\subsubsection{Energy estimates: proof of \eqref{approxbound2}}
The point here is to derive the estimates of Theorem \ref{theorem2}, uniformly in $\eps$.  More precisely, one can show that for any $t \mapsto \tau(t)  \in   [\frac{\tilde \tau_0}{2}, \tilde \tau_0]$,
\begin{equation}  \label{energyeps}
 \pa_t {\cal E}(t,\tau(t)) \: \le C \:  \pa_\tau {\cal E}(t,\tau(t)) 
 \end{equation}
 with a constant $C$ that does not depend on $\eps$ (but depends on $M$, the initial data, and $\inf \tau = \frac{\tilde \tau_0}{2}$). Then,  
$$ \frac{d}{dt} {\cal E}(\alpha,t,\tau(t)) = \pa_t {\cal E}(\alpha,t,\tau(t)) + \tau'(t) \pa_\tau  {\cal E}(\alpha,t,\tau(t))  \leq  (C + \tau'(t))   {\cal E}(\alpha,t,\tau(t)). $$
Then, choosing $\tau(t) = \tilde \tau_0  - C  t$ leads to \eqref{approxbound2} (up to reduce $T^0$ so that the constraint $\tau \ge \frac{\tilde \tau_0}{2}$ remains satisfied).   

\medskip
The proof of \eqref{energyeps}, on the approximate system \eqref{Prandtl-eps},  mimics the proof of the {\it a priori} estimate in Theorem \ref{theorem2}, on the exact Prandtl system. Indeed, the additional term $-\eps \pa^2_x u$ does not raise any problem. There are only two noticeable differences: 
\begin{itemize}
\item  In the energy estimate for $\dot{E}_\omega$, the estimate of the boundary term 
$$ E_j \:  :=  \:  \sum_{\substack{J = (j_1,j_2) \in \N^2 \\\, |J| = j, \, 0 <  j_2 \le s}}   \int_{\T \times \{0\}} \pa_y \omega_J \, \omega_J $$
changes slightly. Indeed, we recall that the estimate of $E_j$ goes through the Lemma \ref{reductionlemma}: it uses the Prandtl equation to reduce 
the number of $y$ derivatives in the expression of $\pa_y^k \omega\vert_{y=0}$, $k$ odd. Due to the regularization, the equation changes, and so the expression: the modified formula is given in \cite[Lemma 5.9]{MW13prep}. The new terms in the formula do not raise any serious difficulty. For the sake   of brevity, we leave the details to the reader.
\item In the energy estimate for $E^2_g$,  a new commutator appears, namely  
$$ 2\eps  \pa_x^{j-5} (\pa_{xy} \omega \, \pa_x^6 u - 2 \eps \pa_x \omega \pa_x^6 \omega ). $$
It can be expanded by Leibniz rule: for brevity, we focus on one of the most difficult terms in the sum, that is:   
$$ {\cal C}_j \:  := \:  -2 \eps \pa_x \omega \pa_x^{j+1} \omega.  $$
Proceeding as in subsection \ref{subsecmono}, we need to evaluate 
\begin{align*}
& \| j^{3/2} \int_{\T \times \R_+} {\cal C}_j \, \tilde g_j \|_{l^1(\tau)} \\
&  \lesssim \:  \eps \| j^{3/4} \| \pa_x^j \omega \|_{L^2} \|_{l^2(\tau)}  \: \| j^{3/4} \| \pa_x \tilde g_j \|_{L^2} \|_{l^2(\tau)} 
\:  + \:  \eps \| j^{1/4}   \| \pa_x^j \omega \|_{L^2} \|_{l^2(\tau)}   \: \| j^{5/4} \|  \tilde g_j \|_{L^2} \|_{l^2(\tau)}  \\
& \: \lesssim \:  C_\eta \eps  \| j^{3/4} \| \pa_x^j \omega \|_{L^2} \|_{l^2(\tau)}^2 \: + \:  \eta \eps  \| j^{3/4} \| \pa_x \tilde g_j \|_{L^2} \|_{l^2(\tau)}^2 \:  + \:  \eps \, \sqrt{\pa_\tau E_\omega(t,\tau)} \, \sqrt{\pa_\tau E^2_g(t,\tau)}. 
  \end{align*}
Note that we have performed an integration by parts to go from the first to the second line. The second term can be absorbed for small $\eta$ by the tangential dissipative term  at the left-hand side.  As regards the first term, we use that     
\begin{align*} 
  \eps & \| j^{3/4} \| \pa_x^j \omega \|_{L^2} \|_{l^2(\tau)}^2  \\
 & \lesssim \: \eps  \Bigl( \| \| j^{3/4} \pa_x h_{j-1} \|_{L^2} \|_{l^2(\tau)}^2 + \| j^{3/4} \pa_x g_{j-1} \|_{L^2} \|_{l^2(\tau)}^2 +\| \| j^{3/4} \pa_x^{j-1} \omega  \|_{L^2} \|_{l^2(\tau)}^2    \Bigr) \\
& \lesssim \: \eps  \Bigl( \| \| \pa_x h_j \|_{L^2} \|_{l^2(\tau)}^2 + \| \pa_x g_j \|_{L^2} \|_{l^2(\tau)}^2 +\| \|  \pa_x^{j} \omega   
\|_{L^2} \|_{l^2(\tau)}^2    \Bigr).
\end{align*}
\end{itemize}
Such inequality can be obtained along the lines of  Lemma \ref{relations}. Eventually, multiplying by a small constant $\alpha$ the energy estimate on $E^2_g$, this term will be absorbed by the tangential  dissipative terms related to $E_h$ and $E^1_g$. For the sake of brevity, we skip the details.

\section{Uniqueness}
 

In this section, we prove the uniqueness of the solution 
constructed in the previous section. As usual, uniqueness can be 
derived by some energy estimate in a space weaker than the 
space where existence was proved. One can hope to do it in 
$L^2$. However, due to the loss of derivative in the equation, 
we have to do it with the same Gevrey regularity.  We will only 
use a small loss of Sobolev correction. 
 
 \mspace
We assume that we have two solutions $u^1$ and  $u^2$ of our system 
\eqref{Prandtl},  such that the  functionals  ${\cal E}^1(\alpha^0,t,\tau^0) $
and  ${\cal E}^2 (\alpha^0,t,\tau^0) $, defined as in \eqref{calE} 
with $\omega$ replaced by $\omega^1$ and $\omega^2$, are bounded (in this whole section, $\omega^2 $ should not be confused with the 
square of $\omega$). Here, we can take $\alpha^0 = 1$ and $\tau^0 = \min(\tau^1,\tau^2 ) $.
We also assume the lower bound and the upper bound \eqref{lowerboundomega}
to hold for both solutions $\omega^1$ and $\omega^2$. Finally, we assume 
the existence of  critical curves $a^1(t,x)$ and $a^2(t,x)$ such 
that for $i=1,2$, we have   
\begin{equation}
 \pa_t a^i(t,x) + \frac{\pa_t \omega^i(t,x,a(t,x))}{\pa_y \omega^i(t,x,a(t,x))} = 0, \quad a^i(0,x) \: = \: a_0(x). 
 \end{equation}
We also recall that  $u^1$ and  $u^2$   should remain convex: $\pa^2_y u^i =
 \pa_y \omega^i  > 0$, for $t$ and $y-a^i(t,x)$ small enough.

 \mspace
Let us denote 
$u = u^1 - u^2$, $v = v^1 - v^2$ and $\om  = \om^1 - \om^2$. Hence, 
we have  
 \begin{equation} \label{u-d}
\pa_t u + u^1 \pa_x u + v^1 \pa_y u + u  \pa_x u^2  +  v  \pa_y u^2 
  - \pa_y^2 u  = 0,
\end{equation}
 \begin{equation} \label{om-d}
\pa_t \omega + u^1 \pa_x \omega + v^1 \pa_y \omega + u  \pa_x \omega^2  +  v  \pa_y \omega^2 
  - \pa_y^2 \om  = 0. 
\end{equation}

We are going to perform  the same type of estimates as 
 in the existence part: the relevant energy, still called 
${\cal E}(\alpha,t,\tau)$,  is defined as in \eqref{calE},   
but  the definition of
$h_j$, $g_j$ and $\tilde g_j$ used to define $E_h$, $ E^1_g$ and 
 $ E^2_g$  need to be changed. Indeed, we take  
\begin{equation} \label{hjn}
 h_j(t,x,y) \: = \: \chi(y-a^2(t,x)) \, \frac{\pa_x^j \omega}{\sqrt{\pa_y 
\omega^2}}(t,x,y),
 \end{equation}
\begin{equation} \label{gjn}
g_j(t,x,y) \: :=  \: \Big(\psi (y) \omega^2(t,x,y) + 1-  \psi(y)  \Big)  \left( \pa_x^j \omega - \frac{\pa_y \omega^2}{\omega^2} \pa_x^j u \right)(t,x,y),
\end{equation} 
\begin{equation}
\tilde g_j(t,x,y) \: := \:  \pa_x^{j-5}\left(  \omega^2 \pa_x^5 \omega  -   \pa_y\omega^2  \pa_x^5 u \right).
\end{equation}
with $\psi = \psi(y) \in C^\infty_c(\R)$ equal to $1$ in an open  neighborhood of  $[0,\sup |a^2|]$. Also, we replace the constraint $j_2 \le s$ in the definition of $\dot{E}_\omega$  by $j_2 \le s-2$. Finally, in all energies, we replace  the Sobolev correction term, namely the factor $(j+1)^{10}$, 
by $(j+1)^8$. 

\mspace
Now, we follow the strategy used 
in the proof of the a priori estimates. We do not give all the details,  but just mention the main changes 
we need to make.

\subsection{Estimate of the vorticity energy}
 

We differentiate  the vorticity equation \eqref{om-d}   $J$ times, 
 $ J = (j_1,j_2) \in \N^2$, \: $0 < j_2 \le s-2$.
We find, for  $\omega_J \: := \: \pa^{J} \omega$:
\begin{align*}
 \pa_t \omega_J  + u^1 \pa_x  \omega_J + v^1 \pa_y \omega_J 
 - \pa^2_y \omega_J  + [\pa^J, u^1] \pa_x \omega + [\pa^J, v^1] \pa_y \omega \\  
+  \pa^{J} ( u \pa_x  \omega^2 + v \pa_y \omega^2  ) 
\end{align*} 
Arguing as in the proof of Proposition  \ref{propvorticity}, we 
deduce that 

\begin{align*}
 \pa_t \dot{E}_\omega(t,\tau) \: + \: \dot{D}_\omega(t,\tau) \:  
&  \lesssim   C \, \pa_\tau E_\omega(t,\tau).  
\end{align*} 
We only point out the presence of the extra terms 
$ \pa^{J} ( u \pa_x  \omega_J^2 + v \pa_y \omega_J^2  )  $
which can be treated along the estimates of 
$A_j$ and $B_j$ in paragraph \ref{subsecvorticity}. Note here 
that  $u$ and $v$ can only be hit by $ J$ derivatives. Moreover,   
due to the difference between the Sobolev correction factors 
 $(j+1)^8$ and $(j+1)^{10}$, we can  estimate $\pa^{J}   \pa_x  \omega_J^2 $ 
 and $ \pa^{J} \pa_y \omega_J^2$ by the new ${\cal E}^2$, despite the extra $x$ and $y$ derivative.   
This is the reason we took a slight loss in Sobolev correction.  Also, thanks to the new restriction $j_2 \le s-2$, not more than $s$  derivatives with respect to $y$ hit $\omega^2$. 


\subsection{Estimate of the hydrostatic energy}

As in subsection \eqref{subsechydro}, we write an 
equation for $h_j$. It is now:

\begin{align*}
\pa_t h_j + u^1 \pa_x h_j + v^1 \pa_yh_j - \pa^2_y h_j \: & = \: 
\frac{\chi^a}{\sqrt{\pa_y \omega^2}} [ \pa^j_x, u^1 ] \pa_x \omega 
\:  + \:  \frac{\chi^a}{\sqrt{\pa_y \omega^2}}   [ \pa^j_x, v^1 ] \pa_y \omega  \\
&  + \: [  \frac{\chi^a}{\sqrt{\pa_y \omega^2}}, \pa_t + u^1 \pa_x  + v^1 \pa_y - \pa^2_y] \pa^j_x \omega  \: -  \:   \frac{\chi^a}{\sqrt{\pa_y \omega^2}} 
  \pa^j_x v \pa_y \omega^2 \\
& -  \frac{\chi^a}{\sqrt{\pa_y \omega^2}}\pa^j_x (u  \pa_x \omega^2)    -   \frac{\chi^a}{\sqrt{\pa_y \omega^2}} [ \pa^j_x,  \pa_y \omega^2  ] v 
\end{align*}
Performing a standard energy estimate on the previous equation,
 multiplying by 
$ \left(\tau^{j} (j!)^{-7/4} \, j^{8}\right)^2$,  summing over $j$, we end up with 
and equation similar to \eqref{Ehd}.  The main term now is 
 \begin{align*} 
E_j \: & := \: \int_{\T \times \R_+} \frac{\chi^a}{\sqrt{\pa_y \omega^2}}  \pa^j_x v \, \pa_y \omega^2 \; h_j 
 \end{align*}
which can be treated exactly as in  subsection \ref{subsechydro}. 
Hence, we deduce that the estimate \eqref{Eh-prop} holds with the 
new definition of $u,\omega$  and with the change of the Sobolev correction in 
the energies.

\subsection{Estimate of the (second) monotonicity energy} 

As in subsection \eqref{subsechydro}, 
 we apply $\omega^2 \pa_x^5$ to the vorticity equation \eqref{om-d}, 
$\pa_y \omega^2 \pa_x^5$ to the velocity equation \eqref{u-d}, 
and substract  one from the other:  
\begin{equation}   \label{g5d}
\begin{aligned}
\pa_t \tilde g_5 + u^1 \pa_x \tilde g_5 + v^1 \pa_y \tilde g_5 - \pa^2_y \tilde g_5 = & - \sum_{k=1}^5 \binom{5}{k} \pa^k_x u^1 \, \bar{g}_{5-k+1} \: - \: \sum_{k=1}^5 \binom{5}{k}
  \pa_x^k v^1 \,  \hat g_{5-k+1} \\
& - \sum_{k=0}^5  \binom{5}{k}\pa^k_x u \, \bar{g}_{5-k+1}^2  \: - \: \sum_{k=0}^4
  \binom{5}{k} \pa_x^k v \,  \hat g_{5-k+1}^2 \\
& - (\pa_t + u^1 \pa_x + v^1 \pa_y - \pa^2_y ) \pa_y \omega^2 \, \pa^5_x u + 2 \pa^2_y \omega^2 \, \pa^5_x \omega - 2 \pa_y \omega^2 \, \pa_y \pa^5_x \omega \\
& +  (\pa_t + u^1 \pa_x + v^1 \pa_y - \pa^2_y )  \omega^2 \, \pa^5_x \omega 
\end{aligned}
\end{equation}
with   $  \hat g_{k}  $, $  \bar g_{k}  $, $  \hat g_{k}^2  $, $  \bar g_{k}^2  $ defined 
as in        \eqref{ggg},  for instance
 $\bar g_k \: := \: \omega^2 \pa^k_x \omega  -  \pa_y \omega^2 
\pa^k_x u$   and 
  $\bar g_k^2 \: := \: \omega^2 \pa^k_x \omega^2  -  \pa_y \omega^2 
\pa^k_x u^2$
for $k \le 5$. Note that the  term
that caused a loss of one $x$ derivative has disappeared. 
As in the proof of Proposition \ref{propmono2}, we apply 
$\pa_x^{j-5}$ to the equation \eqref{g5d} and we conclude in a  similar way
that \eqref{Eg-prop} holds.  

\subsection{Estimate of the (first) monotonicity energy}
As in the proof of Proposition \ref{propmono1}, we write an 
equation for $ g_j$ defined in \eqref{gjn}: 

  \begin{align*} 
& \pa_t g_j  + u^1 \pa_x g_j + v^1 \pa_y g_j - \pa^2_y g_j  \\
& = -\sum_{k=1}^j \binom{j}{k} \pa^k_x u^1 \, g_{j-k+1} - \sum_{k=1}^{j} \binom{j}{k} \pa^k_x v^1 \,  \tilde \psi \left( \pa_y \pa^{j-k}_x \omega  - 
 \frac{\pa_y \omega^2}{\omega^2}  \pa^{j-k}_x \omega \right) \\ 
&  -\sum_{k=0}^j \binom{j}{k} \pa^k_x u \, g_{j-k+1}^2 - \sum_{k=0}^{j-1} \binom{j}{k} \pa^k_x v  \,  \tilde \psi \left( \pa_y \pa^{j-k}_x \omega^2  - 
 \frac{\pa_y \omega^2}{\omega^2}  \pa^{j-k}_x \omega^2 \right) \\ 
 & + \left( \pa_t + u^1 \pa_x + v^1 \pa_y - \pa^2_y \right)\tilde \psi \, \pa^j_x \omega   -   \left( \pa_t + u^1 \pa_x + v^1 \pa_y - \pa^2_y \right)
 \left( \tilde \psi \frac{\pa_y \omega^2}{\omega^2} \right) \pa^j_x u \\
& - 2 \pa_y \tilde \psi \pa_y \pa^j_x \omega \: + \: 2 \pa_y \left( \tilde \psi \frac{\pa_y \omega}{\omega} \right) \pa^j_x \omega \:  
\end{align*}
where $ \tilde \psi  =   \Big(\psi (y) \omega^2(t,x,y) + 1-  \psi(y)  \Big) $. Again the main 
thing here is that in the third line the term that had 
a potential loss of $x$ derivative, namely the one 
 involving $\pa^j_x v  $
was canceled.  As above, we get an estimate similar to 
\eqref{Eg1-prop}.

Putting all these estimates together, we can easily conclude,  
as in the a priori estimate section,  by a Gronwall lemma that 
$u = 0$.

\section{Conclusions}

We have proved in this paper the local wellposedness of  the Prandtl system 
for data that are  in the  Gevrey  class $7/4$ in the horizontal 
 variable $x$.
We would like to discuss here some possible extension of this work:

\begin{itemize} 

\item  First, it is unlikely that the  Gevrey  class $7/4$ is optimal. In particular, the analysis and numerics performed in \cite{GD10} (on a simple linearization) suggest that the optimal exponent may be $s=2$.  However, we believe that 
extending our result  to the  Gevrey  class $2$ or even $2-\eps$ for all $\eps >0$ would require some new ideas. Moreover, it remains uncertain  that   the index $s=2$ is the critical one since  there may be more severe instabilities than those constructed in \cite{GD10}.

\item For simplicity, we have considered here homogeneous source terms: $U=0$ and $P=0$. However, our work should  extend  to the case where $U$ is a function of 
$t,x$, along the lines of \cite{KV13}. 
 Also, we have considered  data that have polynomial decay 
when $y$ goes to infinity. One can adapt the proof to  data with exponential decay, 
using some Hardy type inequality with exponential weights.

\item An important open question is to prove global existence 
 for small data for the Prandtl system. 

\item Finally, we recall that  the study of the Prandtl system is motivated 
by the zero viscosity limit in the presence of a boundary. 
The present work opens the road to  a full justification of the Euler-Prandtl 
description, without analyticity.

\end{itemize}

\appendix

\section{Sobolev  and Hardy inequalities}


We  recall   classical Sobolev  and Hardy inequalities, see \cite{KMP07}:
 For $f = f(x,y)$, 
we have   
\begin{equation} \label{sobolev} 
\| f \|_{L^\infty}  \: \lesssim \: \| f \|_{L^2} \: + \: \| \pa_x f \|_{L^2} \: + \:  \| \pa_y f \|_{L^2} \: + \:     \| \pa_y \pa_x f \|_{L^2}.
\end{equation} 
If $\lambda > -\frac{1}{2}$ and $\lim_{y\rightarrow +\infty} f(x,y) = 0$, then 
\begin{equation} \label{hardy1}
\| (1+y)^\lambda f \|_{L^2} \: \le \: \frac{2}{2\lambda+1} \| (1+y)^{\lambda+1} \pa_y f \|_{L^2}  
\end{equation}
If $\lambda < -\frac{1}{2}$, we have 
 \begin{equation} \label{hardy2}
\| (1+y)^\lambda f \|_{L^2} \: \le \: \sqrt{-\frac{1}{2\lambda+1}}\| f\vert_{y=0} \|_{L^2} - \frac{2}{2 \lambda+1} \| (1+y)^{\lambda+1} \pa_y f \|_{L^2}  
\end{equation}
\mspace
 


\section{Proof of technical lemmas}
\subsection{Proof of Lemma \ref{relations}}
The point is to prove that 
\begin{equation}
\| \pa_x^j \omega \|_{L^2_\gamma} \: \lesssim \: \| h_j \|_{L^2(\T \times \R_+)} + \| g_j \|_{L^2_\gamma} \: \lesssim \: \| \pa_x^j \omega \|_{L^2_\gamma}. 
\end{equation}
{\em Inequality on the right.} Clearly, $\| h_j \|_{L^2(\T \times \R_+)}  \lesssim \: \| \pa_x^j \omega \|_{L^2_\gamma}$. Then, by assumption  \eqref{lowerboundomega}, 
$$ (1+y)^\gamma |g_j| \: \lesssim (1+y)^\gamma |\pa_x^j \omega| \: + \: (1+y)^{\gamma-1}  |\pa_x^j u|
$$
and from Hardy inequality \eqref{hardy1}, we deduce: $ \| g_j \|_{L^2_\gamma} \: \lesssim \: \| \pa_x^j \omega \|_{L^2_\gamma}$. 

\mspace
{\em Inequality on the left.} First, from the definition of $h_j$, 
\begin{equation} \label{appendix1}
 \| \pa_x^j \omega \|_{L^2(\{ |y-a| < \eps\})} \: \lesssim \: \| h_j \|_{L^2(\T \times \R_+)} 
 \end{equation}
for $\eps > 0$ small enough so that $\{|y-a| < \eps \} \subset \{ \chi = 1 \}$. Then, we use the identities stated in Lemma \ref{decompo}: 
\begin{equation} \label{decompodxju}
 \pa^j_x u(t,x,y) = \left\{ \begin{aligned} & \omega(t,x,y)  \int_3^y \left(\psi + \frac{(1-\psi)}{\omega}\right)^{-1} \frac{g_j}{\omega^2} \: \: + \:C_j(t,x) \omega(t,x,y), \quad y > a(t,x),
  \\ &  \omega(t,x,y)  \int_0^y \left(\psi + \frac{(1-\psi)}{\omega}\right)^{-1} \frac{g_j}{\omega^2} , \quad y < a(t,x)  \end{aligned} \right. \end{equation} 
  We  integrate the first equality for $x \in \T$,  $\frac{\eps}{2} < y-a < \eps$.  As $\omega$ does not vanish in this region, it allows to obtain a control of $C_j$:   
\begin{align}
\nonumber
 \| C_j \|_{L^2(\T)} \: & \le \: C_\eps \left( \| \pa_x^j u \|_{L^2(\{\frac{\eps}{2} < y-a < \eps\})}    + \| g_j \|_{L^2(\{\frac{\eps}{2} < y-a < \eps\}\})} \right) \:   \\
\nonumber
& \le \: C'_\eps \left( \| \pa_x^j \omega \|_{L^2(\{0 <  y < \eps + a\})}    + \| g_j \|_{L^2(\{y  < 3 \})} \right)  \\
\nonumber
& \le \: C'_\eps  \left( \| \pa_x^j \omega \|_{L^2(\{ |y -a|  < \eps \})}    +   \| \pa_x^j \omega \|_{L^2(\{ y -a  < -\eps \})}  \: + \:  \|  g_j \|_{L^2(\{y  < 3 \})} \right) \\
  \label{estimCjlemma2} 
&   \le \: C''_\eps \left( \| h_j \|_{L^2(\T \times \R_+)} + \| g_j \|_{L^2(L^2_\gamma)}  \right). 
 \end{align}
  Finally, considering identities \eqref{decompodxju} for $|y-a| \ge \eps$, we get 
 $$ \| (1+y)^\gamma \pa^j_x \omega \|_{L^2(\{|y-a| \ge \eps\})} \: \le \: C_\eps \left(   \| C_j \|_{L^2(\T)} 
  + \| g_j \|_{L^2(L^2_\gamma)}  \right)  \le  C_\eps \left( \| h_j \|_{L^2(\T \times \R_+)} + \| g_j \|_{L^2(L^2_\gamma)}  \right)  $$
where the last inequality comes from \eqref{estimCjlemma2}. 
Combining with \eqref{appendix1} yields the result.

\subsection{Proof of Lemma \ref{lem-ujg}}
{\em Bound on $u^j_g$.} We denote $\psi_\omega := \left( \psi  + (1-\psi)/\omega\right)^{-1}$. For  $a - \eps \le y  \le 3$  
\begin{multline*}
\int_0^y u^j_g \: 
 = \: \int_0^a \omega \int_0^y  \frac{g \psi_\omega}{\omega^2} \: + \: \int_a^y \omega \int_3^y \frac{g \psi_\omega}{\omega^2}  \: = \: \int_0^a \frac{u\vert_{y=a}- u}{\omega^2}  \, g   \psi_\omega \:  + \:  \int_a^y  
\frac{u\vert_{y=a}-u}{\omega^2}  \, g   \psi_\omega \\
 + \: (u\vert_{y=a}-u(\cdot,y)) \int_3^y \frac{g \psi_\omega}{\omega^2} 
\end{multline*} 
after integration by parts. The third term at the right-hand side does not raise any problem, as the integral only involves values of $y$ $\eps$-away from $y=a$. For the third and second term, we notice that 
$$ \left| \frac{u - u\vert_{y=a}}{\omega^2}\right| \: \le \: C \quad \mbox{for } \: 0 \le y \le  3.   $$ 
Finally, we get 
$$ \left| \int_0^y u^j_g \right| \: \lesssim  \: \int_0^3 | g |, $$ 
and the bound follows easily. 

\mspace
{\em Bound on $C_j$.} From the proof of Lemma \ref{relations} ({\it cf} \eqref{estimCjlemma2}), we know that 
$$ \| C_j \|_{L^2(\T)} \: \lesssim \: \| h_j \|_{L^2(\T \times \R_+)} + \| g_j \|_{L^2(L^2_\gamma)}  \: \lesssim \: \| \pa^j_x \omega \|_{L^2(L^2_\gamma)}. $$

\mspace
{\em Bound on $\pa_x C_j$}:  
We write: for $y > a$, 
\begin{equation}
\omega \pa_x C_j = \pa_x^{j+1} u \: - \:  \pa_x \omega \, C_j   \: - \:    \pa_x \left( \omega \int_3^y \frac{g_j \psi_\omega}{\omega^2} \right)
\end{equation}
As for $C_j$, we integrate for $\frac{\eps}{2} \le y \le \eps$, $\eps > 0$ small enough. We end up with 
\begin{equation*}
\| \pa_x C_j \|_{L^2(\T)} \: \lesssim \: \|   C_j \|_{L^2(\T)} \: + \: \| \pa^{j+1}_x \omega \|_{L^2(L^2_\gamma)} + \| g_j \|_{L^2(\{ y < 3 \})} \: + \: \| \pa_x g_j \|_{L^2(\{ y < 3 \})}
\end{equation*}
Easily, the last two terms are bounded by $\| \pa^j_x \omega \|_{L^2(L^2_\gamma)} + \| \pa^{j+1}_x \omega \|_{L^2(L^2_\gamma)}$. The result follows. 

\subsection{Proof of Lemma \ref{lemmatildegjgj}}
It is enough to prove that  for all $\alpha \ge 0$, 
\begin{equation}
\| \, \| j^{\alpha/4} (\tilde g_j - g_j) \|_{L^2(\{ y \le 3\})}   \|_{l^2(\tau)}  \: \lesssim \:   \|  j^{(\alpha-3)/4} \| \omega \|_{{\cal H}^j_\gamma} \|_{l^2(\tau)} 
\end{equation}  
For $y \in [0,3]$, we have 
\begin{align*}
 \tilde g_j - g_j \: &  = \: [ \pa_x^{j-5} , \omega ] \pa_x^5 \omega  - [\pa_x^{j-5}, \pa_y \omega] \pa^5_x u \\
 & = \: \sum_{k=1}^{j-5} \binom{j-5}{k} \pa^k_x \om \, \pa^{j-k}_x \om - \pa^k_x \pa_y \om \, \pa^{j-k}_x u. 
 \end{align*}
 
This kind of sums has been considered several times throughout the paper. We get here
\begin{multline*}
\| \tilde g_j - g_j \|_{L^2(\{ y \le 3\})} \: \le \: \sum_{k=1}^{\frac{j-5}{2}} \binom{j-5}{k}  \sum_{m=0}^3 \| \om \|_{{\cal H}^{k+m}_\gamma} \, \| \omega \|_{{\cal H}^{j-k}_\gamma}   \\
+  \sum_{k=\frac{j-5}{2}}^{j-5}   \binom{j-5}{k} \left( \| \omega \|_{{\cal H}^{k}_\gamma} + \| \omega \|_{{\cal H}^{k+1}_\gamma} \right) \sum_{m=0}^2 \| \omega \|_{{\cal H}^{j-k+m}_\gamma} \: := \: A_j \: + \:  B_j.
\end{multline*}
We focus on $A_j$, as the other term $B_j$ is much better. Through the change of index $k' := k-1$, we find that 
\begin{align*}
A_j \: & \lesssim \: \sum_{k=0}^{\frac{j-5}{2}-1} \binom{j-5}{k+1} \sum_{m=1}^4 
\| \om \|_{{\cal H}^{k+m}_\gamma} \| \om \|_{{\cal H}^{j-1-k}_\gamma} \\
 & \lesssim \: \sum_{k=0}^{\frac{j-1}{2}} \binom{j-1}{k}   \sum_{m=1}^4 \| \om \|_{{\cal H}^{k+m}_\gamma}  \, 
(j-1- k) \| \om \|_{{\cal H}^{j-1-k}_\gamma}
\end{align*}
using the inequality $\binom{j}{k} \lesssim \binom{j-1}{k} j \:  \lesssim \binom{j-1}{k} (j-1-k)$. It follows that 
$$ \| j^{\gamma/4} A_j \|_{l^2(\tau)} \:  \lesssim \:   \| j^{\gamma/4} j  \| \omega \|_{{\cal H}^{j-1}_\gamma} \|_{l^2(\tau)} \: \le \:  \| j^{(\gamma-3)/4}   \| \omega \|_{{\cal H}^j_\gamma} \|_{l^2(\tau)}. $$
This leads to the expected bound. 

\subsection{Proof of Lemma \ref{lemmatildegj}}
We only treat the case $l=1$: 
$$ \|  \| j^\alpha \pa_y \bar{g}_{j-1} \|_{L^2(\T \times \R_+)} \|_{l^2(\tau)} \: \lesssim \:  \| j^\alpha \| \omega \|_{{\cal H}^j_\gamma} \|_{l^2(\tau)}.$$
The other cases $l=0$ and $l=2$ are handled in the same way. We compute
\begin{align*}
\pa_y \tilde g_{j-1} & = \pa_x^{j-6} \left( \omega \, \pa_y \pa_x^5 \omega - \pa^2_y \omega \, \pa_x^5 u \right) \\
& = \sum_{k=0}^{j-6} \binom{j-6}{k} \pa_x^k \omega \, \pa_y \pa_x^{j-k-1} \omega \: - \: \sum_{k=0}^{j-6}  \binom{j-6}{k}  \pa^k_x \pa^2_y \omega \, \pa_x^{j-k-1} u 
\end{align*}
By now standard computations:
\begin{multline*} 
 \| \pa_y \tilde g_{j-1} \|_{L^2(\T \times \R_+)} \: \lesssim \: \sum_{k=0}^{\frac{j-6}{2}} \binom{j-6}{k} \sum_{m=0}^4 \| \omega \|_{{\cal H}_\gamma^{k+m}} \, \left(  \| \omega \|_{{\cal H}^{j-k}_\gamma} +  \| \omega \|_{{\cal H}^{j-1-k}_\gamma} \right) \\
 + \sum_{\frac{j-6}{2}}^{j-6}  \binom{j-6}{k}  \left(  \| \omega \|_{{\cal H}^{k}_\gamma} +  \| \omega \|_{{\cal H}^{k+2}_\gamma} \right) \sum_{m=-1}^2
\| \omega \|_{{\cal H}^{j-k+m}}  \: := \: A_j + B_j.
\end{multline*}
On one hand, as $\binom{j-6}{k} \lesssim \binom{j}{k}$, Lemma \ref{binom} clearly implies 
$$ \| j^\alpha A_j \|_{l^2(\tau)} \: \lesssim \: \| j^\alpha \| \omega \|_{{\cal H}^{j}_\gamma} \|_{l^2(\tau)}. $$
On the other hand, after the change of index $k' := k+2$ 
$$ B_j \: \lesssim \: \sum_{k=\frac{j}{2}}^j \binom{j}{k} \left(  \| \omega \|_{{\cal H}^{k-2}_\gamma} +  \| \omega \|_{{\cal H}^k_\gamma} \right) \sum_{m=1}^4 \| \omega \|_{{\cal H}^{j-k+m}} $$
(notice that $\binom{j-6}{k-2} \lesssim \binom{j}{k}$). Still by Lemma \ref{binom},
$$ \| j^\alpha B_j \|_{l^2(\tau)} \: \lesssim \: \| j^\alpha \| \omega \|_{{\cal H}^{j}_\gamma} \|_{l^2(\tau)}. $$ 
Combining the bounds on $A_j$ and $B_j$ yields the result. 

\subsection{Proof of Lemma \ref{dytildegjgj}}
For $(x,y)$ in the support of $\chi^{\eps,a}$, we have 
$$ \pa_y \left( \tilde g_j - g_j \right) = \sum_{k=1}^{j-5} \binom{j-5}{k} \left( \pa_x^k \omega \pa_x^{j-k} \, \pa_y \omega -  \pa_x^k \pa^2_y \omega \, \pa_x^{j-k} u \right).  $$
We proceed as we did several times in the paper: we get 
\begin{align*}
& \| \sqrt{\chi^{\eps,a}}   \pa_y \left( \tilde g_j - g_j \right) \|_{L^2(\T \times \R_+)} \\ 
&  \lesssim \: \sum_{k=1}^{\frac{j-5}{2}}  \sum_{m=0}^4 \| \omega \|_{{\cal H}^{k+m}_\gamma}  \left( \|  \sqrt{\chi^{\eps,a}} \pa^{j-k}_x \pa_y \omega \|_{L^2(\T \times \R_+)} \: + \: \| \omega \|_{{\cal H}^{j-k}_\gamma} \right) \\
& + \: \sum_{k=\frac{j-5}{2}}^{j-5} \binom{j-5}{k} \left( \| \omega \|_{{\cal H}^k_\gamma}  \sum_{m=1}^3 \| \omega \|_{{\cal H}^{j-k+m}_\gamma} \: + \:   \| \omega \|_{{\cal H}^{k+2}_\gamma}   \sum_{m=0}^2 \| \omega \|_{{\cal H}^{j-k+m}_\gamma}  \right)  \\ 
&  := \: A_j + B_j. 
\end{align*}
Note that we did not bound  $\|  \sqrt{\chi^{\eps,a}} \pa^{j-k}_x \pa_y \omega \|_{L^2}$ by $\| \pa_y \omega \|_{{\cal H}^{j-k}_\gamma}$:  we need to keep track of the localization to handle this term. 

\mspace
{\em Treatment of $A_j$.} We set $k':= k-1$. Dropping the prime, 
\begin{align*}
A_j \: & \lesssim \: \sum_{k=0}^{\frac{j-5}{2}-1} \binom{j-5}{k+1}   \sum_{m=1}^5 \| \omega \|_{{\cal H}^{k+m}_\gamma}  \left( \|  \sqrt{\chi^{\eps,a}} \pa^{j-1-k}_x \pa_y \omega \|_{L^2(\T \times \R_+)} \: + \: \| \omega \|_{{\cal H}^{j-1-k}_\gamma} \right) \\
 & \lesssim \: \sum_{k=0}^{\frac{j-1}{2}} \binom{j-1}{k} \sum_{m=0}^5 \| \omega \|_{{\cal H}^{k+m}_\gamma} (j-1-k) \left( \|  \sqrt{\chi^{\eps,a}} \pa^{j-1-k}_x \pa_y \omega \|_{L^2(\T \times \R_+)} \: + \: \| \omega \|_{{\cal H}^{j-1-k}_\gamma} \right)
 \end{align*}
 Here, we have used $\binom{j-5}{k+1} \: \lesssim \:  j \, \binom{j-1}{k} \: \lesssim  \: (j-1-k) \,  \binom{j-1}{k}$. By lemma \ref{binom}, we find 
 \begin{align*}
 \| j^{3/4} A_j \|_{l^2(\tau)} \: & \lesssim \: \| \| \omega \|_{{\cal H}^j_\gamma} \|_{l^2(\tau)} \left( \| j^{7/4} \|  \sqrt{\chi^{\eps,a}} \pa^{j-1}_x \pa_y \omega \|_{L^2(\T \times \R_+)}  \|_{l^2(\tau)} +  \| j^{7/4} \| \omega \|_{{\cal H}^{j-1}_\gamma} \|_{l^2(\tau)} \right)  \\
 & \lesssim \:  \| \| \omega \|_{{\cal H}^j_\gamma} \|_{l^2(\tau)}  \left( \| \|  \sqrt{\chi^{\eps,a}} \pa^j_x \pa_y \omega \|_{L^2(\T \times \R_+)}  \|_{l^2(\tau)} +  \| \| \omega \|_{{\cal H}^{j}_\gamma} \|_{l^2(\tau)} \right). 
 \end{align*}
 Finally, we notice that for $\eps$ small enough:
 \begin{equation*} 
   \| \sqrt{\chi^{\eps,a}} \pa^j_x \pa_y\omega \|_{L^2(\T \times \R_+)} \: \lesssim \: \| \pa_y h_j \|_{L^2(\T \times \R_+)} \: + \: \| \pa^j_x \omega \|_{L^2(\T \times \R_+)} 
 \end{equation*}
 so that 
 $$ \| j^{3/4} A_j \|_{l^2(\tau)} \: \lesssim \: \sqrt{E_\omega(t,\tau)} \left( \sqrt{D_h(t,\tau)} \: + \: \sqrt{E_\omega(t,\tau)} \right) \: \lesssim  \:  \sqrt{D_h(t,\tau)} \: + \: \sqrt{E_\omega(t,\tau)}$$
 
 \mspace
 {\em Treatment of $B_j$.} Through the change  of index $k' := k+1$ (first term) and $k' := k+3$ (second term), and using that $\binom{j-5}{k-d} \lesssim \binom{j}{k}$, $d=1$ or $3$, we get:
 \begin{equation*}
 B_j \: \lesssim \: \sum_{k=\frac{j}{2}}^j \binom{j}{k}  \| \omega \|_{{\cal H}^{k-2}_\gamma}  \sum_{m=2}^5 \| \omega \|_{{\cal H}^{j-k+m}_\gamma}.
 \end{equation*}
Thus, 
\begin{align*}
\| j^{3/4} B_j \|_{l^2(\tau)} \: & \lesssim \: \| \| \omega   \|_{{\cal H}^j_\gamma} \|_{l^2(\tau)} \:  \| j^{3/4} \|  \omega \|_{{\cal H}^{j-1}_\gamma} \|_{l^2(\tau)} \\
  & \lesssim \: \| \| \omega   \|_{{\cal H}^j_\gamma} \|_{l^2(\tau)} \: \| j^{-1} \| \omega \|_{{\cal H}^j_\gamma} \|_{l^2(\tau)}  \: \lesssim \sqrt{E_\omega(t,\tau)}. 
  \end{align*}
Combining the bounds on $A_j$ and $B_j$ yields the result. 

\bibliographystyle{abbrv}
\bibliography{biblio}

\end{document}